\documentclass{article}
\usepackage[utf8]{inputenc}

\usepackage{euscript}

\usepackage{amsmath,amsthm,amssymb}

\usepackage{enumerate}

\usepackage{xcolor}

\usepackage[colorlinks=true, pdfstartview=FitV, linkcolor=blue, citecolor=blue, urlcolor=blue]{hyperref}

\usepackage[english]{babel}
\usepackage{amsfonts}

\usepackage{comment}

\definecolor{darkred}{rgb}{0.7,0,0} 
\newcommand{\defn}[1]{{\color{darkred}\emph{#1}}} 

\usepackage{todonotes}

\newcommand{\PER}{\operatorname{PER}}
\newcommand{\PD}{\operatorname{PD}}
\newcommand{\proj}{\operatorname{proj}}

\newcommand{\RR}{\mathbb{R}}

\newcommand{\NN}{\mathbb{N}}
\newcommand{\ZZ}{\mathbb{Z}}

\newcommand{\PPos}{\RR_{>0}}

\newcommand{\set}[1]{\left\{ #1 \right\}}
\newcommand{\abs}[1]{\left| #1 \right|}
\newcommand{\tup}[1]{\left( #1 \right)}
\newcommand{\ive}[1]{\left[ #1 \right]}
\newcommand{\floor}[1]{\left\lfloor #1 \right\rfloor}

\newcommand{\calF}{\mathcal{F}}
\newcommand{\calB}{\mathcal{B}}

\newcommand{\IN}{\ive{N}}
\newcommand{\Eh}{\widehat{E}}
\newcommand{\wh}{\widehat{w}}
\renewcommand{\dh}{\widehat{d}}
\newcommand{\Ch}{\widehat{C}}

\newcommand{\timesu}{\mathbin{\underline{\times}}}

\newcommand{\verts}[1]{\operatorname{V}\left( #1 \right)}
\newcommand{\edges}[1]{\operatorname{E}\left( #1 \right)}

\renewcommand{\leq}{\leqslant}
\renewcommand{\geq}{\geqslant}

\renewcommand{\subset}{\subseteq}

\newcommand{\nuo}{\nu^{\circ}}

\newtheoremstyle{plainsl}
  {8pt plus 2pt minus 4pt}
  {8pt plus 2pt minus 4pt}
  {\slshape}
  {0pt}
  {\bfseries}
  {.}
  {5pt plus 1pt minus 1pt}
  {}

\theoremstyle{plainsl}
\newtheorem{theorem}{Theorem}[section]
\newtheorem{lemma}[theorem]{Lemma}
\newtheorem{corollary}[theorem]{Corollary}
\newtheorem{proposition}[theorem]{Proposition}
\newtheorem{example}[theorem]{Example}
\newtheorem{remark}[theorem]{Remark}
\newtheorem{definition}[theorem]{Definition}
\newtheorem{question}[theorem]{Question}

\newenvironment{statement}{\begin{quote}}{\end{quote}}

\newenvironment{verlong}{}{}

\excludecomment{verlong}
\includecomment{vershort}
\excludecomment{noncompile}

\title{A greedoid and a matroid inspired by Bhargava's $p$-orderings}
\author{Darij Grinberg\thanks{DG thanks the Mathematisches Forschungsinstitut Oberwolfach for its hospitality.} \\
\small Mathematics Department\\[-0.8ex]
\small Drexel University\\[-0.8ex]
\small Philadelphia, PA 19104, U.S.A.\\[-0.8ex]
\small \url{darijgrinberg@gmail.com} / \url{http://www.cip.ifi.lmu.de/~grinberg/}
\and
Fedor Petrov\\
\small St. Petersburg State University and\\[-0.8ex]
\small St. Petersburg Department\\[-0.8ex]
\small Steklov Mathematical Institute of Russian Academy of Sciences\\[-0.8ex]
\small St. Petersburg, Russia.\\[-0.8ex]
\small \url{fedyapetrov@gmail.com} / \url{http://math-cs.spbu.ru/people/petrov-f-v/}
}
\date{29 March 2021}

\begin{document}

\maketitle

\begin{abstract}
Consider a finite set $E$.
Assume that each $e \in E$ has a
``weight'' $w \left(e\right) \in \mathbb{R}$ assigned to it,
and any two distinct $e, f \in E$ have a ``distance''
$d \left(e, f\right) = d \left(f, e\right) \in \mathbb{R}$
assigned to them,
such that the distances satisfy the ultrametric
triangle inequality
$d(a,b)\leqslant \max \left\{d(a,c),d(b,c)\right\}$.
We look for a subset of $E$ of given size with maximum perimeter
(where the perimeter is defined by summing the weights of
all elements and their pairwise distances).
We show that any such subset can be found by a greedy algorithm
(which starts with the empty set, and then adds new elements one
by one, maximizing the perimeter at each step).
We use this to define numerical invariants, and also to
show that the maximum-perimeter subsets of all
sizes are the feasible sets of a strong greedoid,
and the maximum-perimeter subsets
of any given size are the bases of a matroid.
This essentially generalizes the ``$P$-orderings'' constructed by
Bhargava in order to define his generalized factorials, and is also
similar to the strong greedoid of maximum diversity subsets in
phylogenetic trees studied by Moulton, Semple and Steel.

We further discuss some numerical invariants of $E, w, d$ stemming
from this construction, as well as an analogue where
maximum-perimeter subsets are replaced by maximum-perimeter
tuples (i.e., elements can appear multiple times).
\end{abstract}

\tableofcontents

\section{Introduction}

In this paper, we study a combinatorial setting
consisting of a finite set $E$ with a
``weight function'' $w : E \to \mathbb{R}$
and a
(symmetric) ``distance function'' $d : E \timesu E \to \mathbb{R}$
(where $E \timesu E = \set{\tup{e, f} \in E \times E \mid e \neq f}$)
satisfying the ultrametric triangle inequality.
This generalizes the notion of an ultrametric space.
Given any finite subset $A$ of $E$, we can define
the perimeter of $A$ to be the sum of the weights and
of the pairwise distances of the elements of $A$.
Given an integer $k \geq 0$ and a finite subset $C$ of $E$,
we show (Theorem~\ref{thm.2a}) that if we want
to construct a $k$-element subset of $C$ having maximum
perimeter, we can do so by a greedy algorithm (i.e., by starting
with the empty set and repeatedly adding new elements that increase
the perimeter as much as possible), and that every
maximum-perimeter $k$-element subset of $C$ can be
constructed through this algorithm (Theorem~\ref{thm.2b}).
We furthermore show that these maximum-perimeter $k$-element
subsets are the bases of a matroid (when $k$ is fixed) and
the feasible sets of a strong greedoid
(when $k$ ranges over all nonnegative integers).
In a followup paper \cite{Grinbe20}, this strong greedoid
is studied from an algebraic viewpoint, which also addresses
questions of its (linear) representability.

Our greedy construction of maximum-perimeter subsets
is inspired by Manjul Bhargava's concept of a $P$-ordering
(\cite[Section 2]{Bharga97}), which laid the foundation for
his theory of generalized factorials
(see \cite[Section 4]{Bharga00} and \cite[Section 2]{Bharga09});
we connect the two notions (in Section~\ref{sect.bhargava})
and obtain new proofs of two results from \cite[Section 2]{Bharga97}.

A similar problem -- also leading to a strong greedoid --
has appeared in the mathematical biology literature:
Given a phylogenetic tree $T$ and an integer $k$, the problem
asks to find a set of $k$ leaves of $T$ having maximum
``phylogenetic diversity'' (i.e., the total weight of the edges
of the subtree that connects these $k$ leaves).
In \cite{MoSeSt06}, Moulton, Semple and Steel show that such
diversity-maximizing $k$-element sets form a strong greedoid,
just as our maximum-perimeter subsets do.
The similarity does not end here:
Phylogenetic trees are close relatives of ultra triples
(and can be translated to and from the latter without much loss
of information).
However, the strong greedoid of Moulton, Semple and Steel
is not the same as ours, since perimeter (when
restated in terms of the phylogenetic tree) is not the same as
phylogenetic diversity\footnote{Roughly speaking, in a
star-shaped phylogenetic tree with $1$ internal vertex and $p$
leaves, the perimeter of a $k$-leaf set
is quadratic in $k$, while its phylogenetic diversity is
linear in $k$.
Also, our Lemma~\ref{lem.strong-greedoid.AB}, while being
an analogue of \cite[Lemma 3.1]{MoSeSt06}, differs from
the latter in that it requires
$\abs{B} = \abs{A} + 1$ rather than $\abs{B} > \abs{A}$
(and indeed, the latter requirement would not suffice).}.
It is an interesting question to what extent these two
problems can be reconciled, and perhaps a more general
class of optimization problems on phylogenetic trees (or
ultra triples) can be shown to lead to a strong greedoid.


\section{The setup}

\subsection{Defining ultra triples}

Let $E$ be a set.
We shall use $E$ as our ground set throughout this paper.

We shall refer to the elements of $E$ as \defn{points}.

For a nonnegative integer $m$, an \defn{$m$-set} means
a subset $A$ of $E$ which consists of $\abs{A}=m$ elements.
If $B \subset E$ is any subset
and $m$ is a nonnegative integer,
then an \defn{$m$-subset of $B$} means an $m$-element
subset of $B$.

Define the set \defn{$E \timesu E$} by
\[
E \timesu E = \set{\tup{e, f} \in E \times E \mid e \neq f} .
\]
Thus, $E \timesu E$ is the set of all ordered pairs
$\tup{e, f}$ of two distinct elements of $E$.

Assume that we are given a function $w : E \to \RR$.
In other words, each point $a\in E$ has a
real-valued \defn{weight} $w\tup{a}$ assigned to it.

Assume further that we are given a function
$d : E \timesu E \to \RR$, which we will call
the \defn{distance function}.
Thus, any two distinct points $a,b\in E$
have a real-valued \defn{distance} $d\tup{a,b}$.
We assume that this distance function has the following properties:

\begin{itemize}

\item It is \defn{symmetric}: that is, $d\tup{a,b} = d\tup{b,a}$
      for any two distinct $a, b \in E$.

\item It satisfies the following inequality:
        \begin{equation}\label{ultra}
        d(a,b)\leqslant \max \set{d(a,c),d(b,c)}
        \end{equation}
      for any three distinct $a, b, c \in E$.
\end{itemize}

\noindent
(The inequality \eqref{ultra} is commonly known as the
\defn{ultrametric triangle inequality}\footnote{It
can be restated as ``the longest two sides
of a triangle are always equal in length''.
(Here, a \defn{triangle} means a $3$-subset of $E$;
its \defn{sides} are its $2$-subsets; the
\defn{length} of a side
$\set{u, v}$ is $d\tup{u, v}$.)}; but unlike
the distance function of an ultrametric space, our
$d$ can take negative values.
The values of $w$ are completely unrestrained.)

Such a structure $\tup{E,w,d}$ will be called an
\defn{ultra triple}.

From now on, we shall always be considering an
ultra triple $\tup{E,w,d}$ (unless stated otherwise).

\subsection{Examples}

We shall now provide a few examples of ultra triples.
In each case, the proof that our triple is indeed an
ultra triple is easy and left to the reader.

\begin{example} \label{exa.01.d}
For this example, we let $E$ be an arbitrary set,
and we define the distances $d\tup{a, b}$ as follows:
\[
d\tup{a, b} = 1
\qquad \text{ for all } \tup{a, b} \in E \timesu E .
\]
We define the weights $w\tup{a}$ arbitrarily.
Then, $(E, w, d)$ is an ultra triple.
\end{example}

\begin{example} \label{exa.5-point-1.d}
For this example, we let $E = \set{1,2,3,4,5}$,
and we define the distances $d\tup{a, b}$ as follows:
\[
d\tup{a, b}
= \begin{cases}
      1, & \text{if } a \equiv b \mod 2 ; \\
      2, & \text{if } a \not\equiv b \mod 2
  \end{cases}
\qquad \text{ for all } \tup{a, b} \in E \timesu E .
\]
We define the weights $w\tup{a}$ arbitrarily.
Then, $(E, w, d)$ is an ultra triple.
\end{example}

\begin{example} \label{exa.mod-or-not-mod.d}
For this example, we fix two reals $\varepsilon$
and $\alpha$ with $\varepsilon \leq \alpha$.
Furthermore, we fix an integer $m$
and a subset $E$ of $\ZZ$.
We define the distance function
$d : E \timesu E \to \RR$ by setting
\[
d\tup{a, b}
= \begin{cases}
      \varepsilon, & \text{if } a \equiv b \mod m ; \\
      \alpha, & \text{if } a \not\equiv b \mod m
  \end{cases}
\qquad \text{ for all } \tup{a, b} \in E \timesu E .
\]
We define the weights $w\tup{a}$ arbitrarily.
Then, $(E, w, d)$ is an ultra triple.
\end{example}

Note that Example~\ref{exa.5-point-1.d} is the
particular case of Example~\ref{exa.mod-or-not-mod.d}
obtained by setting $\varepsilon = 1$,
$\alpha = 2$, $m = 2$ and $E = \set{1,2,3,4,5}$.

\begin{example} \label{exa.p-adic.d}
For this example, we fix a prime number $p$
and a subset $E$ of $\ZZ$,
and we define the distances $d\tup{a, b}$ as follows:
\[
d\tup{a, b}
= p^{-v_p\tup{a-b}}
\qquad \text{ for all } \tup{a, b} \in E \timesu E .
\]
Here, for any nonzero integer $m$, we let $v_p\tup{m}$
denote the \defn{$p$-adic valuation} of $m$
(that is, the largest nonnegative integer $k$ such that $p^k \mid m$).
The distance function $d : E \timesu E \to \RR$ is
called the \defn{$p$-adic metric}.
We define the weights $w\tup{a}$ arbitrarily.
Then, $(E, w, d)$ is an ultra triple.
\end{example}

\begin{example} \label{exa.p-adic.d'}
For this example, we fix a prime number $p$
and a subset $E$ of $\ZZ$.

For any nonzero integer $m$, we define $v_p\tup{m}$
as in Example~\ref{exa.p-adic.d}.
We define a map $d' : E \timesu E \to \RR$ by setting
\[
d' \tup{a, b}
= -v_p\tup{a-b}
\qquad \text{ for all } \tup{a, b} \in E \timesu E .
\]
We define the weights $w\tup{a}$ arbitrarily.
Then, $(E, w, d')$ is an ultra triple.
\end{example}

%

Most of the examples above are particular cases of
a more general construction:

\begin{example} \label{exa.n-adic}

Let $\NN = \set{0, 1, 2, \ldots}$.
Let $c : \NN \to \RR$ be a weakly decreasing function.

Fix a sequence $\mathbf{r} = \tup{r_0, r_1, r_2, \ldots}$
of integers such that $r_0 \mid r_1 \mid r_2 \mid \cdots$.

For each $x \in \ZZ$,
define an element $v_{\mathbf{r}}\tup{x} \in \NN$
by
\[
v_{\mathbf{r}}\tup{x}
= \max\set{ i \in \NN \text{ such that } r_i \mid x } ,
\]
assuming that this maximum exists.
(Otherwise, leave $v_{\mathbf{r}}\tup{x}$ undefined.)

Let $E$ be a subset of $\ZZ$.
Define a distance function
$d : E \timesu E \to \RR$ by setting
\[
d\tup{a, b}
= c\tup{v_{\mathbf{r}}\tup{a-b}}
\qquad \text{ for all } \tup{a, b} \in E \timesu E .
\]
Assume that this is well-defined
(i.e., the values $v_{\mathbf{r}}\tup{a-b}$ are
well-defined for all $\tup{a, b} \in E \timesu E$).
We define the weights $w\tup{a}$ arbitrarily.
Then, $(E, w, d)$ is an ultra triple.
\end{example}

Example~\ref{exa.mod-or-not-mod.d}
is obtained from Example~\ref{exa.n-adic}
by setting $r_0 = 1$ and
$r_1 = m$ and $r_2 = r_3 = r_4 = \cdots = 0$
and $c\tup{0} = \alpha$ and
$c\tup{1} = \varepsilon$
(defining the remaining values of $c$
arbitrarily to be weakly decreasing).
Example~\ref{exa.p-adic.d}
is obtained from Example~\ref{exa.n-adic}
by setting $r_i = p^i$ and
$c\tup{n} = p^{-n}$
(indeed, if we set $r_i = p^i$,
then $v_{\mathbf{r}}\tup{m} = v_p\tup{m}$
for each nonzero $m \in \ZZ$).
Likewise, Example~\ref{exa.p-adic.d'}
is obtained from Example~\ref{exa.n-adic}
by setting $r_i = p^i$ and
$c\tup{n} = -n$.

An even more general (and simpler) example of
an ultra triple (more precisely, of a distance
function satisfying \eqref{ultra})
can be obtained from a hierarchy of equivalence
relations:

\begin{example} \label{exa.eqrels}
Let $E$ be a set.
Let $\underset{0}{\sim}, \underset{1}{\sim}, \underset{2}{\sim}, \ldots$
be equivalence relations on $E$.
Assume that:\footnote{These three assumptions can be
restated in terms of the \defn{partition lattice} on
$E$, which is the lattice of all set partitions of $E$
(see, e.g., \cite[Section 1.7]{Oxley11}).
Indeed, it is well-known that the equivalence relations
on $E$ are in bijection with the set partitions of $E$.
Our sequence $\underset{0}{\sim}, \underset{1}{\sim}, \underset{2}{\sim}, \ldots$
of equivalence relations thus corresponds to a
sequence $P_0, P_1, P_2, \ldots$ of set partitions of $E$.
The three assumptions below thus say that
$P_0$ is the trivial partition; that we have
$P_0 \geq P_1 \geq P_2 \geq \cdots$ (meaning that each
partition $P_i$ refines $P_{i-1}$);
and that the meet $\bigwedge_{i=0}^\infty P_i$ is
the partition of $E$ into singletons.}

\begin{itemize}

\item[(A)] Every $e, f \in E$ satisfy $e \underset{0}{\sim} f$.

\item[(B)] If some $e, f \in E$ and $i > 0$ satisfy $e \underset{i}{\sim} f$,
then $e \underset{i-1}{\sim} f$.

\item[(C)] If $e, f \in E$ are distinct,
then there exists some $i \geq 0$
such that we don't have $e \underset{i}{\sim} f$.

\end{itemize}

Let $c : \NN \to \RR$ be a weakly decreasing function.

Define a distance function $d : E \timesu E \to \RR$
by
\[
d\tup{e, f}
= c\tup{\max\set{ i \geq 0 \mid e \underset{i}{\sim} f }}
\qquad
\text{ for all } \tup{e, f} \in E \timesu E .
\]

We define the weights $w\tup{a}$ arbitrarily.
Then, $(E, w, d)$ is an ultra triple.
\end{example}

We can obtain Example~\ref{exa.n-adic} from
Example~\ref{exa.eqrels} by defining the relation
$\underset{i}{\sim}$ to be congruence modulo $r_i$.
(The assumption that the $v_{\mathbf{r}}\tup{a-b}$ are
well-defined in Example~\ref{exa.n-adic} ensures
that assumptions (A), (B) and (C) of
Example~\ref{exa.eqrels} are satisfied.)

Hierarchical taxonomies can be viewed as sets $E$
equipped with sequences
$\underset{0}{\sim}, \underset{1}{\sim}, \underset{2}{\sim}, \ldots$
of equivalence relations (usually finite, however)
satisfying assumptions (A), (B) and (C) of
Example~\ref{exa.eqrels}.
For example:

\begin{example} \label{exa.taxon}
Let $E$ be the set of all living organisms.
Define equivalence relations
$\underset{0}{\sim}, \underset{1}{\sim}, \underset{2}{\sim}, \ldots$
on $E$ as follows:
\begin{align*}
\tup{e \underset{0}{\sim} f} & \text{ always holds;} \\
\tup{e \underset{1}{\sim} f} & \iff \tup{e \text{ and } f \text{ belong to the same domain}}; \\
\tup{e \underset{2}{\sim} f} & \iff \tup{e \text{ and } f \text{ belong to the same kingdom}}; \\
\vdots \\
\tup{e \underset{7}{\sim} f} & \iff \tup{e \text{ and } f \text{ belong to the same genus}}; \\
\tup{e \underset{8}{\sim} f} & \iff \tup{e \text{ and } f \text{ belong to the same species}}; \\
\tup{e \underset{i}{\sim} f} & \iff \tup{e = f} \qquad \text{ for all } i \geq 9
\end{align*}
(following the taxonomic ranks of biology).
Then, assumptions (A), (B) and (C) of
Example~\ref{exa.eqrels} are satisfied.
\end{example}

Example~\ref{exa.taxon} yields not so much a genuine
biological application of our theory as it does a
helpful mental model for it.
A less naive (if still simplified) model of the
interrelation of organisms is provided by phylogenetic
trees \cite{MoSeSt06}
-- rooted trees (in the combinatorial sense)
whose vertices correspond to organisms or species,
and whose edges signify relationships of ancestry.
Often the edges are equipped with weights (or, better,
lengths) encoding the evolutionary distance between
parent and child vertices.
This model, too, leads to an ultra triple.
At the mathematical heart of this construction
is the following example:

\begin{example} \label{exa.phylo}
A \defn{tree} is a connected finite undirected
graph that has no cycles.
(Thus, our trees are unrooted and have no order-like
structures assigned to them.)

Let $T$ be a tree.
For each edge $e$ of $T$, let $\lambda\tup{e}$ be a
nonnegative real.
We shall call this real the \defn{weight} of $e$.

For any vertices $u$ and $v$ of $T$, let
\defn{$\lambda\tup{u, v}$} denote the sum of the
weights of all edges on the (unique) path
from $u$ to $v$.
Note that this $\lambda\tup{u, v}$ generalizes
the usual (graph-theoretical) distance between
$u$ and $v$; indeed,
if $\lambda\tup{e} = 1$ for each edge $e$ of $T$,
then $\lambda\tup{u, v}$ is the distance between
$u$ and $v$ (that is, the length of the unique
path from $u$ to $v$).

A known result (the ``four-point condition'') says
that if $x, y, z, w$ are four vertices of $T$,
then the two largest of the three numbers
\[
\lambda\tup{x, y} + \lambda\tup{z, w}, \qquad
\lambda\tup{x, z} + \lambda\tup{y, w}, \qquad \text{ and }
\lambda\tup{x, w} + \lambda\tup{y, z}
\]
are equal.
In the particular case when each edge of $T$
has weight $1$, this
is a standard exercise in graph theory
(see, e.g., \url{https://math.stackexchange.com/questions/2899278});
most of its solutions generalize to the
case of arbitrary weights.
We can use this to define an ultra triple as follows:

Fix any vertex $r$ of $T$.
Let $E$ be any subset of the vertex set of $T$.
We define a map $w : E \to \RR$ by setting
\[
w\tup{x} = \lambda\tup{x, r}
\qquad \text{ for each } x \in E.
\]
We define a map $d : E \timesu E \to \RR$ by
setting
\[
d\tup{x, y} = \lambda\tup{x, y} - \lambda\tup{x, r} - \lambda\tup{y, r}
\qquad \text{ for each } \tup{x, y} \in E \timesu E .
\]
Then, $\tup{E, w, d}$ is an ultra triple.

Note that our special choice of $w$ was not
necessary for this (any function $w : E \to \RR$
would have worked),
but it has the advantage that
$\lambda\tup{x, y} = w\tup{x} + w\tup{y} + d\tup{x, y}$
for any distinct $x, y \in E$.
The right hand side of this equality will later
be called the perimeter $\PER\set{x, y}$ of the
set $\set{x, y}$.
\end{example}

\subsection{Projections}

Let us now return to the setting of an arbitrary ultra
triple $\tup{E, w, d}$.

\begin{definition}
Let $C \subset E$ be a non-empty subset.
Let $v \in E$ be any point.

We define a subset \defn{$\proj_C\tup{v}$} of $C$ as follows:
\begin{itemize}
\item If $v \in C$, then we define $\proj_C\tup{v}$
      to be the one-element set $\set{v}$.
\item If $v \notin C$, then we define $\proj_C\tup{v}$
      to be the set of all $c \in C$ that minimize
      the distance $d\tup{v, c}$\ \ \ \ %
      \footnote{This distance $d \tup{v, c}$ is well-defined,
      since $v \neq c$ (because $v \notin C$ and $c \in C$).}.
\end{itemize}

The elements of $\proj_C\tup{v}$ will be called the
\defn{projections of $v$ onto $C$}.
\end{definition}

The following is easy to see and will be used
without explicit mention:

\begin{proposition} \label{prop.proj-ex}
Let $C \subset E$ be a finite non-empty subset.
Let $v \in E$ be any point.
Then, there exists at least one projection of $v$
onto $C$.
\end{proposition}

\begin{proof}
We are in one of the following two cases:

\textit{Case 1:} We have $v \in C$.

\textit{Case 2:} We have $v \notin C$.

Let us consider Case 1.
In this case, we have $v \in C$.
Hence, the definition of $\proj_C\tup{v}$
yields $\proj_C\tup{v} = \set{v}$ and thus
$v \in \set{v} = \proj_C\tup{v}$.
Thus, $v$ itself is a projection of $v$ onto $C$.
Hence, there exists at least one projection of $v$
onto $C$.
This proves Proposition~\ref{prop.proj-ex} in Case 1.

Let us now consider Case 2.
In this case, we have $v \notin C$.
Hence, the elements of $\proj_C\tup{v}$ are
the $c \in C$ that minimize
the distance $d\tup{v, c}$
(by the definition of $\proj_C\tup{v}$).
Clearly, there exists at least one such $c$
(since $C$ is finite and non-empty).
Thus, there exists at least one element of $\proj_C\tup{v}$.
In other words, there exists at least one projection
of $v$ onto $C$.
This proves Proposition~\ref{prop.proj-ex} in Case 2.

Hence, Proposition~\ref{prop.proj-ex} is proven in both Cases 1 and 2.
\end{proof}

\begin{example} \label{exa.5-point-1.proj}
Let $(E, w, d)$ be as in Example~\ref{exa.5-point-1.d}.
Then, the projections of $2$ onto $\set{1, 3}$ are
$1$ and $3$,
while the only projection of $2$ onto $\set{1, 3, 4}$
is $4$.
\end{example}

In Example~\ref{exa.phylo}, a projection of a
$v \notin C$ onto a subset $C$ is usually called
a ``closest relative of $v$ in $C$''.

The crucial property of projections is the following:

\begin{lemma}\label{lem.projection}
Assume that $C\subset E$ is a non-empty subset
and $v\in E$ is any point.
Let $u$ be a projection of $v$ onto $C$.

\begin{enumerate}

\item[\textbf{(a)}] If $v \in C$, then $u = v$.

\item[\textbf{(b)}] If $x \in C$ satisfies $x \neq u$,
then $x \neq v$.

\item[\textbf{(c)}]
Let $x \in C$ be such that $x \neq u$.
Then, $d(u,x)\leqslant d(v,x)$.

\end{enumerate}
\end{lemma}

\begin{proof}
We have $u \in \proj_C\tup{v}$ (since $u$ is
a projection of $v$ onto $C$).

\textbf{(a)} Assume that $v \in C$.
Hence, the definition of $\proj_C\tup{v}$
yields $\proj_C\tup{v} = \set{v}$.
Therefore,
$u \in \proj_C\tup{v} = \set{v}$.
In other words, $u = v$.
This proves Lemma~\ref{lem.projection} \textbf{(a)}.

\textbf{(b)} Let $x \in C$ satisfy $x \neq u$.
Hence, $u \neq x$.
If we had $x = v$, then we would have
$v = x \in C$ and therefore $u = v$
(by Lemma~\ref{lem.projection} \textbf{(a)}),
which would contradict $u \neq x = v$.
Hence, we cannot have $x = v$.
Thus, $x \neq v$.
This proves Lemma~\ref{lem.projection} \textbf{(b)}.

\textbf{(c)}
Lemma~\ref{lem.projection} \textbf{(b)} shows that
$x \neq v$. Hence, $d\tup{v, x}$ is well-defined.
Also, $d\tup{u, x}$ is well-defined (since $x \neq u$).

If $v \in C$, then Lemma~\ref{lem.projection} \textbf{(a)}
yields $u = v$ and therefore $d\tup{u, x} = d\tup{v, x}$.
Thus, Lemma~\ref{lem.projection} \textbf{(c)} is proven
if $v \in C$.
Hence, for the rest of the proof, we WLOG assume that
$v \notin C$.
Hence, the elements of $\proj_C\tup{v}$ are
the $c \in C$ that minimize
the distance $d\tup{v, c}$
(by the definition of $\proj_C\tup{v}$).
Thus, from $u \in \proj_C\tup{v}$,
we conclude that $u$ is a $c \in C$ that minimizes
the distance $d\tup{v, c}$.
Hence, $d\tup{v, u} \leq d\tup{v, x}$
(since $x \in C$).

Now, the points $u$ and $x$ belong to $C$,
while the point $v$ does not (since $v \notin C$).
Hence, $u$ and $x$ are distinct from $v$.
Therefore, the three points $u$, $x$ and $v$
are distinct (since $x \neq u$).
Hence,
\eqref{ultra} (applied to $a = u$, $b = x$ and $c = v$) yields
$d(u,x) \leqslant \max \set{d(u,v),d(x,v)} = \max \set{d(v,u),d(v,x)} = d(v,x)$
(since $d(v, u) \leq d(v, x)$).
This proves Lemma~\ref{lem.projection} \textbf{(c)}.
\end{proof}

\section{Perimeters and greedy $m$-permutations}

\subsection{The perimeter of an $m$-set}

For any finite subset $A\subset E$, we define its \defn{perimeter $\PER\tup{A}$} by
\[
\PER\tup{A} := \sum_{a\in A} w(a)+\sum_{\substack{\set{a,b}\subset A; \\ a \neq b}} d(a,b).
\]
The second sum here is taken over all
\textbf{unordered} pairs $a\ne b$ of distinct elements of $A$.
(This is well-defined, since $d\tup{a,b} = d\tup{b,a}$ for any distinct $a,b \in E$.)

\begin{example} \label{exa.5-point-1.per}
Let $(E, w, d)$ be as in Example~\ref{exa.5-point-1.d}.
Then,
\[
\PER\set{1, 2, 3} = w\tup{1} + w\tup{2} + w\tup{3} + \underbrace{d\tup{1,2}}_{=2} + \underbrace{d\tup{1,3}}_{=1} + \underbrace{d\tup{2,3}}_{=2} . 
\]
\end{example}

\subsection{Defining greedy $m$-permutations}

\begin{definition}
Let $C \subset E$ be any subset, and let $m$ be a nonnegative
integer.

A \defn{greedy $m$-permutation} of $C$ is a list $\tup{c_1, c_2, \ldots, c_m}$
of $m$ distinct elements of $C$
such that for each $i \in \set{1,2,\ldots,m}$
and each $x \in C \setminus \set{c_1, c_2, \ldots, c_{i-1}}$,
we have
\begin{align}
\PER\set{c_1, c_2, \ldots, c_i} \geqslant \PER\set{c_1, c_2, \ldots, c_{i-1}, x} .
\label{eq.def.greedy-enum.geq}
\end{align}
\end{definition}

Thus, roughly speaking, a greedy $m$-permutation is an
ordered sample of $m$ distinct elements of $C$ such that
at each step of the sampling procedure, the new element
is chosen in such a way as to maximize the perimeter
of the sample.
This procedure can be viewed as a greedy algorithm to
construct an $m$-subset of $C$ that has maximum
perimeter. As we shall see in Theorem~\ref{thm.2a},
this algorithm indeed succeeds at constructing such a
subset.

\subsection{Examples of greediness}

\begin{example} \label{exa.5-point-1.greedy}
Let $(E, w, d)$ be as in Example~\ref{exa.5-point-1.d}.
Assume that $w\tup{a} = 0$ for all $a \in E$.

Then, $\tup{1, 2}$, $\tup{2, 1}$ and $\tup{5, 4}$ (and
several others) are greedy $2$-permutations of $E$.
Actually, a pair $\tup{i, j}$ of elements of $E$
is a greedy $2$-permutation of $E$ if and only if
$i \not\equiv j \mod 2$.

Also, $\tup{1, 3}$ is a greedy $2$-permutation of
$\set{1, 3, 5}$,
but not of $E$ (since $\PER\set{1, 3} < \PER\set{1, 2}$).

Also, $\tup{1, 2, 3, 4, 5}$ is a greedy $5$-permutation
of $E$, but $\tup{1, 2, 3, 5, 4}$ is not
(since $\PER\set{1, 2, 3, 5} < \PER\set{1, 2, 3, 4}$).
\end{example}

\begin{example} \label{exa.6-point-1.greedy}
Let $E$ be the set $\set{1, 2, 3, 4, 5, 6}$.
Fix five reals
$\alpha, \lambda, \kappa, \varepsilon, \delta$
such that $\lambda$ and $\kappa$ are both smaller than $\alpha$
and both larger than each of $\varepsilon$ and $\delta$.
For any distinct $a, b \in E$, we define the distance $d\tup{a, b}$ by the following rule:
\begin{itemize}
\item If $a \not\equiv b \mod 2$, then $d\tup{a, b} = \alpha$.
\item If $a = 1$ and $b \in \set{3, 5}$, then $d\tup{a, b} = \lambda$.
\item If $a = 2$ and $b \in \set{4, 6}$, then $d\tup{a, b} = \kappa$.
\item If $a = 3$ and $b = 5$, then $d\tup{a, b} = \varepsilon$.
\item If $a = 4$ and $b = 6$, then $d\tup{a, b} = \delta$.
\item Otherwise, $d \tup{a, b} = d \tup{b, a}$.
\end{itemize}
Set $w\tup{a} = 0$ for all $a \in E$.

It is easy to check that $\tup{E, w, d}$ is an ultra triple.

The pair $\tup{1, 2}$ is always a greedy $2$-permutation.

The $4$-tuple $\tup{1, 2, 3, 4}$ is a greedy $4$-permutation if and only if $\lambda \geq \kappa$.
The $4$-tuple $\tup{1, 2, 4, 3}$ is a greedy $4$-permutation if and only if $\kappa \geq \lambda$.

The $5$-tuple $\tup{1, 2, 3, 4, 5}$ is a greedy $5$-permutation if and only if $\lambda \geq \kappa$ and $\lambda + \varepsilon \geq \kappa + \delta$.

The $5$-set $\set{1, 2, 3, 4, 5}$ has maximum perimeter among all $5$-sets if and only if $\lambda + \varepsilon \geq \kappa + \delta$.

This example illustrates that greedy permutations and maximum-perimeter sets depend not just on the order relations between the distances of the points, but also on the order relations between sums of these distances.
\end{example}

\begin{example} \label{exa.p-adic.greedy}
For this example, we fix a prime number $p$
and a nonnegative integer $m$.
We let $E$ be any subset of $\ZZ$ that contains
$1, 2, \ldots, m$.
We define $d : E \timesu E \to \RR$ as in
Example~\ref{exa.p-adic.d}.
We define
$d' : E \timesu E \to \RR$ as in
Example~\ref{exa.p-adic.d'}.
We define $w : E \to \RR$ by setting
$w \tup{e} = 0$ for all $e \in E$.

Then, $\tup{1, 2, \ldots, m}$ is a greedy
$m$-permutation of $E$ both for the ultra triple
$\tup{E, w, d}$ and for the ultra triple $\tup{E, w, d'}$.
\end{example}

We relegate the proof of this claim to
Section~\ref{sect.app-p}, as we shall not use
it in what follows.

\begin{example} \label{exa.p-adic.greedy-weird}
Example~\ref{exa.p-adic.greedy} might suggest that
the ultra triples $\tup{E, w, d}$ and $\tup{E, w, d'}$
(defined in that example) have the same greedy
$m$-permutations in general.
This is not the case.
For instance, set $p = 2$ and
$E = \set{0, 1, 2, 9, 17, 128}$.
Define $d$, $d'$ and $w$ as in
Example~\ref{exa.p-adic.greedy}.

Now it is easy to check that
$\tup{2, 9, 17, 0, 1}$ is a greedy $5$-permutation
for $\tup{E, w, d'}$ but not for $\tup{E, w, d}$,
while $\tup{2, 9, 17, 0, 128}$ is a greedy $5$-permutation
for $\tup{E, w, d}$ but not for $\tup{E, w, d'}$.

Moreover, the $5$-set $\set{2, 9, 17, 0, 1}$ has
maximum perimeter for $\tup{E, w, d'}$ but not for
$\tup{E, w, d}$,
while the $5$-set
$\set{2, 9, 17, 0, 128}$ has maximum perimeter
for $\tup{E, w, d}$ but not for $\tup{E, w, d'}$.
\end{example}

\subsection{Basic properties of greediness}

We will use the following shorthand notations:
If $S$ is a subset of $E$, and if $e \in E$, then
\defn{$S \cup e$} and \defn{$S \setminus e$} will stand
for the subsets $S \cup \set{e}$
and $S \setminus \set{e}$, respectively.
Set operations like $\cup$ and $\setminus$ shall be
read in a left-associative way; thus, e.g., the
expression ``$S \cup e \setminus f$'' shall be
understood as $\tup{S \cup e} \setminus f$.

Let us observe some easy consequences of the definition
of greedy permutations
(which will be later used without mention):

\begin{proposition} \label{prop.greed.trivial}
Let $C$ be a subset of $E$.
Let $m$ be a nonnegative integer.

\begin{enumerate}

\item[\textbf{(a)}] If a greedy $m$-permutation of $C$
    exists,
    then $m \leq \abs{C}$.

\item[\textbf{(b)}] If
    $\tup{c_1, c_2, \ldots, c_m}$
    is a greedy $m$-permutation of $C$,
    then $\set{c_1, c_2, \ldots, c_k}$ is a $k$-subset of $C$
    for each $k \in \set{0, 1, \ldots, m}$.

\item[\textbf{(c)}] If
    $\tup{c_1, c_2, \ldots, c_m}$
    is a greedy $m$-permutation of $C$,
    then $\set{c_1, c_2, \ldots, c_{k-1}}
    = \set{c_1, c_2, \ldots, c_k} \setminus c_k$
    for each $k \in \set{1, 2, \ldots, m}$.

\item[\textbf{(d)}] If $\abs{C} = m$,
    then any greedy $m$-permutation of $C$
    must be a list of all the $m$ elements of $C$.

\item[\textbf{(e)}] If $C$ is finite and satisfies $m \leq \abs{C}$,
    then there exists a greedy $m$-permutation
    of $C$.

\end{enumerate}
\end{proposition}

\begin{proof}
\textbf{(a)}
A greedy $m$-permutation consists of $m$ distinct elements of $C$
(by definition).
Thus, if it exists, then $C$ must have at least $m$ elements,
so that $m \leq \abs{C}$.
This proves part \textbf{(a)}.

\textbf{(b)} Let
$\tup{c_1, c_2, \ldots, c_m}$ be a greedy $m$-permutation of $C$.
Then, $c_1, c_2, \ldots, c_m$ are distinct (by the definition of
a greedy $m$-permutation).
Hence, $\set{c_1, c_2, \ldots, c_k}$ is a $k$-element set
whenever $k \in \set{0, 1, \ldots, m}$.
This $k$-element set is furthermore a subset of $C$
(since $c_1, c_2, \ldots, c_m \in C$),
hence a $k$-subset of $C$.
This proves part \textbf{(b)}.

\textbf{(c)}
This is obvious, since $c_1, c_2, \ldots, c_m$ are distinct.

\textbf{(d)} Assume that $\abs{C} = m$.
Any greedy $m$-permutation of $C$ must be a list of $m$
distinct elements of $C$,
and therefore must be a list of all the $m$ elements of $C$
(since $C$ has only $m$ elements in total).
This proves part \textbf{(d)}.

\textbf{(e)}
Assume that $C$ is finite and satisfies $m \leq \abs{C}$.
We can then construct a greedy $m$-permutation
$\tup{c_1, c_2, \ldots, c_m}$ of $C$ according to the following
recursive procedure:

\begin{itemize}
\item For each $i = 1, 2, \ldots, m$, we assume that
      $c_1, c_2, \ldots, c_{i-1}$ have already been defined;
      we then choose an element
      $c_i \in C \setminus \set{c_1, c_2, \ldots, c_{i-1}}$
      that maximizes the perimeter $\PER\set{c_1, c_2, \ldots, c_i}$.
      (If there are several such elements, then we choose any
      of them.)
\end{itemize}

This procedure can be carried out, since at each step
we can find an element $c_i \in C \setminus \set{c_1, c_2, \ldots, c_{i-1}}$
that maximizes the perimeter $\PER\set{c_1, c_2, \ldots, c_i}$.
(Indeed, the set $C \setminus \set{c_1, c_2, \ldots, c_{i-1}}$
is nonempty because $\abs{C} \geq m \geq i > i-1 \geq \abs{\set{c_1, c_2, \ldots, c_{i-1}}}$;
furthermore, this set is finite, and thus at least one of its
elements will maximize the perimeter in question.)

Clearly, the result of this procedure is an $m$-tuple
$\tup{c_1, c_2, \ldots, c_m}$ of elements of $C$.
The entries $c_1, c_2, \ldots, c_m$ of this $m$-tuple are
distinct (since each $c_i$ is chosen to be an element of
$C \setminus \set{c_1, c_2, \ldots, c_{i-1}}$, and thus
is distinct from all of $c_1, c_2, \ldots, c_{i-1}$),
and furthermore it satisfies
\eqref{eq.def.greedy-enum.geq} for each
$i \in \set{1,2,\ldots,m}$ and each
$x \in C \setminus \set{c_1, c_2, \ldots, c_{i-1}}$
(due to how $c_i$ is chosen).
Thus, this $m$-tuple $\tup{c_1, c_2, \ldots, c_m}$ is
a greedy $m$-permutation of $C$.

Hence, a greedy $m$-permutation of $C$ exists.
This proves part \textbf{(e)}.
\end{proof}

The procedure used in the proof of
Proposition~\ref{prop.greed.trivial} \textbf{(e)}
also works for
infinite $C$ as long as the maxima exist.


Proposition~\ref{prop.greed.trivial} \textbf{(e)}
can be generalized further:
Any greedy $n$-permutation with $n \leq m$
can be extended to a greedy $m$-permutation:

\begin{proposition} \label{prop.greedy.extend}
Let $m$ and $n$ be integers such that $m \geq n \geq 0$.
Let $C$ be a finite subset of $E$ such that $\abs{C} \geq m$.

If $\tup{c_1, c_2, \ldots, c_n}$ is a greedy $n$-permutation of $C$,
then we can find $m-n$ elements $c_{n+1}, c_{n+2}, \ldots, c_m$ of $C$
such that $\tup{c_1, c_2, \ldots, c_m}$ is a greedy $m$-permutation of $C$.
\end{proposition}

\begin{proof}
Recall the procedure we used to construct
a greedy $m$-permutation $\tup{c_1, c_2, \ldots, c_m}$
in the proof of
Proposition~\ref{prop.greed.trivial} \textbf{(e)}.
The same procedure can be used here, as long as
we replace ``For each $i = 1, 2, \ldots, m$''
by ``For each $i = n+1, n+2, \ldots, m$''
(since $c_1, c_2, \ldots, c_n$ are already defined).
\end{proof}

Conversely, of course, we have the following obvious fact:

\begin{proposition} \label{prop.greedy.prefix}
Let $C$ be a subset of $E$.
Let $m$ and $n$ be integers such that $m \geq n \geq 0$.

If $\tup{c_1, c_2, \ldots, c_m}$ is a greedy $m$-permutation of $C$,
then $\tup{c_1, c_2, \ldots, c_n}$ is a greedy $n$-permutation of $C$.
\end{proposition}

\section{The main theorems}

We now state two central theorems for this paper:

\begin{theorem} \label{thm.2a}

Let $C\subset E$ be any subset,
and let $m$ be a nonnegative integer.

Let $\tup{c_1, c_2, \ldots, c_m}$ be
any greedy $m$-permutation of $C$.

Then, for each $k \in \set{0,1,\dots,m}$, the
set $\set{c_1, c_2, \ldots, c_k}$ has maximum perimeter
among all $k$-subsets of $C$.

\end{theorem}

\begin{theorem} \label{thm.2b}

Let $C\subset E$ be any finite subset,
and let $m$ be a nonnegative integer such that $\abs{C} \geq m$.
Let $k \in \set{0,1,\dots,m}$.

Let $A$ be a $k$-subset of $C$ having maximum perimeter (among the $k$-subsets of $C$).
Then, there exists a greedy $m$-permutation
$\tup{v_1, v_2, \ldots, v_m}$ of $C$
such that $A = \set{v_1, v_2, \ldots, v_k}$.

\end{theorem}

\begin{proof}[Proof of Theorem \ref{thm.2a}.]
The points $c_1, c_2, \ldots, c_m$ are distinct
(by the definition of a greedy $m$-permutation).

Fix $k \in \set{0, 1, \ldots, m}$.
Then, $\set{c_1, c_2, \ldots, c_k}$ is a $k$-subset of $C$
(by Proposition~\ref{prop.greed.trivial} \textbf{(b)}).
It remains to prove that every $k$-subset $A$ of $C$
satisfies $\PER\tup{A} \leq \PER\set{c_1, c_2, \ldots, c_k}$.

Let $A$ be any $k$-subset of $C$.
We shall show that $\PER\tup{A} \leq \PER\set{c_1, c_2, \ldots, c_k}$.
This will clearly prove Theorem~\ref{thm.2a}.

We define $k$ elements $v_1, v_2, \ldots, v_k \in A$
recursively as follows:
For each $i=1,2,\ldots,k$, we let $v_i$ be a projection
of $c_i$ onto $A \setminus \set{v_1,v_2,\ldots, v_{i-1}}$
(assuming that $v_1, v_2, \ldots, v_{i-1}$ have already
been constructed).\footnote{Thus, in particular, $v_1$ is
a projection of $c_1$ onto $A$.}
(These projections $v_i$ exist because of Proposition~\ref{prop.proj-ex}%
\footnote{In more detail:
Let $i \in \set{1,2,\ldots,k}$, and assume that
$v_1, v_2, \ldots, v_{i-1}$ have already
been constructed; we must prove that $v_i$ exists.
We have $\abs{\set{v_1,v_2,\ldots,v_{i-1}}} \leq i-1
< i \leq k = \abs{A}$; thus, the set
$A \setminus \set{v_1,v_2,\ldots, v_{i-1}}$ is nonempty.
Since this set is furthermore finite, we thus conclude
(by Proposition~\ref{prop.proj-ex}) that there exists
a projection of $c_i$ onto $A \setminus \set{v_1,v_2,\ldots, v_{i-1}}$.
In other words, $v_i$ exists.};
they may be non-unique, but any choice is fine.)

Thus, we get $k$ elements $v_1, v_2, \ldots, v_k$ of $A$.
These $k$ elements are distinct (since each $v_i$ has been constructed
to belong to $A \setminus \set{v_1,v_2,\ldots, v_{i-1}}$).
Since $\abs{A} = k$, these $k$ distinct elements must cover the whole set $A$.
Hence, $\tup{v_1,v_2,\ldots,v_k}$ is an enumeration of $A$
(that is, a list of distinct elements of $A$ such that
$A = \set{v_1, v_2, \ldots, v_k}$).

Let $j \in \set{1, 2, \ldots, k}$. Then,
\begin{align}
v_j \notin \set{c_1, c_2, \ldots, c_{j-1}} .
\label{pf.thm.2a.notin}
\end{align}

[\textit{Proof of \eqref{pf.thm.2a.notin}:}
Assume the contrary. Thus, $v_j \in \set{c_1, c_2, \ldots, c_{j-1}}$,
so that $v_j = c_i$ for some $i<j$. Consider this $i$.
Hence, $c_i = v_j \in \set{v_i,v_{i+1},\dots,v_k}
= A \setminus \set{v_1, v_2, \ldots, v_{i-1}}$
(since $\tup{v_1,v_2,\ldots,v_k}$ is an enumeration of $A$).
But our recursive definition of $v_i$ shows that
$v_i$ is a projection of $c_i$ onto the set
$A \setminus \set{v_1, v_2, \ldots, v_{i-1}}$.
Hence, Lemma~\ref{lem.projection} \textbf{(a)}
(applied to $A \setminus \set{v_1, v_2, \ldots, v_{i-1}}$,
$c_i$ and $v_i$ instead of $C$, $v$ and $u$)
yields $v_i = c_i$
(since $c_i \in A \setminus \set{v_1, v_2, \ldots, v_{i-1}}$).
Hence, $v_i = c_i = v_j$, whence $i = j$
(since $v_1, v_2, \ldots, v_k$ are distinct).
But this contradicts $i < j$.
This contradiction shows that our assumption was false,
and thus \eqref{pf.thm.2a.notin} is proven.]

Hence,
$v_j \in
A \setminus \set{c_1, c_2, \ldots, c_{j-1}}
\subseteq
C \setminus \set{c_1, c_2, \ldots, c_{j-1}}$
and therefore
$\PER\set{c_1, c_2, \ldots, c_{j-1}, v_j}
\leqslant
\PER\set{c_1, c_2, \ldots, c_j}$
by the definition of a greedy $m$-permutation
(specifically, by \eqref{eq.def.greedy-enum.geq}).

But $c_1, c_2, \ldots, c_{j-1}, v_j$ are distinct
(since $c_1, c_2, \ldots, c_m$ are distinct,
and since $v_j \notin
\set{c_1, c_2, \ldots, c_{j-1}}$),
and thus
\begin{align*}
&\PER\set{c_1, c_2, \ldots, c_{j-1}} +
w(v_j)+\sum_{i=1}^{j-1} d(c_i,v_j)
= \PER\set{c_1, c_2, \ldots, c_{j-1}, v_j} \\
& \leqslant \PER\set{c_1, c_2, \ldots, c_j}
= \PER\set{c_1, c_2, \ldots, c_{j-1}} +
w(c_j)+\sum_{i=1}^{j-1} d(c_i,c_j)
\end{align*}
(since $c_1, c_2, \ldots, c_j$ are distinct).
After cancelling equal terms, this rewrites as
\begin{equation}\label{algo}
w(v_j)+\sum_{i=1}^{j-1} d(c_i,v_j) \leqslant
w(c_j)+\sum_{i=1}^{j-1} d(c_i,c_j) .
\end{equation}

Furthermore, for each $i \in \set{1, 2, \ldots, j-1}$,
we have $j > i$ and thus
$v_j \in \set{v_i, v_{i+1}, \ldots, v_k}
= A \setminus \set{v_1, v_2, \ldots, v_{i-1}}$
(since $\tup{v_1,v_2,\ldots,v_k}$ is an enumeration of $A$)
and $v_j \neq v_i$ (since $v_1, v_2, \ldots, v_k$
are distinct).
Hence,
Lemma \ref{lem.projection} \textbf{(c)} (applied to
$A \setminus \set{v_1, v_2, \ldots, v_{i-1}}$, $c_i$, $v_i$ and $v_j$
instead of $C$, $v$, $u$ and $x$) yields
\begin{align}
\label{algo2}
d(v_i,v_j)\leqslant d(c_i,v_j)
\text{ for each } i \in \set{1, 2, \ldots, j-1}
\end{align}
(since $v_i$ is a projection of $c_i$ onto the set
$A \setminus \set{v_1, v_2, \ldots, v_{i-1}}$,
whereas $v_j \in A \setminus \set{v_1, v_2, \ldots, v_{i-1}}$
and $v_j \neq v_i$).

Now, forget that we fixed $j$.
We thus have proven \eqref{algo} and \eqref{algo2}
for each $j \in \set{1,2,\ldots,k}$.

But $\tup{v_1, v_2, \ldots, v_k}$ is an enumeration of $A$;
thus,
\begin{align*}
\PER\tup{A}
&= \sum_{j=1}^k w(v_j) + \sum_{1\leq i<j\leq k} d(v_i, v_j)
 =\sum_{j=1}^k \left(w(v_j)+\sum_{i=1}^{j-1} d(v_i,v_j)\right) \\
& \leqslant
\sum_{j=1}^k \left(w(v_j)+\sum_{i=1}^{j-1} d(c_i,v_j)\right)
\qquad \tup{\text{by \eqref{algo2}}}
\\
& \leqslant
\sum_{j=1}^k \left(w(c_j)+\sum_{i=1}^{j-1} d(c_i,c_j)\right)
\qquad \tup{\text{by \eqref{algo}}}
\\
&= \sum_{j=1}^k w(c_j) + \sum_{1\leq i<j\leq k} d(c_i, c_j)
 = \PER\set{c_1, c_2, \ldots, c_k}
\end{align*}
(since $c_1, c_2, \ldots, c_k$ are distinct).
This proves Theorem~\ref{thm.2a}.
\end{proof}

\begin{proof}[Proof of Theorem \ref{thm.2b}.]
Proposition~\ref{prop.greed.trivial} \textbf{(e)} shows that
there exists a greedy $m$-permutation of $C$
(since $C$ is finite and since $\abs{C} \geq m$).
Choose such a greedy $m$-permutation $\tup{c_1,c_2,\dots,c_m}$ of $C$.
Then, Theorem~\ref{thm.2a} shows that
the set $\set{c_1, c_2, \ldots, c_k}$
has maximum perimeter among all $k$-subsets of $C$.
Hence, $\PER\tup{A} = \PER\set{c_1, c_2, \ldots, c_k}$
(since the set $A$ also has maximum perimeter among them).

Construct an enumeration $\tup{v_1, v_2, \ldots, v_k}$ of $A$
as in the proof of Theorem~\ref{thm.2a} above.
In our above proof of Theorem~\ref{thm.2a}, we have
proven the inequalities \eqref{algo} and \eqref{algo2} for all
$j \in \set{1, 2, \ldots, k}$.
But by adding together all these inequalities, we have
obtained the inequality $\PER\tup{A} \leq \PER\set{c_1, c_2, \ldots, c_k}$,
which must be an equality (since $\PER\tup{A} = \PER\set{c_1, c_2, \ldots, c_k}$).
Thus, all the inequalities \eqref{algo} and \eqref{algo2}
must become equalities
(because if adding together a bunch of inequalities
produces an equality, then all the inequalities
must themselves be equalities).

Hence, for all $j \in \set{1,2,\dots,k}$, we have
\begin{align}
\label{zamena}
d(v_i,v_j)= d(c_i,v_j) \quad \text{for each } i \in \set{1,2,\ldots,j-1}
\end{align}
(since the inequalities \eqref{algo2} become equalities),
and thus
\begin{align}
w(v_j)+\sum_{i=1}^{j-1} d(v_i,v_j)
&=
w(v_j)+\sum_{i=1}^{j-1} d(c_i,v_j)
\nonumber \\
&=
w(c_j)+\sum_{i=1}^{j-1} d(c_i,c_j)
\label{zamena1}
\end{align}
(since the inequality \eqref{algo} becomes an equality).

Now, fix $p \in \set{1, 2, \ldots, k}$.
Hence, $p \leq k \leq m$.
The points $c_1, c_2, \ldots, c_m$ are distinct
(since $\tup{c_1, c_2, \ldots, c_m}$ is a
greedy $m$-permutation).
Thus, the points $c_1, c_2, \ldots, c_p$ are distinct.
Also, the points $v_1, v_2, \ldots, v_p$ are distinct
(since $v_1, v_2, \ldots, v_k$ are distinct); hence,
the definition of a perimeter yields
\begin{align}
\PER\set{v_1, v_2, \ldots, v_p}
&= \sum_{j=1}^p w\tup{v_j} + \sum_{1 \leq i < j \leq p} d\tup{v_i, v_j} \nonumber\\
&= \sum_{j=1}^p \tup{w\tup{v_j} + \sum_{i=1}^{j-1} d\tup{v_i, v_j}} \nonumber\\
&= \sum_{j=1}^p \tup{w(c_j)+\sum_{i=1}^{j-1} d(c_i,c_j)}
 \qquad \tup{\text{by \eqref{zamena1}}} \nonumber\\
&= \sum_{j=1}^p w\tup{c_j} + \sum_{1 \leq i < j \leq p} d\tup{c_i, c_j} \nonumber\\
&= \PER\set{c_1, c_2, \ldots, c_p}
\label{pf.thm.2b.PER=PER}
\end{align}
(since the points $c_1, c_2, \ldots, c_p$ are distinct).

But Theorem~\ref{thm.2a} (applied to $p$ instead of $k$)
shows that
the set $\set{c_1, c_2, \ldots, c_p}$ has maximum perimeter
among all $p$-subsets of $C$.
Hence, the set $\set{v_1, v_2, \ldots, v_p}$ must also
have maximum perimeter among all $p$-subsets of $C$
(because \eqref{pf.thm.2b.PER=PER} shows that this latter
set has the same perimeter as the former set).
Hence, for each $x \in C \setminus \set{v_1, v_2, \ldots, v_{p-1}}$,
we have
$\PER\set{v_1, v_2, \ldots, v_p} \geqslant
\PER\set{v_1, v_2, \ldots, v_{p-1}, x}$.

Now, forget that we fixed $p$.
We thus have shown that for each $p \in \set{1, 2, \ldots, k}$
and each $x \in C \setminus \set{v_1, v_2, \ldots, v_{p-1}}$,
we have
\[
\PER\set{v_1, v_2, \ldots, v_p} \geqslant
\PER\set{v_1, v_2, \ldots, v_{p-1}, x} .
\]
In other words, $\tup{v_1, v_2, \ldots, v_k}$ is a greedy
$k$-permutation of $C$
(since $v_1, v_2, \ldots, v_k$ are distinct).
Hence, Proposition~\ref{prop.greedy.extend}
(applied to $k$ and $v_i$ instead of $n$ and $c_i$)
shows that we can find $m-k$ elements
$v_{k+1}, v_{k+2}, \ldots, v_m$ of $C$ such that
$\tup{v_1, v_2, \ldots, v_m}$ is a greedy
$m$-permutation of $C$.
This proves Theorem~\ref{thm.2b}
(since $A = \set{v_1, v_2, \ldots, v_k}$).
\end{proof}

\section{The $\nuo_k\tup{C}$ invariants}

We shall next prove two corollaries of the above results
that resemble (and, as we will later see,
generalize) \cite[Theorem 1]{Bharga97}
and \cite[Lemma 2]{Bharga97}.

\begin{corollary} \label{cor.bh-t1}
Let $C \subset E$ be any subset.
Let $m$ be a nonnegative integer.
Let $k \in \set{1, 2, \ldots, m}$.
If $\tup{c_1, c_2, \ldots, c_m}$ is a greedy
$m$-permutation of $C$, then the number
\[
w\tup{c_k} + \sum_{i=1}^{k-1} d\tup{c_i, c_k}
\]
does not depend on the choice of this $m$-permutation
(but only depends on $k$ and on $C$).
\end{corollary}

\begin{proof}
Let $\tup{c_1, c_2, \ldots, c_m}$ be a greedy
$m$-permutation of $C$.
Hence, Theorem~\ref{thm.2a} shows that the set
$\set{c_1, c_2, \ldots, c_k}$
has maximum perimeter among all $k$-subsets of $C$.
In other words, $\PER\set{c_1, c_2, \ldots, c_k}$
equals the maximum possible perimeter of a $k$-subset of $C$.
Likewise,
$\PER\set{c_1, c_2, \ldots, c_{k-1}}$
equals the maximum possible perimeter of a $\tup{k-1}$-subset of $C$.
Hence, both numbers $\PER\set{c_1, c_2, \ldots, c_k}$ and
$\PER\set{c_1, c_2, \ldots, c_{k-1}}$
depend only on $k$ and $C$.
Thus, so does their difference
$\PER\set{c_1, c_2, \ldots, c_k} - \PER\set{c_1, c_2, \ldots, c_{k-1}}$.
In view of
\[
\PER\set{c_1, c_2, \ldots, c_k} - \PER\set{c_1, c_2, \ldots, c_{k-1}}
= w\tup{c_k} + \sum_{i=1}^{k-1} d\tup{c_i, c_k}
\]
(which is because $c_1, c_2, \ldots, c_k$ are distinct),
this rewrites as follows:
The number
$w\tup{c_k} + \sum_{i=1}^{k-1} d\tup{c_i, c_k}$
depends only on $k$ and $C$.
\end{proof}

From now on, the number
\[
w\tup{c_k} + \sum_{i=1}^{k-1} d\tup{c_i, c_k}
\]
in Corollary~\ref{cor.bh-t1} will be denoted by
\defn{$\nuo_k\tup{C}$}.

\begin{corollary} \label{cor.bh-l2}
Let $C \subset E$ be any subset.
Let $m$ be a nonnegative integer.
Let $k \in \set{1, 2, \ldots, m}$.
Let $\tup{c_1, c_2, \ldots, c_m}$ be a greedy
$m$-permutation of $C$.
Let $j \in \set{1, 2, \ldots, k}$.
Then,
\begin{equation}
\nuo_k\tup{C}
\leq w\tup{c_j} + \sum_{i \in \set{1, 2, \ldots, k} \setminus \set{j}} d\tup{c_i, c_j} .
\label{eq.cor.bh-l2.claim}
\end{equation}
\end{corollary}

\begin{proof}
We can see (as in the proof of Corollary~\ref{cor.bh-t1})
that $\PER\set{c_1, c_2, \ldots, c_{k-1}}$
equals the maximum possible perimeter of a $\tup{k-1}$-subset of $C$.
Thus,\footnote{Here, the hat over the $c_j$ signifies that
$c_j$ is omitted from the list.}
\[
\PER\set{c_1, c_2, \ldots, c_{k-1}}
\geq \PER\set{c_1, c_2, \ldots, \widehat{c_j}, \ldots, c_k}
\]
(here, we have used the fact that
$\set{c_1, c_2, \ldots, \widehat{c_j}, \ldots, c_k}$
is a $\tup{k-1}$-subset of $C$,
which is because $c_1, c_2, \ldots, c_m$ are distinct).
In view of
\begin{align*}
\PER\set{c_1, c_2, \ldots, c_{k-1}}
&= \PER\set{c_1, c_2, \ldots, c_k}
   - \underbrace{\tup{w\tup{c_k} + \sum_{i=1}^{k-1} d\tup{c_i, c_k}}}_{\substack{= \nuo_k\tup{C} \\ \tup{\text{by the definition of $\nuo_k\tup{C}$}}}} \\
& \qquad \tup{\text{since $c_1, c_2, \ldots, c_k$ are distinct}} \\
&= \PER\set{c_1, c_2, \ldots, c_k} - \nuo_k\tup{C}
\end{align*}
and
\begin{align*}
& \PER\set{c_1, c_2, \ldots, \widehat{c_j}, \ldots, c_k} \\
&= \PER\set{c_1, c_2, \ldots, c_k} - \tup{w\tup{c_j} + \sum_{i \in \set{1, 2, \ldots, k} \setminus \set{j}} d\tup{c_i, c_j}} \\
& \qquad \tup{\text{since $c_1, c_2, \ldots, c_k$ are distinct}} ,
\end{align*}
this rewrites as
\begin{align*}
& \PER\set{c_1, c_2, \ldots, c_k} - \nuo_k\tup{C} \\
&\geq \PER\set{c_1, c_2, \ldots, c_k} - \tup{w\tup{c_j} + \sum_{i \in \set{1, 2, \ldots, k} \setminus \set{j}} d\tup{c_i, c_j}} .
\end{align*}
Subtracting this inequality from the obvious equality
$\PER\set{c_1, c_2, \ldots, c_k} = \PER\set{c_1, c_2, \ldots, c_k}$,
we obtain precisely \eqref{eq.cor.bh-l2.claim}.
\end{proof}

\section{The greedoid}

Throughout this section, we assume that the set $E$ is finite.

\subsection{Defining greedoids and strong greedoids}

We shall now recall the definition of a ``greedoid'':

A collection\footnote{The word ``\defn{collection}''
just means ``set'', but will be used exclusively
for sets of sets.}
$\calF \subset 2^E$ of subsets of a finite set $E$ is
said to be a
\defn{greedoid}\footnote{More precisely, the sets in
the collection are said to be the \defn{feasible sets}
of a greedoid. We will, however, just say that the
collection is a greedoid.} (on the ground set $E$)
if it satisfies the following three axioms:

\begin{enumerate}

\item[\textbf{(i)}] We have $\varnothing \in \calF$.

\item[\textbf{(ii)}] If $B\in \calF$ satisfies $\abs{B}>0$,
then there exists
$b\in B$ such that $B\setminus b\in \calF$.

\item[\textbf{(iii)}] If
$A,B\in \calF$ satisfy
$\abs{B} = \abs{A} + 1$, then
there exists $b\in B\setminus A$
such that $A\cup b\in \calF$.

\end{enumerate}


We refer to \cite{KoLoSc91} for a book-length treatment of
greedoids.
Our above definition of a greedoid appears implicitly in
\cite[Section IV.1]{KoLoSc91}
(indeed, our axioms \textbf{(i)} and \textbf{(iii)} correspond to the
conditions (1.4) and (1.6) in \cite[Section IV.1]{KoLoSc91},
while our axioms \textbf{(i)} and \textbf{(ii)} make $\tup{E, \calF}$ into
what is called an \defn{accessible set system} in \cite{KoLoSc91}).

There are several classes of greedoids having additional
properties besides the above three axioms.
(See \cite{KoLoSc91} for an overview.)
Let us define one of these classes -- that of ``strong greedoids''
(also known as ``Gauss greedoids''):

A greedoid $\calF$ on a ground set $E$ is said to be
a \defn{strong greedoid} if it satisfies the following
axiom:

\begin{enumerate}

\item[\textbf{(iv)}]
If $A, B \in \calF$ satisfy $\abs{B} = \abs{A} + 1$,
then there exists some $x \in B \setminus A$ such that
$A \cup x \in \calF$ and $B \setminus x \in \calF$.

\end{enumerate}

This definition of strong greedoids appears in
\cite{BrySha99} (where the above axiom \textbf{(iv)}
appears as property G(3)$'$).
Note that axiom \textbf{(iv)} is clearly stronger than
axiom \textbf{(iii)}.
The theorem in Section 2 of \cite{BrySha99} says that
strong greedoids are the same as Gauss greedoids (one of
the classes of greedoids studied in \cite{KoLoSc91}).
See \cite[Section IX.4]{KoLoSc91} for further properties
and characterizations of Gauss greedoids.

\subsection{The Bhargava greedoid}

The following theorem shows that a greedoid can be
obtained from any ultra triple $\tup{E, w, d}$:

\begin{theorem}\label{thm.greedoid}
Let $\calF$ denote the collection
of subsets $A\subset E$ that have maximum
perimeter among all $\abs{A}$-sets:
\[
\calF
= \set{ A \subset E \ \mid \  \PER\tup{A} \geqslant \PER\tup{B}
         \text{ for all } B \subset E \text{ satisfying } \abs{B} = \abs{A} }.
\]
Then $\calF$ is a strong greedoid on the ground set $E$.
\end{theorem}

We call this $\calF$ the \defn{Bhargava greedoid}
of the ultra triple $(E,w,d)$.

\begin{example} \label{exa.5-point-1.greedoid}
Let $(E, w, d)$ be as in Example~\ref{exa.5-point-1.d}.
Assume that $w\tup{a} = 0$ for all $a \in E$.

Then, the collection $\calF$ in Theorem~\ref{thm.greedoid}
contains $\set{1, 2, 3}$ and $\set{1, 2, 3, 4, 5}$
but not $\set{1, 2, 3, 5}$.

Theorem~\ref{thm.greedoid} says that this collection is a
strong greedoid; hence, axiom \textbf{(iii)} in the definition
of a greedoid yields that for any $A,B\in \calF$
satisfying
$\abs{B} = \abs{A} + 1$, there exists $b\in B\setminus A$
such that $A\cup b\in \calF$.
For example, if we pick $A = \set{1, 2, 5}$ and
$B = \set{2, 3, 4, 5}$, then this says that there exists
$b \in \set{3, 4}$ such that $\set{1, 2, 5, b} \in \calF$.
And indeed, $b = 4$ works (though $b = 3$ does not).
\end{example}

\begin{example} \label{exa.non-transp-greedoid}
Let $p = 3$ and $E = \set{0, 1, 2, 3, 4, 5, 6, 12}$.
Define the distance function $d : E \timesu E \to \RR$ as
in Example~\ref{exa.p-adic.d}.
Set $w\tup{e} = 0$ for all $e \in E$.

Then, the collection $\calF$ in Theorem~\ref{thm.greedoid}
contains $\set{0, 1, 2}$ and $\set{0, 1, 2, 3}$ and
$\set{0, 1, 2, 6}$ and $\set{0, 1, 2, 4, 5, 6, 12}$
but not $\set{0, 1, 2, 3, 6}$ and not
$\set{0, 1, 2, 3, 4, 5, 12}$.
\end{example}

For readers familiar with the alternative description
of greedoids as hereditary languages (see, e.g.,
\cite[Section IV.1]{KoLoSc91}), we note in passing
that the language corresponding to the greedoid
$\calF$ in Theorem~\ref{thm.greedoid} is precisely
the set of greedy $m$-permutations for $m \geq 0$.
This observation will not be used in what follows,
but helps illuminate the proofs.

Our proof of Theorem~\ref{thm.greedoid} will rely on the
following lemma (inspired by \cite[Theorem 3.2]{MoSeSt06}):

\begin{lemma} \label{lem.strong-greedoid.AB}
Let $A$ and $B$ be two subsets of $E$ such that
$\abs{B} = \abs{A} + 1$.

Then, there exists a $u \in B \setminus A$
satisfying
\begin{align}
\PER\tup{B \setminus u} + \PER\tup{A \cup u}
\geq \PER\tup{A} + \PER\tup{B} .
\label{eq.lem.strong-greedoid.AB.per}
\end{align}
\end{lemma}

\begin{proof}[Proof of Lemma~\ref{lem.strong-greedoid.AB}.]
Let $k = \abs{A}$; thus, $\abs{B} = \abs{A} + 1 = k + 1$.
Let $\tup{a_1, a_2, \ldots, a_k}$ be a list of all $k$ elements of $A$
(with no repetitions).

We define $k$ elements $b_1, b_2, \ldots, b_k \in B$
recursively as follows:
For each $i = 1,2,\ldots,k$, we let $b_i$ be a projection
of $a_i$ onto $B \setminus \set{b_1,b_2,\ldots, b_{i-1}}$
(assuming that $b_1, b_2, \ldots, b_{i-1}$ have already
been constructed).%
\footnote{Thus, in particular, $b_1$ is a projection
of $a_1$ onto $B$.}
Thus, $b_1, b_2, \ldots, b_k$ are $k$ distinct\footnote{The
distinctness of $b_1, b_2, \ldots, b_k$ follows from
$b_i \in B \setminus \set{b_1,b_2,\ldots, b_{i-1}}$.}
elements of $B$.
Thus, $\set{b_1, b_2, \ldots, b_k}$ is a $k$-element
subset of $B$.
Hence, its complement $B \setminus \set{b_1, b_2, \ldots, b_k}$
has size $\abs{B} - k = 1$ (since $\abs{B} = k + 1$).
In other words,
there is a unique element
$u \in B \setminus \set{b_1, b_2, \ldots, b_k}$.
Consider this $u$.
Hence,
$B \setminus \set{b_1, b_2, \ldots, b_k} = \set{u}$,
so that $B \setminus u = \set{b_1, b_2, \ldots, b_k}$.
From $u \in B \setminus \set{b_1, b_2, \ldots, b_k}$,
we obtain $u \notin \set{b_1, b_2, \ldots, b_k}$.

We have $u \notin A$.

[\textit{Proof:} Assume the contrary.
Thus, $u \in A = \set{a_1, a_2, \ldots, a_k}$.
Hence, $u = a_i$ for some $i \in \set{1, 2, \ldots, k}$.
Consider this $i$.
But $a_i = u \in B \setminus \set{b_1,b_2,\ldots, b_k}
\subset B \setminus \set{b_1,b_2,\ldots, b_{i-1}}$.
Hence, Lemma~\ref{lem.projection} \textbf{(a)} (applied
to $B \setminus \set{b_1,b_2,\ldots, b_{i-1}}$, $a_i$
and $b_i$ instead of $C$, $v$ and $u$)
yields $b_i = a_i$
(because $b_i$ is defined as a projection of
$a_i$ onto $B \setminus \set{b_1,b_2,\ldots, b_{i-1}}$).
Hence, $u = a_i = b_i$, which contradicts
$u \notin \set{b_1, b_2, \ldots, b_k}$.
This contradiction shows that our assumption was false.
Hence, $u \notin A$ is proven.]

Combining $u \in B$ with $u \notin A$, we find
$u \in B \setminus A$.

For each $i \in \set{1, 2, \ldots, k}$, we have
\begin{align}
d \tup{a_i, u} \geq d \tup{b_i, u} .
\label{pf.prop.strong-greedoid.1}
\end{align}

[\textit{Proof:} Let $i \in \set{1, 2, \ldots, k}$.
Then, from $u \notin \set{b_1, b_2, \ldots, b_k}$,
we obtain $u \neq b_i$.
Also, $u \in B \setminus \set{b_1,b_2,\ldots, b_k}
\subset B \setminus \set{b_1,b_2,\ldots, b_{i-1}}$,
whereas $b_i$ is a projection of $a_i$ onto
$B \setminus \set{b_1,b_2,\ldots, b_{i-1}}$.
Hence, Lemma~\ref{lem.projection} \textbf{(c)} (applied to
$B \setminus \set{b_1,b_2,\ldots, b_{i-1}}$,
$a_i$, $b_i$ and $u$ instead of $C$, $v$, $u$ and $x$)
shows that $d \tup{b_i, u} \leq d \tup{a_i, u}$
(since $u \neq b_i$).
This proves \eqref{pf.prop.strong-greedoid.1}.]

We have
$B \setminus u = \set{b_1, b_2, \ldots, b_k}$
and thus
\begin{align}
\sum_{b \in B \setminus u} d\tup{b, u}
= \sum_{i=1}^k d\tup{b_i, u}
\label{pf.prop.strong-greedoid.2a}
\end{align}
(since $b_1, b_2, \ldots, b_k$ are distinct).

From $u \in B$, we obtain
\begin{align*}
\PER\tup{B}
&= \PER\tup{B \setminus u} + w\tup{u} + \sum_{b \in B \setminus u} d\tup{b, u} \\
&= \PER\tup{B \setminus u} + w\tup{u} + \sum_{i=1}^k d\tup{b_i, u}
\end{align*}
(by \eqref{pf.prop.strong-greedoid.2a}).
Solving this for $\PER\tup{B \setminus u}$, we obtain
\begin{align}
\PER\tup{B \setminus u}
= \PER\tup{B} - w\tup{u} - \sum_{i=1}^k d\tup{b_i, u} .
\label{pf.prop.strong-greedoid.PERBu}
\end{align}

We have $A = \set{a_1, a_2, \ldots, a_k}$ and thus
\begin{align}
\sum_{a \in A} d\tup{a, u} = \sum_{i=1}^k d\tup{a_i, u}
\label{pf.prop.strong-greedoid.2b}
\end{align}
(since $a_1, a_2, \ldots, a_k$ are distinct).

From $u \notin A$, we obtain
\begin{align*}
\PER\tup{A \cup u}
&= \PER\tup{A} + w\tup{u} + \sum_{a \in A} d\tup{a, u} \\
&= \PER\tup{A} + w\tup{u} + \sum_{i=1}^k \underbrace{d\tup{a_i, u}}_{\substack{\geq d\tup{b_i, u} \\ \tup{\text{by \eqref{pf.prop.strong-greedoid.1}}}}}
 \qquad \tup{\text{by \eqref{pf.prop.strong-greedoid.2b}}} \\
&\geq \PER\tup{A} + w\tup{u} + \sum_{i=1}^k d\tup{b_i, u} .
\end{align*}
Adding this inequality to the equality
\eqref{pf.prop.strong-greedoid.PERBu}, we obtain
\begin{align*}
\PER\tup{B \setminus u} + \PER\tup{A \cup u}
&\geq \PER\tup{B} + \PER\tup{A} \\
&= \PER\tup{A} + \PER\tup{B} .
\end{align*}
This is precisely the inequality \eqref{eq.lem.strong-greedoid.AB.per}.

Thus, we have found a $u \in B \setminus A$
satisfying \eqref{eq.lem.strong-greedoid.AB.per}.
Hence, such a $u$ exists.
This proves Lemma~\ref{lem.strong-greedoid.AB}.
\end{proof}

\begin{proof}[Proof of Theorem \ref{thm.greedoid}.]
We only need to prove the two axioms \textbf{(i)}
and \textbf{(ii)} from the definition of a greedoid
and the axiom \textbf{(iv)} from the definition of
a strong greedoid (because axiom \textbf{(iii)}
will follow from axiom \textbf{(iv)}).

Axiom \textbf{(i)} is obvious.

Next, let us prove axiom \textbf{(iv)}.
So let $A, B \in \calF$ be such that $\abs{B} = \abs{A} + 1$.
We must prove that there exists some $x \in B \setminus A$ such that
$A \cup x \in \calF$ and $B \setminus x \in \calF$.

Lemma~\ref{lem.strong-greedoid.AB} shows that
there exists a $u \in B \setminus A$
satisfying \eqref{eq.lem.strong-greedoid.AB.per}.
Consider this $u$.

Let $k = \abs{A}$; thus, $\abs{B} = \abs{A} + 1 = k + 1$.
But $u \in B \setminus A \subset B$, so that
$\abs{B \setminus u} = \abs{B} - 1 = k$
(since $\abs{B} = k+1$).
Thus, $B \setminus u$ is a $k$-set.
But $A$ is a $k$-set in $\calF$,
and thus has the largest perimeter among all
$k$-sets.
Hence, $\PER\tup{A} \geq \PER\tup{B \setminus u}$.

Furthermore, $u \in B \setminus A$, thus $u \notin A$, so that
$\abs{A \cup u} = \abs{A} + 1 = k + 1$.
Hence, $A \cup u$ is a $\tup{k+1}$-set.
But $B$ is a $\tup{k+1}$-set in $\calF$,
and thus has the largest perimeter among all $\tup{k+1}$-sets.
Hence, $\PER\tup{B} \geq \PER\tup{A \cup u}$.
Adding this inequality to
$\PER\tup{A} \geq \PER\tup{B \setminus u}$,
we obtain
\[
\PER\tup{A} + \PER\tup{B}
\geq \PER\tup{B \setminus u} + \PER\tup{A \cup u} .
\]
Contrasting this inequality with
the opposite inequality \eqref{eq.lem.strong-greedoid.AB.per}
(which, as we know, is satisfied),
we conclude that it must be an equality.
Hence, both inequalities
$\PER\tup{A} \geq \PER\tup{B \setminus u}$
and
$\PER\tup{B} \geq \PER\tup{A \cup u}$
(which we added to obtain it) must be equalities as well.
In other words,
$\PER\tup{A} = \PER\tup{B \setminus u}$
and
$\PER\tup{B} = \PER\tup{A \cup u}$.
Hence, $B \setminus u$ is a $k$-set of maximum perimeter
(since $A$ is a $k$-set of maximum perimeter, but
$\PER\tup{A} = \PER\tup{B \setminus u}$), and thus belongs
to $\calF$; in other words, $B \setminus u \in \calF$.
Likewise, from the other inequality, we obtain
$A \cup u \in \calF$.
Hence, there exists some $x \in B \setminus A$ such that
$A \cup x \in \calF$ and $B \setminus x \in \calF$
(namely, $x = u$).
Thus, axiom \textbf{(iv)} is proven.

Let us now prove axiom \textbf{(ii)}.
So let $B\in \calF$ satisfy $\abs{B}>0$.
Then, $\abs{B} - 1 \in \set{0, 1, \ldots, \abs{E}}$.
Hence, there exists at least one $\tup{\abs{B}-1}$-subset
of $E$.
Since $E$ is finite, we can thus find a $\tup{\abs{B}-1}$-subset
of $E$ having maximum perimeter
(among all $\tup{\abs{B}-1}$-subsets of $E$).
Choose such a subset, and denote it by $A$.
Thus, $A \in \calF$ (by the definition of $\calF$,
since $A$ has maximum perimeter)
and $\abs{B} = \abs{A} + 1$ (since $A$ is a
$\tup{\abs{B}-1}$-subset).
Hence, axiom \textbf{(iv)} (which we have already proved)
shows that there exists some $x \in B \setminus A$ such that
$A \cup x \in \calF$ and $B \setminus x \in \calF$.
Consider this $x$. Thus, $x \in B \setminus A \subset B$.
Hence, there exists $b \in B$ such that $B \setminus b \in \calF$
(namely, $b = x$).
This proves axiom \textbf{(ii)}.

This shows that $\calF$ is a strong greedoid.
\end{proof}

We now know that the Bhargava greedoid $\calF$
of an ultra triple is a strong greedoid.
It is natural to inquire which other known classes of
greedoids $\calF$ belongs to.
However, for many of these classes (including
interval greedoids), the answer is negative,
because $\calF$ is (in general) not a transposition
greedoid.
We refer to \cite[Chapter X]{KoLoSc91} for the
definition of transposition greedoids (and for why many
classes of greedoids are subclasses of transposition
greedoids); let us merely remark that the Bhargava
greedoid $\calF$ fails to be a transposition greedoid
in Example~\ref{exa.non-transp-greedoid}, since the
transposition property \cite[(1.1) in Section X.1]{KoLoSc91}
is violated for $A = \set{0, 1, 2}$, $x = 3$, $y = 6$
and $B = \set{4, 5, 12}$.

The Bhargava greedoid $\calF$ also fails to be a
transversal greedoid in the sense of
\cite{Brooks97}\footnote{Transversal greedoids are the
same as medieval marriage greedoids in the sense of
\cite[Section IV.2.14]{KoLoSc91}.}.
Indeed, the ultra triple
$\tup{E, w, d}$ constructed in Example~\ref{exa.p-adic.d}
for $p = 2$ and $E = \set{1, 2, 3, 4}$ provides
a counterexample\footnote{The easiest way to check
this is to observe that it violates the condition
(M3)$^\dagger$ from \cite[Theorem 2.1]{Brooks97}.
(Note that there is a typo in \cite[Theorem 2.1]{Brooks97}:
In Condition (M3)$^\dagger$, replace ``$Z \neq \varnothing$''
by ``$X \neq \varnothing$''.)}.

Another class of greedoids that the Bhargava
greedoid $\calF$ does not belong to is that
of twisted matroids (\cite[Section IV.2.18]{KoLoSc91}).
Indeed, \cite[Proposition 3.1]{Kloock03} shows
that every twisted matroid is a $\Delta$-matroid
(see \cite[Section 2.4]{Kloock03} for a definition of the
latter concept); but $\calF$ is not in general
a $\Delta$-matroid\footnote{For an example, use
the ultra triple $\tup{E, w, d}$ constructed in
Example~\ref{exa.p-adic.d}
for $p = 2$ and $E = \set{1, 2, 4, 8}$.
Here, the axiom defining a $\Delta$-matroid
fails for $X = \set{1, 2, 4, 8}$, $Y = \varnothing$
and $x = 1$.}.

In \cite{Grinbe20}, it is shown that the
Bhargava greedoid $\calF$ is
a Gaussian elimination greedoid
(see \cite[Sections IV.2 and IX.4]{KoLoSc91}
for this concept).

\begin{question}
Is $\calF$ a linking greedoid?
(This is yet another subclass of Gauss greedoids,
and can in some sense be understood as
``Gaussian elimination greedoids over the
field with one element'';
see again \cite[Sections IV.2 and IX.4]{KoLoSc91}.)
\end{question}

\section{The matroid}

Throughout this section, we assume that the set $E$ is finite.

\subsection{Defining matroids}

We shall now recall one of the many definitions of
a matroid.
Namely, if $E$ is a finite set, $k$ is a nonnegative
integer, and $\calB$ is a collection of
$k$-subsets of $E$, then we say that $\calB$
is \defn{the collection of bases of a matroid}
if and only if $\calB$ is nonempty and
satisfies the following axiom:\footnote{This
axiom is condition (1.4) in \cite[Section II.1]{KoLoSc91}.
See \cite[Theorem II.1.1]{KoLoSc91} for a proof of its
equivalence to other definitions of a matroid.
See also \cite{Oxley11} for much more about matroids.}

\begin{itemize}
    \item For any two $k$-subsets
    $B_1, B_2 \in \calB$ and any
    $x \in B_1 \setminus B_2$, there exists
    a $y \in B_2 \setminus B_1$ such that
    $B_1 \cup y \setminus x \in \calB$.
\end{itemize}

\subsection{Matroids from strong greedoids}

We now get to the main result of this section:

\begin{theorem}\label{thm.matroid}
The Bhargava greedoid $\calF$ has the following property:
Fix a nonnegative integer $k \leqslant \abs{E}$.
All sets $A \in \calF$ having size $k$
form the collection of bases of a matroid.
\end{theorem}

Not all greedoids enjoy this property.
For example, if $\set{a,b,c,d}$ is a poset
with two inequalities $a<b$ and $c<d$, then
the greedoid of lower ideals of this poset
contains the subsets $\set{a,b}$ and $\set{c,d}$,
but $a$ in the set
$\set{a,b}$ cannot be replaced by any of
$c$ and $d$.

However, all strong greedoids (i.e., Gauss greedoids)
enjoy this property:

\begin{theorem} \label{thm.s-g-matroid}
Let $\calF$ be a strong greedoid on the ground set $E$.
Let $B_1 \in \calF$ and $B_2 \in \calF$ satisfy
$\abs{B_1} = \abs{B_2}$.
Let $x \in B_1 \setminus B_2$.
Then, there exists some $y \in B_2 \setminus B_1$ such that
$B_1 \cup y \setminus x \in \calF$.
\end{theorem}

Theorem~\ref{thm.s-g-matroid} is (implicitly) proven
in the third paragraph of \cite[Proof of the Theorem]{BrySha99}.
For the sake of completeness, we shall present this proof
in a slightly modified form below.
First, we need two lemmas about greedoids:

\begin{lemma} \label{lem.greedoid-bigger}
Let $\calF$ be a greedoid on the ground set $E$.
Let $A, B \in \calF$ satisfy $\abs{B} > \abs{A}$.
Then, there exists some $b \in B \setminus A$ such that
$A \cup b \in \calF$.
\end{lemma}

\begin{proof}[Proof of Lemma~\ref{lem.greedoid-bigger}.]
A \defn{nice set} will mean a subset $C$ of $B$ such that
$\abs{C} > \abs{A}$ and $C \in \calF$.
There exists at least one nice set (namely, $B$ is a nice set).
Thus, there exists a nice set of smallest possible size.
Let $D$ be such a set.
Thus, $D$ is a subset of $B$ such that $\abs{D} > \abs{A}$
and $D \in \calF$ (since $D$ is a nice set).
Hence, $\abs{D} > \abs{A} \geq 0$.
Thus, axiom \textbf{(ii)} in the definition of a greedoid
(applied to $D$ instead of $B$) shows that there exists
a $b \in D$ such that $D \setminus b \in \calF$.
Pick such a $b$ and denote it by $d$.
Thus, $d \in D$ and $D \setminus d \in \calF$.
Note that $D \setminus d$ is a subset of $B$ (since $D$ is),
and has smaller size than $D$ (since $d \in D$).
Hence, if we had $\abs{D \setminus d} > \abs{A}$, then
$D \setminus d$ would be a nice set of smaller size than $D$;
but this would contradict the fact that $D$ is a nice set
of smallest possible size.
Thus, we must have $\abs{D \setminus d} \leq \abs{A}$.
Since $d \in D$, we have $\abs{D \setminus d} = \abs{D} - 1$,
so that $\abs{D} - 1 = \abs{D \setminus d} \leq \abs{A}$,
and therefore $\abs{D} \leq \abs{A} + 1$.
Combining this with $\abs{D} > \abs{A}$, we obtain
$\abs{D} = \abs{A} + 1$.
Hence, axiom \textbf{(iii)} in the definition of a greedoid
(applied to $D$ instead of $B$) shows that there exists
a $b \in D \setminus A$ such that $A \cup b \in \calF$.
Consider this $b$.
We have $b \in D \setminus A \subset B \setminus A$
(since $D \subset B$).
Thus, we have found a $b \in B \setminus A$ such that
$A \cup b \in \calF$.
This proves Lemma~\ref{lem.greedoid-bigger}.
\end{proof}

\begin{lemma} \label{lem.s-g-Dxyz}
Let $\calF$ be a strong greedoid on the ground set $E$.
Let $D$ be a subset of $E$, and let $x, y, z$ be three
elements of $E \setminus D$.
Assume that $D \cup \set{x, z} \in \calF$ and
$D \cup y \in \calF$ and $D \cup z \notin \calF$.
Then, we have $D \cup \set{y, z} \in \calF$.
\end{lemma}

\begin{proof}[Proof of Lemma~\ref{lem.s-g-Dxyz}.]
We have $D \cup \set{x, z} \neq D \cup z$
(since $D \cup \set{x, z} \in \calF$ but $D \cup z \notin \calF$).
Hence, $x \neq z$.
Furthermore,
none of the elements $x, y, z$ belongs to $D$ (since
they all belong to $E \setminus D$).
Hence, $\abs{D \cup \set{x, z}} = \abs{D \cup y} + 1$
(since $x \neq z$).
Consequently,
axiom \textbf{(iv)} in the definition of a strong greedoid
(applied to $A = D \cup y$ and $B = D \cup \set{x, z}$)
yields that there exists some
$t \in \tup{D \cup \set{x, z}} \setminus \tup{D \cup y}$
such that
$\tup{D \cup y} \cup t \in \calF$ and
$\tup{D \cup \set{x, z}} \setminus t \in \calF$.
Consider this $t$.

Combining $x \neq z$ with $x \notin D$, we obtain
$x \notin D \cup z$.
If we had $t = x$, then we would have
$\underbrace{\tup{D \cup \set{x, z}}}_{= \tup{D \cup z} \cup x} \setminus \underbrace{t}_{=x}
= \tup{D \cup z} \cup x \setminus x
= D \cup z$ (since $x \notin D \cup z$)
and therefore
$D \cup z
= \tup{D \cup \set{x, z}} \setminus t \in \calF$,
which would contradict $D \cup z \notin \calF$.
Hence, we must have $t \neq x$.

We have
$t \in \tup{D \cup \set{x, z}} \setminus \tup{D \cup y}
\subset \set{x, z}$,
so that either $t = x$ or $t = z$.
Thus, $t = z$ (since $t \neq x$).
Hence, $z = t$, so that
$D \cup \set{y, z}
= \tup{D \cup y} \cup \underbrace{z}_{=t}
= \tup{D \cup y} \cup t
\in \calF$.
\end{proof}

\begin{proof}[Proof of Theorem~\ref{thm.s-g-matroid}.]
From $x \in B_1 \setminus B_2$, we obtain $x \in B_1$ and
$x \notin B_2$.
Hence, $\abs{B_1 \setminus x} = \abs{B_1} - 1$.

A \defn{free set} will mean a subset $A$ of $B_1 \setminus x$
such that $A \in \calF$.
Clearly, a free set exists (indeed, $\varnothing$ is a free
set, since axiom \textbf{(i)} in the definition of a greedoid
yields $\varnothing \in \calF$).
Hence, there exists a free set of largest size.
Pick such a free set, and denote it by $A$.
Thus, $A$ is a subset of $B_1 \setminus x$ and satisfies
$A \in \calF$ (since $A$ is a free set).
Since $A$ is a subset of $B_1 \setminus x$, we have
$\abs{A} \leq \abs{B_1 \setminus x} = \abs{B_1} - 1
< \abs{B_1} = \abs{B_2}$.
Thus, Lemma~\ref{lem.greedoid-bigger} (applied to $B = B_2$)
yields that there exists some $b \in B_2 \setminus A$ such that
$A \cup b \in \calF$. Consider this $b$, and denote it by $y$.
Thus, $y \in B_2 \setminus A$ and $A \cup y \in \calF$.

Next, we claim that $A \cup x \in \calF$.

[\textit{Proof:} Assume the contrary.
Thus, $A \cup x \notin \calF$.
Recall that $\abs{A} < \abs{B_1}$.
Thus, Lemma~\ref{lem.greedoid-bigger} (applied to $B = B_1$)
yields that there exists some $b \in B_1 \setminus A$ such that
$A \cup b \in \calF$. Consider this $b$.
Clearly, $b \notin A$.
We have $A \cup b \neq A \cup x$
(since $A \cup b \in \calF$ but $A \cup x \notin \calF$),
and thus $b \neq x$.
Hence, $b \in B_1 \setminus x$ (since $b \in B_1 \setminus A \subset B_1$).
Clearly, the set $A \cup b$ has larger size than $A$
(since $b \notin A$).
Now, $A \cup b$ is a subset of $B_1 \setminus x$
(since $A \subset B_1 \setminus x$ and $b \in B_1 \setminus x$),
and thus is a free set (since $A \cup b \in \calF$)
of larger size than $A$.
This contradicts the fact that $A$ is a free set of largest size.
This contradiction shows that our assumption was wrong.
Hence, we have shown that $A \cup x \in \calF$.]

From $y \in B_2 \setminus A$, we obtain $y \in B_2$ and
$y \notin A$. Hence, the set $A \cup y$ has larger size
than $A$ (since $y \notin A$).
If we had $y \in B_1 \setminus x$, then $A \cup y$ would be
a subset of $B_1 \setminus x$ (since $A \subset B_1 \setminus x$),
and therefore $A \cup y$ would be a free set (since $A \cup y \in \calF$)
of larger size than $A$;
this would contradict the fact that $A$ is a free set
of largest size.
Hence, $y \notin B_1 \setminus x$.
Since $y \neq x$ (because $y \in B_2$
but $x \notin B_2$), we thus obtain $y \notin B_1$.
Hence, $y \in B_2 \setminus B_1$ (since $y \in B_2$).

Thus, if $B_1 \cup y \setminus x \in \calF$, then
Theorem~\ref{thm.s-g-matroid} is proven.
Hence, for the sake of contradiction, we assume that
$B_1 \cup y \setminus x \notin \calF$.

A \defn{useful set} will mean a set $C \subset E$
such that $A \subset C \subset B_1 \setminus x$
and $C \cup x \in \calF$ and $C \cup y \in \calF$.
The set $A$ is a useful set (since $A \subset A
\subset B_1 \setminus x$ and $A \cup x \in \calF$
and $A \cup y \in \calF$).
Hence, there exists a useful set.
Thus, there exists a useful set of maximum size.
Let $D$ be such a set.
Thus, $D$ is a useful set; that is, $D \subset E$
and $A \subset D \subset B_1 \setminus x$ and
$D \cup x \in \calF$ and $D \cup y \in \calF$.

We have $y \neq x$ and thus $\tup{B_1 \setminus x} \cup y
= B_1 \cup y \setminus x \notin \calF$.
Hence, the set $B_1 \setminus x$ is not a useful set.
Thus, $D \neq B_1 \setminus x$ (since $D$ is a useful
set).
Therefore, $D$ is a proper subset of $B_1 \setminus x$
(since $D \subset B_1 \setminus x$).
Hence, $\abs{D} < \abs{B_1 \setminus x} = \abs{B_1} - 1$
(since $x \in B_1$).
Thus, $\abs{D} + 1 < \abs{B_1}$, so that
$\abs{D \cup x} \leq \abs{D} + 1 < \abs{B_1}$.

The two sets $D \cup x$ and $B_1$ belong to $\calF$
and satisfy $\abs{D \cup x} < \abs{B_1}$.
Hence, Lemma~\ref{lem.greedoid-bigger}
(applied to $D \cup x$ and $B_1$ instead of $A$ and $B$)
yields that there exists some $b \in B_1 \setminus \tup{D \cup x}$
such that $D \cup x \cup b \in \calF$.
Consider this $b$, and denote it by $z$.
Thus, $z \in B_1 \setminus \tup{D \cup x}$ and
$D \cup x \cup z \in \calF$.
Hence, $D \cup \set{x, z} = D \cup x \cup z \in \calF$.
Furthermore, $x \notin D$ (since $D \subset B_1 \setminus x$)
and $y \notin D$ (since $y \notin B_1$ but $D \subset B_1 \setminus x \subset B_1$)
and $z \notin D$ (since $z \in B_1 \setminus \tup{D \cup x}$ and thus $z \notin D \cup x$, so that $z \notin D$).
Hence, all of $x, y, z$ are elements of $E \setminus D$.

The set $D \cup z$ has larger size than $D$
(since $z \notin D$), and thus has larger size
than $A$ (since $A \subset D$ entails $\abs{A} \leq \abs{D}$).
Combining $D \subset B_1 \setminus x$
and $z \in B_1 \setminus \tup{D \cup x} \subset B_1 \setminus x$,
we obtain
$D \cup z \subset B_1 \setminus x$.
Hence, if we had $D \cup z \in \calF$, then $D \cup z$
would be a free set of larger size than $A$.
This would contradict the fact that $A$ is a free
set of largest size.
Hence, $D \cup z \notin \calF$.
Thus, Lemma~\ref{lem.s-g-Dxyz} shows that
$D \cup \set{y, z} \in \calF$.
Now, the set $D \cup z$ has larger size than $D$
and satisfies $A \subset D \cup z \subset B_1 \setminus x$
(since $A \subset D \subset D \cup z$ and $D \cup z \subset B_1 \setminus x$)
and $\tup{D \cup z} \cup x = D \cup \set{x, z} \in \calF$
and $\tup{D \cup z} \cup y = D \cup \set{y, z} \in \calF$.
Hence, $D \cup z$ is a useful set of larger size than $D$.
This contradicts the fact that $D$ is a useful set of maximum
size.
This contradiction shows that our assumption (that
$B_1 \cup y \setminus x \notin \calF$) was wrong.
Hence, $B_1 \cup y \setminus x \in \calF$.
This proves Theorem~\ref{thm.s-g-matroid}.
\end{proof}

We note that the condition ``$\abs{B_1} = \abs{B_2}$''
in Theorem~\ref{thm.s-g-matroid} could be replaced by
the weaker condition ``$\abs{B_1} \leq \abs{B_2}$''.
Indeed, our proof of Theorem~\ref{thm.s-g-matroid} only
used the latter condition.

\begin{proof}[Proof of Theorem~\ref{thm.matroid}.]
The assumption $k \leq \abs{E}$ shows that there exist
$k$-sets.
Some of them have maximum perimeter (since $E$ is finite).
Hence, the collection
of all sets $A \in \calF$ having size $k$ is nonempty.

Theorem~\ref{thm.greedoid} shows that
$\calF$ is a strong greedoid.
Hence, Theorem~\ref{thm.s-g-matroid}
shows that for any two sets $B_1 \in \calF$ and
$B_2 \in \calF$ satisfying
$\abs{B_1} = \abs{B_2}$
and for any $x \in B_1 \setminus B_2$,
there exists some $y \in B_2 \setminus B_1$ such that
$B_1 \cup y \setminus x \in \calF$.
This yields that all sets $A \in \calF$ having size $k$
form the collection of bases of a matroid
(because if $B_1$ and $B_2$ are two sets of size $k$,
and if $x \in B_1 \setminus B_2$ and $y \in B_2 \setminus B_1$,
then $B_1 \cup y \setminus x$ is a set of size $k$
as well).
This is precisely the claim of Theorem~\ref{thm.matroid}.
\end{proof}

\begin{verlong}

\subsection{Old proof of Theorem~\ref{thm.matroid}}

Next, we will show a different proof of Theorem~\ref{thm.matroid},
which was written before we were aware of
Theorem~\ref{thm.greedoid}.
It is relatively laborious and only presented here
for the sake of completeness.

Before we show this proof, let us
milk our proofs of Theorem~\ref{thm.2a} and Theorem~\ref{thm.2b}
for some further consequences:

\begin{proposition} \label{prop.u1}
Let $C \subset E$ be any subset.
Let $m$ be a nonnegative integer.
Let $\tup{c_1, c_2, \ldots, c_m}$ be a greedy $m$-permutation of $C$.
Let $k \in \set{0, 1, \ldots, m}$.
Let $A$ be a $k$-subset of $C$ with maximum perimeter.
We define $k$ elements $v_1, v_2, \ldots, v_k \in A$
recursively as follows:
For each $i=1,2,\ldots,k$, we let $v_i$ be a projection
of $c_i$ onto $A\setminus \set{v_1,v_2,\ldots, v_{i-1}}$.

Then:

\begin{enumerate}

\item[\textbf{(a)}] The equalities \eqref{zamena}
and \eqref{zamena1} hold for all $j \in \set{1, 2, \ldots, k}$.

\item[\textbf{(b)}] We have
$\PER\set{c_1, c_2, \ldots, c_j} = \PER\set{v_1, v_2, \ldots, v_j}$
for each $j \in \set{0, 1, \ldots, k}$.

\item[\textbf{(c)}] We have
$\PER\set{c_1, c_2, \ldots, c_p, v_{p+1}, v_{p+2}, \ldots, v_k}
= \PER\tup{A}$
for each $p \in \set{0, 1, \ldots, k}$.

\item[\textbf{(d)}] If $k \geq 1$ and $v_1 \neq c_1$,
then $\PER\tup{A \cup c_1 \setminus v_1} = \PER\tup{A}$.

\end{enumerate}

\end{proposition}

\begin{proof}[Proof of Proposition~\ref{prop.u1}.]
Our notations are consistent with the notations used
in the proof of Theorem~\ref{thm.2a} and also
with those used in the proof of Theorem~\ref{thm.2b}.
Hence, the equalities \eqref{zamena},
\eqref{zamena1} and \eqref{zamena2anon} can be shown as in the
proofs of these two theorems.
In particular, we thus see that
the equalities \eqref{zamena}
and \eqref{zamena1} hold for all $j \in \set{1, 2, \ldots, k}$.
Hence, Proposition~\ref{prop.u1} \textbf{(a)} is proven.

\textbf{(b)} Let $j \in \set{0, 1, \ldots, k}$.
As we just observed, \eqref{zamena2anon} holds.
In other words,
$\PER\set{c_1, c_2, \ldots, c_j} = \PER\set{v_1, v_2, \ldots, v_j}$.
This proves Proposition~\ref{prop.u1} \textbf{(b)}.

Furthermore, $A = \set{v_1, v_2, \ldots, v_k}$ (as we have
seen in the proof of Theorem~\ref{thm.2a}).
The points
$v_1, v_2, \ldots, v_k$ are distinct (by construction);
the points
$c_1, c_2, \ldots, c_m$ are distinct
(since $\tup{c_1, c_2, \ldots, c_m}$ is a greedy $m$-permutation).

\textbf{(c)}
Let $p \in \set{0, 1, \ldots, k}$.

Recall that \eqref{zamena} holds for all $j \in \set{1, 2, \ldots, k}$.
Hence, for each $i \in \set{1, 2, \ldots, p}$
and each $j \in \set{p+1, p+2, \ldots, k}$, we have
\begin{align}
d\tup{v_i, v_j} = d\tup{c_i, v_j}
\label{pf.prop.u1.c.i<j}
\end{align}
(since $j \in \set{p+1, p+2, \ldots, k}$ entails
$j \geq p+1$, thus $j-1 \geq p$,
but $i \in \set{1, 2, \ldots, p}$ entails $i \leq p$,
and therefore $i \leq p \leq j-1$, so that
$i \in \set{1, 2, \ldots, j-1}$).

Each $j \in \set{1, 2, \ldots, k}$ satisfies
$v_j \notin \set{c_1, c_2, \ldots, c_{j-1}}$
(as we showed in the proof of Theorem~\ref{thm.2a}\footnote{Indeed, this is precisely \eqref{pf.thm.2a.notin}.});
thus, the points $v_{p+1}, v_{p+2}, \ldots, v_k$
are distinct from the points $c_1, c_2, \ldots, c_p$.
Hence, $c_1, c_2, \ldots, c_p, v_{p+1}, v_{p+2}, \ldots, v_k$
are altogether $k$ distinct points
(since we know that the points
$v_1, v_2, \ldots, v_k$ are distinct, and that
the points
$c_1, c_2, \ldots, c_m$ are distinct).
Thus,
\begin{align*}
& \PER\set{c_1, c_2, \ldots, c_p, v_{p+1}, v_{p+2}, \ldots, v_k} \\
&= \underbrace{\PER\set{c_1, c_2, \ldots, c_p}}_{\substack{= \PER\set{v_1, v_2, \ldots, v_p}\\ \text{(by Proposition~\ref{prop.u1} \textbf{(b)},} \\ \text{ applied to $j = p$)}}}
    + \sum_{i=1}^p \sum_{j=p+1}^k \underbrace{d(c_i, v_j)}_{\substack{= d(v_i, v_j) \\ \text{(by \eqref{pf.prop.u1.c.i<j})}}}
    + \PER\set{v_{p+1}, v_{p+2}, \ldots, v_k} \\
&= \PER\set{v_1, v_2, \ldots, v_p} + \sum_{i=1}^p \sum_{j=p+1}^k d(v_i, v_j) + \PER\set{v_{p+1}, v_{p+2}, \ldots, v_k} \\
&= \PER\set{v_1, v_2, \ldots, v_k} \qquad \tup{\text{since the elements $v_1, v_2, \ldots, v_k$ are distinct}} \\
&= \PER\tup{A}
\end{align*}
(since $A = \set{v_1, v_2, \ldots, v_k}$).
This proves Proposition~\ref{prop.u1} \textbf{(c)}.

\textbf{(d)} Assume that $k \geq 1$ and $v_1 \neq c_1$.
Then, from $A = \set{v_1, v_2, \ldots, v_k}$, we obtain
$A \cup c_1 \setminus v_1
= \set{c_1, v_2, v_3, \ldots, v_k}$
(since $v_1 \neq c_1$, and since $v_1, v_2, \ldots, v_k$
are distinct).
But Proposition~\ref{prop.u1} \textbf{(c)} (applied to
$p = 1$) yields
$\PER\set{c_1, v_2, v_3, \ldots, v_k}
= \PER\tup{A}$.
In view of
$A \cup c_1 \setminus v_1
= \set{c_1, v_2, v_3, \ldots, v_k}$,
this rewrites as $\PER\tup{A \cup c_1 \setminus v_1} = \PER\tup{A}$.
This proves Proposition~\ref{prop.u1} \textbf{(d)}.
\end{proof}

For the next corollary (which will be used in our
proof of Theorem \ref{thm.matroid}), let us define a simple
piece of notation:
If $C$ is any subset of $E$, then a \defn{starter} of $C$
shall mean an element $c \in C$ of maximum weight (among
the elements of $C$).
Note that any nonempty subset $C$ of $E$ has at least
one starter (since $C$ is finite).

\begin{corollary} \label{cor.starters}
Let $C \subset E$ be an $m$-set.
Let $k \in \set{1, 2, \ldots, m}$.
Let $A$ be a $k$-subset of $C$ with maximum perimeter.
Then:

\begin{enumerate}

\item[\textbf{(a)}] Every starter of $A$ is also a starter
of $C$.

\item[\textbf{(b)}] Let $u$ be a starter of $C$.
Let $v$ be a projection of $u$ onto $A$.
Then, $v$ is a starter of $A$
and satisfies $w(v) = w(u)$.

\item[\textbf{(c)}] With the notations of
part \textbf{(b)}, we have
$\PER\tup{A \cup u \setminus v} = \PER\tup{A}$
if $u \notin A$.

\end{enumerate}

\end{corollary}

\begin{proof}[Proof of Corollary~\ref{cor.starters}.]

\textbf{(b)} The $1$-subset $\set{u}$ of $C$ has maximum
perimeter among all $1$-subsets of $C$ (since $u$ is a starter,
i.e., has maximum weight).
Thus, Theorem \ref{thm.2b} (applied to $1$ instead
of $k$) shows that there exists some greedy $m$-permutation
$\tup{c_1, c_2, \ldots, c_m}$ of $C$ such that
$\set{u} = \set{c_1}$.
Consider this greedy $m$-permutation.
From $\set{u} = \set{c_1}$, we obtain $c_1 = u$.
Hence, $v$ is a projection of $c_1$ onto $A$
(since $v$ is a projection of $u$ onto $A$).

Define $k$ elements $v_1, v_2, \ldots, v_k \in A$
as in Proposition~\ref{prop.u1}, making sure to pick
$v_1$ to be $v$.
(Indeed, the definition of these $k$ elements sets
$v_1$ to be any projection of $c_1$ onto $A$;
we use this freedom to set $v_1 = v$, since $v$ is a
projection of $c_1$ onto $A$.)

Proposition \ref{prop.u1} \textbf{(a)} shows that
the equalities \eqref{zamena}
and \eqref{zamena1} hold for all $j \in \set{1, 2, \ldots, k}$.
Thus, in particular, \eqref{zamena1} holds.
Applying \eqref{zamena1} to $j = 1$, we find
$w(v_1) = w(c_1)$ (since the sums on both sides
are empty).
In view of $v_1 = v$ and $c_1 = u$, this becomes
$w(v) = w(u)$.
But $u$ is a starter of $C$; thus, $w(u) \geq w(p)$ for
each $p \in C$.
Hence, $w(u) \geq w(p)$ for each $p \in A$
(since $A \subset C$).
Thus, $w(v) = w(u) \geq w(p)$ for each $p \in A$.
Thus, $v$ is a starter of $A$ (since $v \in A$).
This completes the proof of
Corollary \ref{cor.starters} \textbf{(b)}.

\textbf{(c)}
Continue where our above proof of
Corollary \ref{cor.starters} \textbf{(b)} left off.
Assume that $u \notin A$. Hence, $v \neq u$ (since $v \in A$).
Thus, $v_1 = v \neq u = c_1$. Therefore,
Proposition \ref{prop.u1} \textbf{(d)} shows that
$\PER\tup{A \cup c_1 \setminus v_1} = \PER\tup{A}$.
In view of $v_1 = v$ and $c_1 = u$, this rewrites as
$\PER(A \cup u \setminus v) = \PER\tup{A}$.
This proves Corollary \ref{cor.starters} \textbf{(c)}.

\textbf{(a)}
It clearly suffices to show that the maximum weight
of an element of $A$ equals the maximum weight of
an element of $C$.

To that aim, we pick any starter $u$ of $C$.
(Such a $u$ exists, since $C$ is nonempty.)
Choose any projection $v$ of $u$ onto $A$.
Then, Corollary \ref{cor.starters} \textbf{(b)}
shows that $v$ is a starter of $A$
and satisfies $w(v) = w(u)$.
Now, $w(v)$ is the maximum weight of an element of
$A$ (since $v$ is a starter of $A$),
while $w(u)$ is the maximum weight of an element
of $C$ (since $u$ is a starter of $C$).
Thus, these two maximum weights are equal
(since $w(v) = w(u)$).
This proves Corollary \ref{cor.starters} \textbf{(a)}.
\end{proof}

\begin{proof}[Alternative proof of Theorem \ref{thm.matroid}.]
Strong induction on $k$.
Fix a nonnegative integer $k$.
Assume that for all smaller values of $k$,
Theorem \ref{thm.matroid} is proved (not only for $(E,w,d)$ but
for any ultra triple).
Let $B_1,B_2$ be two $k$-subsets of $E$ having
maximum perimeter, and let $x\in B_1\setminus B_2$.
We should find an element
$y\in B_2\setminus B_1$ such that
\begin{equation}\label{equal-perim}
\PER(B_1\cup y\setminus x)=\PER(B_1) .
\end{equation}

If $A:=B_1\cap B_2\ne \varnothing $,
then we can consider the following ultra triple
on the set $E\setminus A$:
The distance is still $d$, and the weight of
$v\in E\setminus A$ equals
$w(v)+\sum_{a\in A} d(v,a)$. The sets
$B_2\setminus B_1$ and $B_1\setminus B_2$
have maximum perimeters among
$(k-\abs{A})$-subsets of $E\setminus A$
(with respect to this new ultra triple);
thus, by the induction hypothesis,
there exists a $y\in B_2\setminus B_1$
satisfying \eqref{equal-perim}.

Hence,
\begin{align}
\text{Theorem \ref{thm.matroid} is proved for our $k$ whenever $B_1 \cap B_2 \neq \varnothing$.}
\label{pf.matroid.intersect}
\end{align}

It remains to consider the case
$B_1\cap B_2=\varnothing$.

Note that $k \geq 1$ (due to the existence
of $x \in B_1 \setminus B_2$).

Consider two cases.

\begin{itemize}

\item \textit{Case 1:}
There exists a starter $z$ of $B_1$ distinct from $x$.

Consider such a $z$, and let
$v$ be a projection of $z$ onto $B_2$.
Note that $z \in B_1$, thus $z \notin B_2$
(since $B_1 \cap B_2 = \varnothing$).
Also, Corollary \ref{cor.starters} \textbf{(a)}
(applied to $m = 2k$, $C = B_1 \cup B_2$ and $A = B_1$)
yields that $z$ is a starter of $B_1 \cup B_2$ (since
$z$ is a starter of $B_1$).
Hence, Corollary \ref{cor.starters} \textbf{(c)}
(applied to $m = 2k$, $C = B_1 \cup B_2$,
$A = B_2$ and $u = z$) yields that
$\PER\tup{B_2 \cup z \setminus v} = \PER \tup{B_2}$
(since $z \notin B_2$).
Since $B_2 \cup z \setminus v$ is a $k$-subset
of $E$ (since $z \notin B_2$ and $v \in B_2$),
we thus conclude that the $k$-set $B_2 \cup z \setminus v$
has maximum perimeter.
But the maximum-perimeter $k$-sets $B_1$ and
$B_3:=B_2\cup z\setminus v$
have a common element $z$; thus, they satisfy
$B_1 \cap B_3 \neq \varnothing$.
Hence, \eqref{pf.matroid.intersect} shows that
Theorem \ref{thm.matroid} holds for $B_3$ instead of $B_2$.
Thus, there exists some
$y\in B_3\setminus B_1$
such that that \eqref{equal-perim} holds.
This $y$ must also belong to $B_2\setminus B_1$
(since $B_3 \setminus B_1 \subset B_2\setminus B_1$),
and thus completes our induction in Case 1.

\item \textit{Case 2:}
The only starter of $B_1$ is $x$.

Corollary \ref{cor.starters} \textbf{(a)}
(applied to $m = 2k$, $C = B_1 \cup B_2$ and $A = B_1$)
yields that $x$ is a starter of $B_1 \cup B_2$ (since
$x$ is a starter of $B_1$).

Let $p$ be a projection of $x$ onto $B_2$.
Then, Corollary \ref{cor.starters} \textbf{(b)}
(applied to $m = 2k$, $C = B_1 \cup B_2$,
$A = B_2$, $u = x$ and $v = p$) yields that
$p$ is a starter of $B_2$ and satisfies $w(p) = w(x)$.

Let $v$ be a projection of $p$ onto $B_1$.
Then, Corollary \ref{cor.starters} \textbf{(b)}
(applied to $m = 2k$, $C = B_1 \cup B_2$,
$A = B_1$ and $u = p$) yields that
$v$ is a starter of $B_1$ and satisfies $w(v) = w(p)$.
Thus, $v = x$ (since the only starter of $B_1$ is $x$).
But $p \in B_2$, thus $p \notin B_1$ (since $B_1 \cap B_2 = \varnothing$).
Hence, Corollary \ref{cor.starters} \textbf{(b)}
(applied to $m = 2k$, $C = B_1 \cup B_2$,
$A = B_1$ and $u = p$) yields that
$\PER\tup{B_1 \cup p \setminus v} = \PER\tup{B_1}$.
In view of $v = x$, this rewrites as
$\PER\tup{B_1 \cup p \setminus x} = \PER\tup{B_1}$.
Note also that $p \in B_2 = B_2\setminus B_1$ (again
because $B_1 \cap B_2 = \varnothing$).
Hence, we have found a $y\in B_2\setminus B_1$ such that
\eqref{equal-perim} holds
(namely, $y = p$).
Thus, the induction is complete in Case 2.

\end{itemize}

We have thus completed our induction in each of the
two cases. Hence, Theorem~\ref{thm.matroid} is proven.
\end{proof}
\end{verlong}




\section{\label{sect.multi}Greedy subsequences}

We shall now study a slight variation of the notion
of greedy $m$-permutations, in which we allow picking the same
point multiple times.
This requires us to consider distances of the form
$d\tup{a, a}$, which our definition of ultra triple
does not support.
Thus, we begin by introducing a somewhat stronger
concept, that of ``full ultra triples''.

\subsection{Full ultra triples}

Consider again a set $E$.

As before, we assume that a function $w : E \to \RR$
is given, which assigns a \defn{weight} $w\tup{a}$
to each point $a \in E$.

Assume further that we are given a function
$d : E \times E \to \RR$, which we will call
the \defn{distance function}.
Thus, any two points $a,b\in E$
have a real-valued \defn{distance} $d\tup{a,b}$.
We assume that this distance function has the following properties:

\begin{itemize}

\item It is \defn{symmetric}: that is, $d\tup{a,b} = d\tup{b,a}$
      for any $a, b \in E$.

\item It satisfies the following inequality:
        \begin{equation}\label{ultra-full}
        d(a,b)\leqslant \max \set{d(a,c),d(b,c)}
        \end{equation}
      for any $a, b, c \in E$.
\end{itemize}

\noindent
(Again, \eqref{ultra-full} is just the
ultrametric triangle inequality; but keep in mind
that $d\tup{a, a}$ can be nonzero and even negative,
unlike in a metric space.)

Such a structure $(E,w,d)$ will be called a
\defn{full ultra triple}.
Thus, the notion of a full ultra triple differs
from that of an ultra triple in that the distance
function $d$ is defined on $E \times E$ rather than
on $E \timesu E$ (so that the distances
$d\tup{a, b}$ are defined for $a = b$ as well).

\begin{example}
Consider the situation of Example~\ref{exa.phylo},
but now define a map $d : E \times E \to \RR$ by
the same formula that was used to define the map
$d : E \timesu E \to \RR$ in Example~\ref{exa.phylo}.
Then, $\tup{E, w, d}$ is a full ultra triple.
\end{example}

It is immediately clear that if $\tup{E, w, d}$
is a full ultra triple, then $\tup{E, w, \underline{d}}$ is
an ultra triple, where $\underline{d} : E \timesu E \to \RR$
is the restriction of $d$ to the subset $E \timesu E$
of $E \times E$.
In other words, any full ultra triple $\tup{E, w, d}$
becomes an ultra triple if we restrict the distance
function $d$ to $E \timesu E$
(that is, if we forget the distances
$d\tup{a, a}$ between each point and itself).
Thus, any concept that was defined for ultra triples
(e.g., the concept of a greedy $m$-permutation)
is automatically defined for any full ultra triple
$\tup{E, w, d}$ as well
(just apply it to $\tup{E, w, \underline{d}}$),
and any proposition that has been proven for all ultra
triples can be applied to all full ultra triples.

Conversely, we can often -- but not always -- transform
an ultra triple into a full ultra triple as follows:

\begin{remark}
\label{rmk.multi.to-full}
Let $\tup{E, w, d}$ be an ultra triple.

Fix an $N \in \RR$ with the property that
\begin{align}
N \leq d\tup{a, b} \qquad \text{ for all } \tup{a, b} \in E \timesu E .
\label{eq.rmk.multi.to-full.ass}
\end{align}
(Such an $N$ always exists when $E$ is finite.)

Define a map $\overline{d} : E \times E \to \RR$
by setting
\[
\overline{d} \tup{a, b}
= \begin{cases}
    d\tup{a, b}, & \text{ if } a \neq b ; \\
    N          , & \text{ if } a = b
  \end{cases}
\qquad \text{for all } a, b \in E .
\]
Then, $\tup{E, w, \overline{d}}$ is a full ultra triple.
This full ultra triple extends the original ultra triple
$\tup{E, w, d}$
(in the sense that the distance function $d$ of the
latter is a restriction of $\overline{d}$).
\end{remark}

\begin{proof}[Proof of Remark~\ref{rmk.multi.to-full} (sketched).]
We only need to check that $\tup{E, w, \overline{d}}$ is a full ultra triple
(since the claim that $d$ is a restriction of $\overline{d}$
is obvious).
To do so, we must prove that $\overline{d}$ is symmetric, and that
it satisfies
\begin{align}
\overline{d}(a,b)\leqslant \max \set{\overline{d}(a,c),\overline{d}(b,c)} .
\label{pf.rmk.multi.to-full.1}
\end{align}

The symmetry of $\overline{d}$ follows trivially from the symmetry of $d$.
Thus it remains to prove \eqref{pf.rmk.multi.to-full.1}.
To do so, we distinguish between three cases:
\begin{itemize}
\item If $a = c$ or $b = c$, then \eqref{pf.rmk.multi.to-full.1} is trivial (since $\overline{d}(c, b) = \overline{d}(b, c)$).
\item If $a = b$, then \eqref{pf.rmk.multi.to-full.1} follows from \eqref{eq.rmk.multi.to-full.ass}.
\item Otherwise, $a, b, c$ are distinct, and thus \eqref{pf.rmk.multi.to-full.1} follows immediately from \eqref{ultra}.
\end{itemize}
\end{proof}

\begin{example}
For this example, we fix a prime number $p$
and a subset $E$ of $\ZZ$.
Define an ultra triple $(E, w, d)$ as in Example~\ref{exa.p-adic.d}.

Then, $0 \leq d\tup{a, b}$ for all $\tup{a, b} \in E \timesu E$.
Hence, we can define a map $\overline{d} : E \times E \to \RR$
as in Remark~\ref{rmk.multi.to-full} (by setting $N = 0$),
and obtain a full ultra triple $\tup{E, w, \overline{d}}$
that extends our ultra triple $(E, w, d)$.
\end{example}

\begin{example}
Let us see an example where the construction in
Remark~\ref{rmk.multi.to-full} does not work.

For this example, we fix a prime number $p$
and a subset $E$ of $\ZZ$.
Define an ultra triple $(E, w, d')$ as in Example~\ref{exa.p-adic.d'}.

If $E$ is infinite, then
there exists no $N \in \RR$ with the property that
$N \leq d'\tup{a, b}$ for all $\tup{a, b} \in E \timesu E$.
(Indeed, for each $m \in \NN$, there exist two distinct
elements $a$ and $b$ of $E$ satisfying $a \equiv b \mod p^m$
and therefore $d'\tup{a, b} \leq -m$.)
Hence, we cannot define $\overline{d}$ as in
Remark~\ref{rmk.multi.to-full}.%
\footnote{This does not mean that
the ultra triple $(E, w, d')$ cannot be obtained
by restricting a full ultra triple.
Sometimes it can (for example, when $E = \set{p^0, p^1, p^2, \ldots}$);
sometimes it cannot (for example, when $E = \ZZ$).}


It is tempting to try fixing this issue by setting
$d' \tup{a, a} = - \infty$ for all $a \in E$.
However, this would require a generalization of
the notion of a full ultra triple, allowing distances to
be $- \infty$;
this, in turn, would cause some complications in our proofs%
\footnote{In our proofs, we used the fact that if
a sum of finitely many inequalities between real numbers
is an equality, then each of the inequalities being summed
must itself be an equality.
(In other words: If
$\tup{a_i}_{i \in I}$ and $\tup{b_i}_{i \in I}$ are two
finite families of reals satisfying $a_i \geq b_i$ for all
$i \in I$ and $\sum_{i \in I} a_i = \sum_{i \in I} b_i$,
then $a_i = b_i$ for all $i \in I$.)
This is no longer true if we allow $- \infty$ as an addend.}.
Thus we are not making this generalization.
\end{example}

Note that every full ultra triple $\tup{E, w, d}$
satisfies
\begin{equation} \label{aa<ac}
d(a,a)\leqslant d(a,c)\quad \text{ for all } a, c \in E.
\end{equation}
In fact, this follows by substituting $b=a$ in \eqref{ultra-full}.

\subsection{Further definitions}

For the rest of Section~\ref{sect.multi}, we shall
fix a full ultra triple $\tup{E, w, d}$.

If $B \subset E$ and if $m$ is a nonnegative integer,
then an \defn{$m$-subsequence of $B$} shall mean an $m$-tuple
of elements of $B$ (not necessarily distinct).

If $\mathbf{a} = \tup{a_1, a_2, \ldots, a_m} \in E^m$ is any
$m$-tuple, then we define its \defn{perimeter $\PER\tup{\mathbf{a}}$} as
\[
\PER\tup{\mathbf{a}} := \sum_{k=1}^m w\tup{a_k} + \sum_{1 \leq i < j \leq m} d\tup{a_i, a_j} .
\]
This generalizes the perimeter of an $m$-set;
in fact, if the entries of the $m$-tuple
$\mathbf{a} = \tup{a_1, a_2, \ldots, a_m} \in E^m$ are distinct,
then
\begin{align}
\PER\tup{\mathbf{a}} = \PER\set{a_1, a_2, \ldots, a_m} .
\label{eq.per.tup=set}
\end{align}

If an $m$-tuple $\mathbf{a} \in E^m$ is a permutation of
an $m$-tuple $\mathbf{b} \in E^m$, then
$\PER\tup{\mathbf{a}} = \PER\tup{\mathbf{b}}$.
(This follows from the requirement $d\tup{a,b} = d\tup{b,a}$
on our distance function.)

\begin{definition}
Let $C \subset E$ be any subset, and let $m$ be a nonnegative
integer.

A \defn{greedy $m$-subsequence} of $C$ is an
$m$-subsequence $\tup{c_1, c_2, \ldots, c_m}$ of $C$
such that for each $i \in \set{1,2,\ldots,m}$ and each $x \in C$, we have
\begin{align}
\PER\tup{c_1, c_2, \ldots, c_i} \geqslant \PER\tup{c_1, c_2, \ldots, c_{i-1}, x} .
\label{eq.def.greedy-seq.geq}
\end{align}

\end{definition}

Thus, this notion differs from the notion of a
greedy $m$-permutation in two aspects:
A greedy $m$-subsequence is allowed to have equal entries,
and the inequality \eqref{eq.def.greedy-seq.geq} is
required to hold for all $x \in C$ (rather than only for
$x \in C \setminus \set{c_1, c_2, \ldots, c_{i-1}}$).
Thus, greedy $m$-subsequences are like greedy $m$-permutations
except that we are sampling with replacement.

\subsection{Main analogues}

We can now state the following analogues of
Theorem~\ref{thm.2a}, Theorem~\ref{thm.2b},
Corollary~\ref{cor.bh-t1} and Corollary~\ref{cor.bh-l2}, respectively:

\begin{theorem}\label{thm.multi-2a}

Let $C\subset E$ be any subset,
and let $m$ be a nonnegative integer.

Let $\tup{c_1, c_2, \ldots, c_m}$ be any greedy $m$-subsequence of $C$.

Then, for each $k \in \set{0,1,\dots,m}$, the $k$-subsequence
$\tup{c_1, c_2, \ldots, c_k}$ has maximum perimeter among
all $k$-subsequences of $C$.

\end{theorem}

\begin{theorem}\label{thm.multi-2b}

Let $C\subset E$ be any finite nonempty subset,
and let $m$ be a nonnegative integer.
Let $k \in \set{0,1,\dots,m}$.

Let $\mathbf{a}$ be any $k$-subsequence of $C$ with
maximum perimeter.
Then, there exists a greedy $m$-subsequence $\tup{c_1, c_2, \ldots, c_m}$ of $C$
such that $\mathbf{a}$ is a permutation of
the $k$-tuple $\tup{c_1, c_2, \ldots, c_k}$.

\end{theorem}

\begin{corollary} \label{cor.multi-bh-t1}
Let $C \subset E$ be any subset.
Let $m$ be a nonnegative integer.
Let $k \in \set{1, 2, \ldots, m}$.
If $\tup{c_1, c_2, \ldots, c_m}$ is a greedy
$m$-subsequence of $C$, then the number
\[
w\tup{c_k} + \sum_{i=1}^{k-1} d\tup{c_i, c_k}
\]
does not depend on the choice of this $m$-subsequence
(but only depends on $k$ and on $C$).
\end{corollary}

From now on, the number
\[
w\tup{c_k} + \sum_{i=1}^{k-1} d\tup{c_i, c_k}
\]
in Corollary~\ref{cor.multi-bh-t1} will be denoted by
\defn{$\nu_k\tup{C}$}.

\begin{corollary} \label{cor.multi-bh-l2}
Let $C \subset E$ be any subset.
Let $m$ be a nonnegative integer.
Let $k \in \set{1, 2, \ldots, m}$.
Let $\tup{c_1, c_2, \ldots, c_m}$ be a greedy
$m$-subsequence of $C$.
Let $j \in \set{1, 2, \ldots, k}$.
Then,
\begin{equation}
\nu_k\tup{C}
\leq w\tup{c_j} + \sum_{i \in \set{1, 2, \ldots, k} \setminus \set{j}} d\tup{c_i, c_j} .
\label{eq.cor.multi-bh-l2.claim}
\end{equation}
\end{corollary}

Note that Corollary~\ref{cor.multi-bh-t1}
(in the particular case when $w\tup{e} = 0$ for all $e \in E$)
is \cite[Conjecture 1]{Grinbe19}, while
Theorem~\ref{thm.multi-2a}
(in the same particular case)
is \cite[Conjecture 2]{Grinbe19}.

\subsection{\label{subsect.multi.clone}The clone construction}

We shall prove Theorem~\ref{thm.multi-2a},
Theorem~\ref{thm.multi-2b},
Corollary~\ref{cor.multi-bh-t1} and
Corollary~\ref{cor.multi-bh-l2} by deriving them
from the corresponding facts we have already proven
about greedy $m$-permutations and maximum-perimeter
subsets.
This derivation will rely on constructing a larger
full ultra triple $\tup{\Eh, \wh, \dh}$ whose ground set
$\Eh$ will contain a sufficiently large number of
``clones'' of each element of $E$.
These ``clones'' will allow us to transform any
$m$-tuple of elements of $E$ into an $m$-tuple of
\textbf{distinct} elements of $\Eh$ without
disturbing properties like greediness and perimeter.

We construct the new full ultra triple $\tup{\Eh, \wh, \dh}$
as follows:

\begin{itemize}

\item We fix a positive integer $N$.
      (For now, $N$ can be arbitrary, but later $N$
      will be assumed large enough.)

\item We let \defn{$\IN$} be the set $\set{1, 2, \ldots, N}$.

\item We define \defn{$\Eh$} to be the set
      $E \times \IN$.
      It consists of all pairs $\tup{e, i}$ with
      $e \in E$ and $i \in \IN$.

\item We define a function \defn{$\wh : \Eh \to \RR$} by setting
      \[
      \wh \tup{e, i} = w \tup{e}
      \qquad \text{ for each } \tup{e, i} \in \Eh .
      \]

\item We define a function \defn{$\dh : \Eh \times \Eh \to \RR$} by setting
      \[
      \dh \tup{\tup{e, i}, \tup{f, j}} = d \tup{e, f}
      \qquad \text{ for each } \tup{e, i} \in \Eh \text{ and } \tup{f, j} \in \Eh .
      \]

\end{itemize}

It is easy to see that $\tup{\Eh, \wh, \dh}$ is again a
full ultra triple.
The intuitive meaning of the construction of this full
ultra triple is that we have replaced each element $e$
of $E$ by $N$ ``clones'' $\tup{e, 1}, \tup{e, 2}, \ldots, \tup{e, N} \in \Eh$;
the weights and the mutual distances of these clones
are copied over from their originals in $E$.
From this point of view, the following lemma should not be surprising:

\begin{lemma} \label{lem.Eh.c}
Let $C$ be a subset of $E$.
Let $\Ch$ be the subset $C \times \IN$ of $\Eh$.
Let $m$ be a nonnegative integer.
Let $c_1, c_2, \ldots, c_m$ be any $m$ elements of $C$.
Let $r_1, r_2, \ldots, r_m$ be any $m$ elements of $\IN$.
Then:

\begin{enumerate}

\item[\textbf{(a)}]
We have
\[
\PER\tup{\tup{c_1, r_1}, \tup{c_2, r_2}, \ldots, \tup{c_m, r_m}}
= \PER\tup{c_1, c_2, \ldots, c_m} .
\]
(Here, the perimeter on the left hand side is computed
with respect to the full ultra triple $\tup{\Eh, \wh, \dh}$,
whereas that on the right hand side is computed with
respect to the full ultra triple $\tup{E, w, d}$.)

\end{enumerate}

From now on, assume that $r_1, r_2, \ldots, r_m$
are distinct.

\begin{enumerate}

\item[\textbf{(b)}]
We have
\[
\PER\set{\tup{c_1, r_1}, \tup{c_2, r_2}, \ldots, \tup{c_m, r_m}}
= \PER\tup{c_1, c_2, \ldots, c_m} .
\]

\item[\textbf{(c)}]
The $m$-tuple $\tup{c_1, c_2, \ldots, c_m}$
is a greedy $m$-subsequence of $C$
if and only if
the $m$-tuple $\tup{\tup{c_1, r_1}, \tup{c_2, r_2}, \ldots, \tup{c_m, r_m}}$
is a greedy $m$-permutation of $\Ch$.

\item[\textbf{(d)}]
The $m$-tuple $\tup{c_1, c_2, \ldots, c_m}$
has maximum perimeter among all $m$-subsequences of $C$
if and only if
the set $\set{\tup{c_1, r_1}, \tup{c_2, r_2}, \ldots, \tup{c_m, r_m}}$
has maximum perimeter among all $m$-subsets of $\Ch$.

\end{enumerate}
\end{lemma}

The proof of this lemma is just bookkeeping:

\begin{proof}[Proof of Lemma~\ref{lem.Eh.c}.]

\textbf{(a)} This follows from the definitions of
$\wh$ and $\dh$.

\textbf{(b)}
The $m$ elements
$\tup{c_1, r_1}, \tup{c_2, r_2}, \ldots, \tup{c_m, r_m}$
of $\Ch$ are distinct (since $r_1, r_2, \ldots, r_m$
are distinct).
Hence, an application of \eqref{eq.per.tup=set} yields
\begin{align*}
& \PER\set{\tup{c_1, r_1}, \tup{c_2, r_2}, \ldots, \tup{c_m, r_m}} \\
& = \PER\tup{\tup{c_1, r_1}, \tup{c_2, r_2}, \ldots, \tup{c_m, r_m}} \\
& = \PER\tup{c_1, c_2, \ldots, c_m}
 \qquad \tup{\text{by Lemma~\ref{lem.Eh.c} \textbf{(a)}}} .
\end{align*}
This proves Lemma~\ref{lem.Eh.c} \textbf{(b)}.

\textbf{(c)}
%
%
$\Longrightarrow$:
Assume that the $m$-tuple $\tup{c_1, c_2, \ldots, c_m}$
is a greedy $m$-subsequence of $C$.
We must prove that the $m$-tuple
$\tup{\tup{c_1, r_1}, \tup{c_2, r_2}, \ldots, \tup{c_m, r_m}}$
is a greedy $m$-permutation of $\Ch$.
Since the entries of this latter $m$-tuple are distinct
(because $r_1, r_2, \ldots, r_m$ are distinct),
this boils down to proving that
\begin{align}
&\PER\set{\tup{c_1, r_1}, \tup{c_2, r_2}, \ldots, \tup{c_i, r_i}} \nonumber\\
&\geq
\PER\set{\tup{c_1, r_1}, \tup{c_2, r_2}, \ldots, \tup{c_{i-1}, r_{i-1}}, x}
\label{pf.lem.Eh.c.c.1}
\end{align}
for each $i \in \set{1, 2, \ldots, m}$ and each
$x \in \Ch \setminus \set{\tup{c_1, r_1}, \tup{c_2, r_2}, \ldots, \tup{c_{i-1}, r_{i-1}}}$.

So let us prove this.

Fix $i \in \set{1, 2, \ldots, m}$, and fix
$x \in \Ch \setminus \set{\tup{c_1, r_1}, \tup{c_2, r_2}, \ldots, \tup{c_{i-1}, r_{i-1}}}$.
Write $x \in \Ch$ in the form $x = \tup{c', r'}$ for some
$c' \in C$ and $r' \in \IN$.
Since $\tup{c_1, c_2, \ldots, c_m}$
is a greedy $m$-subsequence of $C$, we have
\begin{align}
\PER\tup{c_1, c_2, \ldots, c_i}
\geq
\PER\tup{c_1, c_2, \ldots, c_{i-1}, c'}
\label{pf.lem.Eh.c.c.2}
\end{align}
(by \eqref{eq.def.greedy-seq.geq}, applied to $c'$ instead of $x$).

Clearly, the $i$ elements
$\tup{c_1, r_1}, \tup{c_2, r_2}, \ldots, \tup{c_{i-1}, r_{i-1}}, x$ are distinct
(since $r_1, r_2, \ldots, r_m$ are distinct, and since
$x \in \Ch \setminus \set{\tup{c_1, r_1}, \tup{c_2, r_2}, \ldots, \tup{c_{i-1}, r_{i-1}}}$).
Hence, an application of \eqref{eq.per.tup=set} yields
\begin{align*}
& \PER\set{\tup{c_1, r_1}, \tup{c_2, r_2}, \ldots, \tup{c_{i-1}, r_{i-1}}, x} \\
& = \PER\tup{\tup{c_1, r_1}, \tup{c_2, r_2}, \ldots, \tup{c_{i-1}, r_{i-1}}, x} \\
& = \PER\tup{\tup{c_1, r_1}, \tup{c_2, r_2}, \ldots, \tup{c_{i-1}, r_{i-1}}, \tup{c', r'}}
 \qquad \tup{\text{since $x = \tup{c', r'}$}} \\
& = \PER\tup{c_1, c_2, \ldots, c_{i-1}, c'}
 \qquad \tup{\text{by an application of Lemma~\ref{lem.Eh.c} \textbf{(a)}}} .
\end{align*}
An analogous computation reveals that
\begin{align*}
\PER\set{\tup{c_1, r_1}, \tup{c_2, r_2}, \ldots, \tup{c_i, r_i}}    
= \PER\tup{c_1, c_2, \ldots, c_i} .
\end{align*}
In light of these two equalities, the inequality
\eqref{pf.lem.Eh.c.c.2} (which we know to be true)
turns out to be the same as the inequality
\eqref{pf.lem.Eh.c.c.1} (which we intended to prove).
Thus, \eqref{pf.lem.Eh.c.c.1} is proven, and with it
the ``$\Longrightarrow$'' direction of
Lemma~\ref{lem.Eh.c} \textbf{(c)}.

$\Longleftarrow$:
Assume that the $m$-tuple
$\tup{\tup{c_1, r_1}, \tup{c_2, r_2}, \ldots, \tup{c_m, r_m}}$
is a greedy $m$-permutation of $\Ch$.
We must show that the $m$-tuple $\tup{c_1, c_2, \ldots, c_m}$
is a greedy $m$-subsequence of $C$.
In other words, we must prove that for each
$i \in \set{1, 2, \ldots, m}$ and each $x \in C$, the
inequality \eqref{eq.def.greedy-seq.geq} holds.

So let us fix $i \in \set{1, 2, \ldots, m}$ and $x \in C$.

The $i$ elements
$r_1, r_2, \ldots, r_i$ are distinct
(since $r_1, r_2, \ldots, r_m$ are distinct).
Hence, the $i$ elements
$\tup{c_1, r_1}, \tup{c_2, r_2}, \ldots, \tup{c_{i-1}, r_{i-1}}, \tup{x, r_i}$
of $\Ch$ are distinct;
thus,
$\tup{x, r_i} \in \Ch \setminus \set{\tup{c_1, r_1}, \tup{c_2, r_2}, \ldots, \tup{c_{i-1}, r_{i-1}}}$.
Hence, we can apply \eqref{eq.def.greedy-enum.geq}
to $\tup{\Eh, \wh, \dh}$, $\Ch$,
$\tup{\tup{c_1, r_1}, \tup{c_2, r_2}, \ldots, \tup{c_m, r_m}}$
and $\tup{x, r_i}$ instead of $\tup{E, w, d}$, $C$,
$\tup{c_1, c_2, \ldots, c_m}$ and $x$ (since
$\tup{\tup{c_1, r_1}, \tup{c_2, r_2}, \ldots, \tup{c_m, r_m}}$
is a greedy $m$-permutation of $\Ch$).
Thus, we find
\begin{align}
&\PER\set{\tup{c_1, r_1}, \tup{c_2, r_2}, \ldots, \tup{c_i, r_i}} \nonumber\\
&\geq
\PER\set{\tup{c_1, r_1}, \tup{c_2, r_2}, \ldots, \tup{c_{i-1}, r_{i-1}}, \tup{x, r_i}} .
\label{pf.lem.Eh.c.c.-1}
\end{align}
But we know that the $i$ elements
$\tup{c_1, r_1}, \tup{c_2, r_2}, \ldots, \tup{c_{i-1}, r_{i-1}}, \tup{x, r_i}$ are distinct.
Hence, an application of \eqref{eq.per.tup=set} yields
\begin{align*}
& \PER\set{\tup{c_1, r_1}, \tup{c_2, r_2}, \ldots, \tup{c_{i-1}, r_{i-1}}, \tup{x, r_i}} \\
& = \PER\tup{\tup{c_1, r_1}, \tup{c_2, r_2}, \ldots, \tup{c_{i-1}, r_{i-1}}, \tup{x, r_i}} \\
& = \PER\tup{c_1, c_2, \ldots, c_{i-1}, x}
 \qquad \tup{\text{by an application of Lemma~\ref{lem.Eh.c} \textbf{(a)}}} .
\end{align*}
An analogous computation reveals that
\begin{align*}
\PER\set{\tup{c_1, r_1}, \tup{c_2, r_2}, \ldots, \tup{c_i, r_i}}    
= \PER\tup{c_1, c_2, \ldots, c_i} .
\end{align*}
In light of these two equalities, the inequality
\eqref{pf.lem.Eh.c.c.-1} (which we know to be true)
turns out to be the same as the inequality
\eqref{eq.def.greedy-seq.geq} (which we intended to prove).
Hence, we have shown that \eqref{eq.def.greedy-seq.geq} holds.
This proves the ``$\Longleftarrow$'' direction of
Lemma~\ref{lem.Eh.c} \textbf{(c)}.

\textbf{(d)}
$\Longrightarrow$:
Assume that the $m$-tuple $\tup{c_1, c_2, \ldots, c_m}$
has maximum perimeter among all $m$-subsequences of $C$.
We must prove that the set \newline
$\set{\tup{c_1, r_1}, \tup{c_2, r_2}, \ldots, \tup{c_m, r_m}}$
has maximum perimeter among all $m$-subsets of $\Ch$.
Since this set is clearly an $m$-subset of $\Ch$
(because $r_1, r_2, \ldots, r_m$ are distinct),
this boils down to proving that
\begin{align}
\PER\set{\tup{c_1, r_1}, \tup{c_2, r_2}, \ldots, \tup{c_m, r_m}}
\geq
\PER\tup{G}
\label{pf.lem.Eh.c.d.1}
\end{align}
for each $m$-subset $G$ of $\Ch$.

So let us prove this.
Let $G$ be an $m$-subset of $\Ch$.
Write $G$ in the form
$G = \set{\tup{d_1, q_1}, \tup{d_2, q_2}, \ldots, \tup{d_m, q_m}}$
with $\tup{d_1, q_1}, \tup{d_2, q_2}, \ldots, \tup{d_m, q_m}$ being
$m$ distinct elements of $\Ch$.
Thus,
\begin{align*}
\PER\tup{G}
& = \PER\set{\tup{d_1, q_1}, \tup{d_2, q_2}, \ldots, \tup{d_m, q_m}} \\
& = \PER\tup{\tup{d_1, q_1}, \tup{d_2, q_2}, \ldots, \tup{d_m, q_m}} \\
& \qquad\tup{\text{by an application of \eqref{eq.per.tup=set}}} \\
& = \PER\tup{d_1, d_2, \ldots, d_m}
 \qquad \tup{\text{by an application of Lemma~\ref{lem.Eh.c} \textbf{(a)}}} \\
&\leq \PER\tup{c_1, c_2, \ldots, c_m}
\end{align*}
(since the $m$-tuple $\tup{c_1, c_2, \ldots, c_m}$
has maximum perimeter among all $m$-subsequences of $C$).
In view of
\begin{align*}
& \PER\set{\tup{c_1, r_1}, \tup{c_2, r_2}, \ldots, \tup{c_m, r_m}} \\
& = \PER\tup{\tup{c_1, r_1}, \tup{c_2, r_2}, \ldots, \tup{c_m, r_m}}
 \qquad\tup{\text{by an application of \eqref{eq.per.tup=set}}} \\
& = \PER\tup{c_1, c_2, \ldots, c_m}
 \qquad \tup{\text{by an application of Lemma~\ref{lem.Eh.c} \textbf{(a)}}} ,
\end{align*}
this rewrites as
\[
\PER\tup{G} \leq \PER\set{\tup{c_1, r_1}, \tup{c_2, r_2}, \ldots, \tup{c_m, r_m}} .
\]
Thus, \eqref{pf.lem.Eh.c.d.1} is proven.
This proves the ``$\Longrightarrow$'' direction of
Lemma~\ref{lem.Eh.c} \textbf{(d)}.

$\Longleftarrow$:
Assume that the set
$\set{\tup{c_1, r_1}, \tup{c_2, r_2}, \ldots, \tup{c_m, r_m}}$
has maximum perimeter among all $m$-subsets of $\Ch$.
We must prove that the $m$-tuple $\tup{c_1, c_2, \ldots, c_m}$
has maximum perimeter among all $m$-subsequences of $C$.
In other words, we must prove that
\begin{align}
\PER\tup{c_1, c_2, \ldots, c_m}
\geq \PER\tup{d_1, d_2, \ldots, d_m}
\label{pf.lem.Eh.c.d.-1}
\end{align}
for any $m$-subsequence $\tup{d_1, d_2, \ldots, d_m}$ of $C$.

So let $\tup{d_1, d_2, \ldots, d_m}$ be an $m$-subsequence of $C$.
Then, the elements
$\tup{d_1, r_1}, \tup{d_2, r_2}, \ldots, \tup{d_m, r_m}$
of $\Ch$ are distinct (since $r_1, r_2, \ldots, r_m$ are distinct),
and thus
$\set{\tup{d_1, r_1}, \tup{d_2, r_2}, \ldots, \tup{d_m, r_m}}$
is an $m$-subset of $\Ch$.
Since the set
$\set{\tup{c_1, r_1}, \tup{c_2, r_2}, \ldots, \tup{c_m, r_m}}$
has maximum perimeter among all such $m$-subsets, we thus
obtain
\begin{align*}
\PER\set{\tup{c_1, r_1}, \tup{c_2, r_2}, \ldots, \tup{c_m, r_m}}
\geq
\PER\set{\tup{d_1, r_1}, \tup{d_2, r_2}, \ldots, \tup{d_m, r_m}} .
\end{align*}
In view of
\[
\PER\set{\tup{c_1, r_1}, \tup{c_2, r_2}, \ldots, \tup{c_m, r_m}}
= \PER\tup{c_1, c_2, \ldots, c_m}
 \qquad \tup{\text{by Lemma~\ref{lem.Eh.c} \textbf{(b)}}}
\]
and
\[
\PER\set{\tup{d_1, r_1}, \tup{d_2, r_2}, \ldots, \tup{d_m, r_m}}
= \PER\tup{d_1, d_2, \ldots, d_m}
 \qquad \tup{\text{similarly}} ,
\]
this rewrites as
$\PER\tup{c_1, c_2, \ldots, c_m}
\geq \PER\tup{d_1, d_2, \ldots, d_m}$.
Thus, \eqref{pf.lem.Eh.c.d.-1} is proven.
This proves the ``$\Longleftarrow$'' direction of
Lemma~\ref{lem.Eh.c} \textbf{(d)}.
\end{proof}

\subsection{Proofs of the analogues}

We are now ready to prove the results promised:

\begin{proof}[Proof of Theorem~\ref{thm.multi-2a}.]
Let $k \in \set{0, 1, \ldots, m}$.

Pick any positive integer $N$ such that $N \geq m$.
Define $\IN$, $\Eh$, $\wh$ and $\dh$ as in Subsection~\ref{subsect.multi.clone}.
Pick any $m$ distinct elements $r_1, r_2, \ldots, r_m$
of $\IN$.
(These exist because $N \geq m$;
for example, we can just set $r_i = i$.)
Let $\Ch$ be the subset $C \times \IN$ of $\Eh$.

Lemma~\ref{lem.Eh.c} \textbf{(c)}
(specifically, its ``$\Longrightarrow$'' direction)
shows that
the $m$-tuple $\tup{\tup{c_1, r_1}, \tup{c_2, r_2}, \ldots, \tup{c_m, r_m}}$
is a greedy $m$-permutation of $\Ch$.
Thus, Theorem~\ref{thm.2a} (applied to $\tup{\Eh, \wh, \dh}$,
$\Ch$ and $\tup{c_i, r_i}$ instead of $\tup{E, w, d}$, $C$ and $c_i$)
shows that
the set $\set{\tup{c_1, r_1}, \tup{c_2, r_2}, \ldots, \tup{c_k, r_k}}$
has maximum perimeter among all $k$-subsets of $\Ch$.
Hence, the ``$\Longleftarrow$'' direction of
Lemma~\ref{lem.Eh.c} \textbf{(d)} (applied to $k$ instead of $m$)
shows that the $k$-tuple $\tup{c_1, c_2, \ldots, c_k}$ has
maximum perimeter among all $k$-subsequences of $C$.
This proves Theorem~\ref{thm.multi-2a}.
\end{proof}

Our next task is to prove Theorem~\ref{thm.multi-2b}.
Before we can do this, let us state a straightforward
analogue of Proposition~\ref{prop.greedy.extend}
for greedy $k$-subsequences instead of greedy
$k$-permutations:

\begin{proposition} \label{prop.greedy.multi-extend-b}
Let $m$ and $n$ be integers such that $m \geq n \geq 0$.
Let $C$ be a finite nonempty subset of $E$.

If $\tup{c_1, c_2, \ldots, c_n}$ is a greedy $n$-subsequence of $C$,
then we can find $m-n$ elements $c_{n+1}, c_{n+2}, \ldots, c_m$ of $C$
such that $\tup{c_1, c_2, \ldots, c_m}$ is a greedy $m$-subsequence of $C$.
\end{proposition}

\begin{proof}[Proof of Proposition~\ref{prop.greedy.multi-extend-b}.]
Analogous to the proof of
Proposition~\ref{prop.greedy.extend}.
The main difference is that instead of choosing
$c_i \in C \setminus \set{c_1, c_2, \ldots, c_{i-1}}$
that maximizes $\PER\set{c_1, c_2, \ldots, c_i}$
(in the recursive procedure),
we now have to choose $c_i \in C$
that maximizes $\PER\tup{c_1, c_2, \ldots, c_i}$.
(Such a $c_i$ can always be chosen, since $C$ is nonempty
and finite.)
\end{proof}

The following converse result is obvious, again:

\begin{proposition} \label{prop.greedy.multi-prefix}
Let $C$ be a subset of $E$.
Let $m$ and $n$ be integers such that $m \geq n \geq 0$.

If $\tup{c_1, c_2, \ldots, c_m}$ is a greedy $m$-subsequence of $C$,
then $\tup{c_1, c_2, \ldots, c_n}$ is a greedy $n$-subsequence of $C$.

\end{proposition}

\begin{proof}[Proof of Theorem~\ref{thm.multi-2b}.]
Pick any positive integer $N$ such that $N \geq m$.
Define $\IN$, $\Eh$, $\wh$ and $\dh$ as in Subsection~\ref{subsect.multi.clone}.
Pick any $m$ distinct elements $r_1, r_2, \ldots, r_m$
of $\IN$.
(These exist because $N \geq m$;
for example, we can just set $r_i = i$.)
Let $\Ch$ be the subset $C \times \IN$ of $\Eh$.
Thus, $\abs{\Ch} = \abs{C \times \IN}
= \underbrace{\abs{C}}_{\geq 1} \cdot \underbrace{\abs{\IN}}_{=N}
\geq N \geq m \geq k$.

Recall that $\mathbf{a}$ is a $k$-subsequence of $C$ with
maximum perimeter.
Write this $k$-subsequence $\mathbf{a}$ in the form
$\tup{a_1, a_2, \ldots, a_k}$.
Thus, $\tup{a_1, a_2, \ldots, a_k}$ has maximum perimeter
among all $k$-subsequences of $C$.
Therefore, the ``$\Longrightarrow$'' direction of
Lemma~\ref{lem.Eh.c} \textbf{(d)} (applied to $k$ and $a_i$ instead of $m$ and $c_i$)
shows that
the set $\set{\tup{a_1, r_1}, \tup{a_2, r_2}, \ldots, \tup{a_k, r_k}}$
has maximum perimeter among all $k$-subsets of $\Ch$.
Let us denote this set by $A$.
Thus, $A$ is a $k$-subset of $\Ch$ having maximum
perimeter.
Therefore, Theorem~\ref{thm.2b} (applied to $\tup{\Eh, \wh, \dh}$,
$\Ch$ and $k$ instead of $\tup{E, w, d}$, $C$ and $m$)
shows that there exists a greedy $k$-permutation
$\tup{\tup{c_1, q_1}, \tup{c_2, q_2}, \ldots, \tup{c_k, q_k}}$
of $\Ch$ such that
$A = \set{\tup{c_1, q_1}, \tup{c_2, q_2}, \ldots, \tup{c_k, q_k}}$.
Consider this greedy $k$-permutation.
Hence,
\begin{align*}
\set{\tup{c_1, q_1}, \tup{c_2, q_2}, \ldots, \tup{c_k, q_k}}
= A = \set{\tup{a_1, r_1}, \tup{a_2, r_2}, \ldots, \tup{a_k, r_k}}
\end{align*}
(by the definition of $A$).
Since the $k$ pairs $\tup{a_1, r_1}, \tup{a_2, r_2}, \ldots, \tup{a_k, r_k}$
on the right hand side of this equality are distinct
(because $r_1, r_2, \ldots, r_m$ are distinct),
we thus conclude that the $k$ pairs
$\tup{c_1, q_1}, \tup{c_2, q_2}, \ldots, \tup{c_k, q_k}$ on the left
hand side must also be distinct, and furthermore the former pairs
must be precisely the latter pairs up to order.

In other words, the $k$-tuple
$\tup{\tup{a_1, r_1}, \tup{a_2, r_2}, \ldots, \tup{a_k, r_k}}$
must be a permutation of the $k$-tuple
$\tup{\tup{c_1, q_1}, \tup{c_2, q_2}, \ldots, \tup{c_k, q_k}}$.
Hence, the $k$-tuple
$\tup{a_1, a_2, \ldots, a_k}$
must be a permutation of the $k$-tuple
$\tup{c_1, c_2, \ldots, c_k}$.
In other words, the $k$-tuple $\mathbf{a}$ is
a permutation of the $k$-tuple
$\tup{c_1, c_2, \ldots, c_k}$
(since $\mathbf{a} = \tup{a_1, a_2, \ldots, a_k}$).

Also, the $k$-tuple $\tup{r_1, r_2, \ldots, r_k}$ is a
permutation of the $k$-tuple $\tup{q_1, q_2, \ldots, q_k}$
(since the $k$-tuple
$\tup{\tup{a_1, r_1}, \tup{a_2, r_2}, \ldots, \tup{a_k, r_k}}$
is a permutation of the $k$-tuple
$\tup{\tup{c_1, q_1}, \tup{c_2, q_2}, \ldots, \tup{c_k, q_k}}$).
Hence, $q_1, q_2, \ldots, q_k$ are distinct
(since $r_1, r_2, \ldots, r_k$ are distinct).
Therefore, the ``$\Longleftarrow$'' direction of
Lemma~\ref{lem.Eh.c} \textbf{(c)} (applied to $k$
and $q_i$ instead of $m$ and $r_i$)
shows that $\tup{c_1, c_2, \ldots, c_k}$ is a greedy
$k$-subsequence of $C$
(since $\tup{\tup{c_1, q_1}, \tup{c_2, q_2}, \ldots, \tup{c_k, q_k}}$ is
a greedy $k$-permutation of $\Ch$).
Since $k \leq m$,
we can extend this greedy $k$-subsequence
to a greedy $m$-subsequence
$\tup{c_1, c_2, \ldots, c_m}$ of $C$
(by Proposition~\ref{prop.greedy.multi-extend-b},
applied to $n = k$).
Hence, we have found a greedy $m$-subsequence
$\tup{c_1, c_2, \ldots, c_m}$ of $C$ such that
$\mathbf{a}$ is
a permutation of the $k$-tuple
$\tup{c_1, c_2, \ldots, c_k}$.
This proves Theorem~\ref{thm.multi-2b}.
\end{proof}

\begin{proof}[Proof of Corollary~\ref{cor.multi-bh-t1}.]
Pick any positive integer $N$ such that $N \geq k$.
Define $\IN$, $\Eh$, $\wh$ and $\dh$ as in Subsection~\ref{subsect.multi.clone}.
Let $\Ch$ be the subset $C \times \IN$ of $\Eh$.

Let $\tup{c_1, c_2, \ldots, c_m}$ be a greedy
$m$-subsequence of $C$.
Then, $\tup{c_1, c_2, \ldots, c_k}$ is a greedy
$k$-subsequence of $C$
(by Proposition~\ref{prop.greedy.multi-prefix},
applied to $n = k$).

Recall the number $\nuo_k\tup{C}$ we defined just
after Corollary~\ref{cor.bh-t1}.
Now, consider the number $\nuo_k\tup{\Ch}$ defined in
the same fashion, but with
respect to the ultra triple $\tup{\Eh, \wh, \dh}$.
We claim that
\begin{align}
w\tup{c_k} + \sum_{i=1}^{k-1} d\tup{c_i, c_k}
= \nuo_k\tup{\Ch} .
\label{pf.cor.multi-bh-t1.goal}
\end{align}
Clearly, proving this will yield
Corollary~\ref{cor.multi-bh-t1}.

Pick any $k$ distinct elements $r_1, r_2, \ldots, r_k$
of $\IN$.
(These exist because $N \geq k$;
for example, we can just set $r_i = i$.)

Recall that the $k$-tuple $\tup{c_1, c_2, \ldots, c_k}$
is a greedy $k$-subsequence of $C$.
Hence,
the $k$-tuple $\tup{\tup{c_1, r_1}, \tup{c_2, r_2}, \ldots, \tup{c_k, r_k}}$
is a greedy $k$-permutation of $\Ch$
(by the ``$\Longrightarrow$'' direction of
Lemma~\ref{lem.Eh.c} \textbf{(c)}, applied to $k$
instead of $m$).
Hence, the definition of $\nuo_k\tup{\Ch}$ yields
\begin{align*}
\nuo_k\tup{\Ch}
&= \underbrace{\wh\tup{c_k, r_k}}_{\substack{= w\tup{c_k} \\ \text{(by the definition of $\wh$)}}} + \sum_{i=1}^{k-1} \underbrace{\dh \tup{\tup{c_i, r_i}, \tup{c_k, r_k}}}_{\substack{= d\tup{c_i, c_k} \\ \text{(by the definition of $\dh$)}}}
= w\tup{c_k} + \sum_{i=1}^{k-1} d\tup{c_i, c_k} .
\end{align*}
This proves \eqref{pf.cor.multi-bh-t1.goal}.
Hence, Corollary~\ref{cor.multi-bh-t1} is proven.
\end{proof}

\begin{proof}[Proof of Corollary~\ref{cor.multi-bh-l2}.]
Pick any positive integer $N$ such that $N \geq k$.
Define $\IN$, $\Eh$, $\wh$ and $\dh$ as in Subsection~\ref{subsect.multi.clone}.
Let $\Ch$ be the subset $C \times \IN$ of $\Eh$.

Recall the number $\nuo_k\tup{C}$ we defined just
after Corollary~\ref{cor.bh-t1}.
Now, consider the number $\nuo_k\tup{\Ch}$ defined in
the same fashion, but with
respect to the ultra triple $\tup{\Eh, \wh, \dh}$.

Pick any $k$ distinct elements $r_1, r_2, \ldots, r_k$
of $\IN$.
(These exist because $N \geq k$;
for example, we can just set $r_i = i$.)

The $k$-tuple $\tup{c_1, c_2, \ldots, c_k}$ is a greedy
$k$-subsequence of $C$
(by Proposition~\ref{prop.greedy.multi-prefix},
applied to $n = k$).
Hence,
the $k$-tuple $\tup{\tup{c_1, r_1}, \tup{c_2, r_2}, \ldots, \tup{c_k, r_k}}$
is a greedy $k$-permutation of $\Ch$
(by the ``$\Longrightarrow$'' direction of
Lemma~\ref{lem.Eh.c} \textbf{(c)}, applied to $k$
instead of $m$).
Hence, Corollary~\ref{cor.bh-l2} (applied to
$\tup{\Eh, \wh, \dh}$, $\Ch$, $k$ and $\tup{c_i, r_i}$
instead of $\tup{E, w, d}$, $C$, $m$ and $c_i$) yields
\begin{align}
\nuo_k\tup{\Ch}
&\leq \underbrace{\wh\tup{c_j, r_j}}_{\substack{= w\tup{c_j} \\ \text{(by the definition of $\wh$)}}} + \sum_{i \in \set{1, 2, \ldots, k} \setminus \set{j}} \underbrace{\dh\tup{\tup{c_i, r_i}, \tup{c_j, r_j}}}_{\substack{= d\tup{c_i, c_j} \\ \text{(by the definition of $\dh$)}}} \nonumber \\
&= w\tup{c_j} + \sum_{i \in \set{1, 2, \ldots, k} \setminus \set{j}} d\tup{c_i, c_j} .
\label{pf.cor.multi-bh-l2.4}
\end{align}
But in the proof of Corollary~\ref{cor.multi-bh-t1},
we have seen that
\begin{align*}
\nuo_k\tup{\Ch}
= w\tup{c_k} + \sum_{i=1}^{k-1} d\tup{c_i, c_k} .
\end{align*}
Finally, the definition of $\nu_k\tup{C}$ yields
\[
\nu_k\tup{C} = w\tup{c_k} + \sum_{i=1}^{k-1} d\tup{c_i, c_k}
= \nuo_k\tup{\Ch}
\leq w\tup{c_j} + \sum_{i \in \set{1, 2, \ldots, k} \setminus \set{j}} d\tup{c_i, c_j}
\]
(by \eqref{pf.cor.multi-bh-l2.4}).
This proves Corollary~\ref{cor.multi-bh-l2}.
\end{proof}

\begin{remark}
Lemma~\ref{lem.Eh.c} \textbf{(c)} essentially says that,
using the full ultra triple $\tup{\Eh, \wh, \dh}$,
we can re-interpret greedy $m$-subsequences
as (a certain subclass of)
greedy $m$-permutations (as long as $N$ is chosen to
satisfy $N \geq m$).

The reverse direction can also be done:
We can re-interpret greedy $m$-permutations of $C$
as greedy $m$-subsequences, as long as $C$ is finite
and satisfies $\abs{C} \geq m$.
To do so, we fix a real number $R$ such that
$R > 2 \abs{\PER\tup{D}}$ for every $D \subseteq C$.
We define a new distance function $d_R : E \times E \to \RR$
on $E$ by setting
\[
d_R \tup{e, f}
=
\begin{cases}
d \tup{e, f} + R, & \text{ if } e \neq f ; \\
d \tup{e, f},     & \text{ if } e = f
\end{cases}
\qquad \text{ for all } e, f \in E .
\]
It is easy to see that $\tup{E, w, d_R}$ is again a
full ultra triple.
Moreover, it is easy to see that any
$m$-subsequence of $C$ containing two equal entries
has smaller perimeter with respect to $\tup{E, w, d_R}$
than any $m$-subset of $C$.
Hence, the maximum perimeter of an $m$-subsequence of $C$
with respect to $\tup{E, w, d_R}$
can only be achieved by an $m$-subsequence with no equal
entries.
Hence, this maximum perimeter is the maximum perimeter
of an $m$-subset of $C$ with respect to $\tup{E, w, d_R}$.
Meanwhile, the perimeter of an $m$-subset of $C$
with respect to $\tup{E, w, d_R}$ equals its perimeter
with respect to the original full ultra triple $\tup{E, w, d}$
plus the constant $\dbinom{m}{2} R$.
Hence, the $m$-subsets of $C$ having maximum perimeter
with respect to $\tup{E, w, d_R}$ are precisely the same
as the ones
that have maximum perimeter with respect to $\tup{E, w, d}$.
From this, it is easy to see that the greedy
$m$-subsequences of $C$ with respect to $\tup{E, w, d_R}$
are precisely the greedy $m$-permutations of $C$ with
respect to $\tup{E, w, d}$.

When $\abs{C} < m$, this reasoning no longer works, since
every $m$-subsequence of $C$ has two equal entries
(and there are no $m$-subsets of $C$).
In this case, the greedy $m$-subsequences of $C$ with
respect to $\tup{E, w, d_R}$ can be informally regarded as greedy
$m$-subsequences of $C$ with respect to $\tup{E, w, d}$
that defer picking identical entries as long as they can
(in a sense).
\end{remark}

\section{\label{sect.bhargava}Relation to Bhargava's $P$-orderings}

Let us now explain the connection between greedy
$m$-permutations and
the concept of $P$-orderings introduced by Manjul
Bhargava in \cite[Section 2]{Bharga97}.
(The notions of $p$-orderings in
\cite[Section 4]{Bharga00} and \cite[Section 2]{Bharga09}
are particular cases.)
This connection was already noticed by Bhargava
(see the paragraph after the proof of Lemma 2 in \cite{Bharga97}),
who, however, never elaborated on it or made
any further inroads into the study of general ultra triples.

We fix a Dedekind ring%
\footnote{See \cite[Chapter 1]{Narkie04}
or \cite[Chapter 3]{Ash-ANT} for an introduction
to Dedekind rings.
In a nutshell, a \defn{Dedekind ring} is an integral
domain in which every nonzero ideal has a unique factorization
into a product of prime ideals.
Other equivalent definitions of Dedekind rings exist.
Dedekind rings are also known as \defn{Dedekind domains}.

For our purposes, it suffices to know that
$\ZZ$ is a Dedekind ring;
the examples it provides are sufficiently rich
in substance that greater generality is not strictly necessary.}
$R$
and a nonzero prime ideal $P$ of $R$.
For any nonzero $a \in R$, we let
\defn{$v_P\tup{a}$} denote the
highest\footnote{Here and in the following, we set
$\NN = \set{0,1,2,\ldots}$.}
$k \in \NN$ that satisfies $a \in P^k$.
(Equivalently, $v_P\tup{a}$ is the exponent with
which $P$ appears in the factorization of the
principal ideal $aR$ into prime ideals.%
\footnote{The equivalence between these two
definitions of $v_P\tup{a}$ follows from
\cite[Corollary 3.3.3]{Ash-ANT}; this also
proves that the first definition is valid
(i.e., there exists a highest $k \in \NN$
that satisfies $a \in P^k$).})
We also set $v_P\tup{0} = +\infty$.
Thus, an element $v_P\tup{a} \in \NN \cup \set{+ \infty}$
is defined for every $a \in R$.
Moreover, the map $v_P : R \to \NN \cup \set{+ \infty}$
satisfies
\[
v_P\tup{ab} = v_P\tup{a} + v_P\tup{b}
\qquad \text{ and } \qquad
v_P\tup{a+b} \geq \min\set{v_P\tup{a}, v_P\tup{b}}
\]
for all $a, b \in R$.

The simplest example for this is when $R = \ZZ$ and
$P = p \ZZ$ for some prime number $p$.
In this case, $v_P\tup{a} = v_p\tup{a}$, where $v_p\tup{a}$
is defined as in Example~\ref{exa.p-adic.d}.
This particular case is the one studied in
\cite[Section 4]{Bharga00} and \cite[Section 2]{Bharga09}.

Furthermore, we fix a nonempty subset $E$ of $R$.
(Bhargava denotes this subset by $X$ instead.)
Now, Bhargava defines a \defn{$P$-ordering} of $E$ to be a
sequence $\tup{a_0, a_1, a_2, \ldots}$ of elements of $E$
defined recursively as follows:
For each $k \in \NN$, we define $a_k$ (assuming that
$a_0, a_1, \ldots, a_{k-1}$ are already defined) to be an
element of $E$ minimizing the quantity
\begin{align}
v_P \tup{\tup{a_k-a_0} \tup{a_k-a_1} \cdots \tup{a_k-a_{k-1}}} .
\label{eq.bhargava.prodai}
\end{align}

Note that the quantity \eqref{eq.bhargava.prodai} indeed
attains its minimum at some (usually non-unique) $a_k \in E$,
since it is an element of the well-ordered set $\NN \cup \set{+ \infty}$.

We now claim that this notion of $P$-ordering is
almost a particular case of the notion of a greedy
$m$-permutation for a certain ultra triple.
Some amount of work is necessary to bridge the
technical discrepancies between these two notions:

First of all, $P$-orderings are infinite sequences,
whereas greedy $m$-permutations are $m$-tuples.
To bring them closer together, we fix an $m \in \NN$,
and we define an
\defn{$\tup{P,m}$-ordering} of $E$ to be an $m$-tuple
$\tup{a_0, a_1, \ldots, a_{m-1}}$ of
elements of $E$ such that for each
$k \in \set{0, 1, \ldots, m-1}$, the element $a_k$
of $E$ minimizes the quantity \eqref{eq.bhargava.prodai}
(where $a_0, a_1, \ldots, a_{k-1}$ are considered fixed).
Clearly, the first $m$ entries of any $P$-ordering
form a $\tup{P, m}$-ordering, and conversely, any
$\tup{P, m}$-ordering can be extended to a $P$-ordering.
Thus, if we want to study (finitary) properties of
$P$-orderings, it suffices to understand $\tup{P, m}$-orderings.
Thus we are back in the realm of finite sequences.

We furthermore notice something simple:

\begin{lemma} \label{lem.bhargava.no-rep}
Let $C$ be a subset of $E$, and let $m$ be a nonnegative
integer such that $\abs{C} \geq m$.
Then, any $\tup{P, m}$-ordering of $C$
is an $m$-tuple of distinct elements.
\end{lemma}

\begin{proof}[Proof of Lemma~\ref{lem.bhargava.no-rep} (sketched).]
Any $\tup{P, m}$-ordering
$\tup{a_0, a_1, \ldots, a_{m-1}}$ of $C$ can be
constructed recursively as follows:
For each $k \in \set{0, 1, \ldots, m-1}$,
we define $a_k$ (assuming that
$a_0, a_1, \ldots, a_{k-1}$ are already defined) to be an
element of $C$ minimizing the quantity
\eqref{eq.bhargava.prodai}.
But this quantity \eqref{eq.bhargava.prodai} is $+ \infty$
when $a_k$ equals one of $a_0, a_1, \ldots, a_{k-1}$, and
otherwise is a nonnegative integer.
Hence, an $a_k$ that equals one of $a_0, a_1, \ldots, a_{k-1}$
cannot minimize this quantity (as long as there is
at least one element of $C$ that does not equal any of
$a_0, a_1, \ldots, a_{k-1}$;
but this is always guaranteed thanks to $\abs{C} \geq m > k$).
Thus, any $a_k$ chosen in the construction of a
$\tup{P, m}$-ordering of $C$
must be distinct from $a_0, a_1, \ldots, a_{k-1}$.
Hence, if $\tup{a_0, a_1, \ldots, a_{m-1}}$ is any
$\tup{P, m}$-ordering of $C$, then
$a_0, a_1, \ldots, a_{m-1}$ are distinct.
This proves Lemma~\ref{lem.bhargava.no-rep}.
\end{proof}

Next, we define an ultra triple $\tup{E, w, d'}$ as follows:
We define the weight function $w : E \to \RR$
by setting $w\tup{e} = 0$ for all $e \in E$.
We define a map $d' : E \timesu E \to \RR$ by setting
\[
d' \tup{a, b}
= -v_P\tup{a-b}
\qquad \text{for all } \tup{a, b} \in E \timesu E .
\]
(This generalizes the map $d'$ from Example~\ref{exa.p-adic.d'}.)

Now, $\tup{E, w, d'}$ is an ultra triple.
Throughout this section,
we shall always be using this ultra triple
(when we speak, e.g., of greedy $m$-permutations).
We claim the following:

\begin{proposition} \label{prop.bhargava.Pmo=gmp}
Let $C$ be a subset of $E$.
Let $m \in \NN$.
Let $c_1, c_2, \ldots, c_m \in C$ be distinct.
Then, the $m$-tuple $\tup{c_1, c_2, \ldots, c_m}$
is a greedy $m$-permutation of $C$ if and only if
it is a $\tup{P, m}$-ordering of $C$.
\end{proposition}

\begin{proof}[Proof of Proposition~\ref{prop.bhargava.Pmo=gmp} (sketched).]
We have $\abs{C} \geq m$
(since $C$ has at least the
$m$ distinct elements $c_1, c_2, \ldots, c_m$).

The entries $c_1, c_2, \ldots, c_m$ of the $m$-tuple
$\tup{c_1, c_2, \ldots, c_m}$
are distinct.
Hence, the definition of a greedy $m$-permutation
yields the following:

\begin{statement}
\textit{Claim 1:}
The $m$-tuple $\tup{c_1, c_2, \ldots, c_m}$
is a greedy $m$-permutation of $C$ if and only if
for each $i \in \set{1,2,\ldots,m}$
and each $x \in C \setminus \set{c_1, c_2, \ldots, c_{i-1}}$,
the inequality \eqref{eq.def.greedy-enum.geq} holds.
\end{statement}

On the other hand,
the definition of a $\tup{P,m}$-ordering
shows that $\tup{c_1, c_2, \ldots, c_m}$
is a $\tup{P,m}$-ordering of $C$ if and only if
for each
$k \in \set{0, 1, \ldots, m-1}$, the element $c_{k+1}$
of $C$ minimizes the quantity
$v_P \tup{\tup{c_{k+1}-c_1} \tup{c_{k+1}-c_2} \cdots \tup{c_{k+1}-c_k}}$
(where $c_1, c_2, \ldots, c_k$ are considered fixed).
Substituting $i-1$ for $k$ in this statement,
we obtain the following:
The $m$-tuple $\tup{c_1, c_2, \ldots, c_m}$
is a $\tup{P,m}$-ordering of $C$ if and only if
for each
$i \in \set{1, 2, \ldots, m}$, the element $c_i$
of $C$ minimizes the quantity
$v_P \tup{\tup{c_i-c_1} \tup{c_i-c_2} \cdots \tup{c_i-c_{i-1}}}$
(where $c_1, c_2, \ldots, c_{i-1}$ are considered fixed).
We can restate this as follows:

\begin{statement}
\textit{Claim 2:}
The $m$-tuple $\tup{c_1, c_2, \ldots, c_m}$
is a $\tup{P,m}$-ordering of $C$ if and only if
for each
$i \in \set{1, 2, \ldots, m}$
and each $x \in C$, the inequality
\begin{align}
& v_P \tup{\tup{c_i-c_1} \tup{c_i-c_2} \cdots \tup{c_i-c_{i-1}}}
\nonumber\\
& \leq
v_P \tup{\tup{x-c_1} \tup{x-c_2} \cdots \tup{x-c_{i-1}}}
\label{pf.prop.bhargava.Pmo=gmp.Pmo1}
\end{align}
holds.
\end{statement}

Note that if $i \in \set{1, 2, \ldots, m}$
and $x \in \set{c_1, c_2, \ldots, c_{i-1}}$,
then the inequality \eqref{pf.prop.bhargava.Pmo=gmp.Pmo1}
automatically holds\footnote{because in this case, we have
$\tup{x-c_1} \tup{x-c_2} \cdots \tup{x-c_{i-1}} = 0$
and thus
$v_P \tup{\tup{x-c_1} \tup{x-c_2} \cdots \tup{x-c_{i-1}}}
= v_P \tup{0} = + \infty$}.
Therefore, if $i \in \set{1, 2, \ldots, m}$ is given,
then the inequality \eqref{pf.prop.bhargava.Pmo=gmp.Pmo1}
holds for each $x \in C$ if and only if it holds for each
$x \in C \setminus \set{c_1, c_2, \ldots, c_{i-1}}$.
Hence, in Claim 2, we can replace
``each $x \in C$'' by ``each
$x \in C \setminus \set{c_1, c_2, \ldots, c_{i-1}}$''.
Thus, Claim 2 rewrites as follows:

\begin{statement}
\textit{Claim 3:}
The $m$-tuple $\tup{c_1, c_2, \ldots, c_m}$
is a $\tup{P,m}$-ordering of $C$ if and only if
for each
$i \in \set{1, 2, \ldots, m}$
and each $x \in C \setminus \set{c_1, c_2, \ldots, c_{i-1}}$,
the inequality \eqref{pf.prop.bhargava.Pmo=gmp.Pmo1}
holds.
\end{statement}

For any $i \in \set{1, 2, \ldots, m}$
and $x \in C \setminus \set{c_1, c_2, \ldots, c_{i-1}}$,
we have the following chain of logical equivalences:
\begin{align}
& \tup{\text{the inequality \eqref{eq.def.greedy-enum.geq} holds}}
\nonumber \\
& \Longleftrightarrow
\tup{\PER\set{c_1, c_2, \ldots, c_i} \geqslant \PER\set{c_1, c_2, \ldots, c_{i-1}, x} }
\nonumber \\
& \Longleftrightarrow
\tup{w\left(c_i\right) + \sum_{j=1}^{i-1} d'\tup{c_i, c_j} \geqslant w\left(x\right) + \sum_{j=1}^{i-1} d'\tup{x, c_j} }
\nonumber \\
& \qquad\qquad \tup{\begin{array}{c}
\text{here, we have subtracted $\PER\set{c_1, c_2, \ldots, c_{i-1}}$} \\
\text{ from both sides of the inequality}
\end{array}}
\nonumber \\
& \Longleftrightarrow
\tup{\sum_{j=1}^{i-1} d'\tup{c_i, c_j} \geqslant \sum_{j=1}^{i-1} d'\tup{x, c_j} }
\qquad \tup{\text{since $w\tup{e} = 0$ for all $e \in E$}}
\nonumber \\
& \Longleftrightarrow
\tup{\sum_{j=1}^{i-1} \tup{- v_P\tup{c_i - c_j}} \geqslant \sum_{j=1}^{i-1} \tup{- v_P\tup{x - c_j}} }
 \qquad \tup{\text{by the definition of $d'$}}
\nonumber \\
& \Longleftrightarrow
\tup{\sum_{j=1}^{i-1} v_P\tup{c_i - c_j} \leq \sum_{j=1}^{i-1} v_P\tup{x - c_j} }
\nonumber \\
& \Longleftrightarrow
\tup{v_P\tup{\prod_{j=1}^{i-1}\tup{c_i - c_j}} \leq v_P\tup{\prod_{j=1}^{i-1}\tup{x - c_j}} }
\nonumber \\
& \qquad\qquad \tup{\begin{array}{c}
    \text{since $\sum_{j \in J} v_P\tup{a_j} = v_P\tup{\prod_{j \in J} a_j}$} \\
    \text{for any finite family $\tup{a_j}_{j\in J}$ of elements of $R$}
    \end{array} } \nonumber \\
& \Longleftrightarrow
\tup{\text{the inequality \eqref{pf.prop.bhargava.Pmo=gmp.Pmo1} holds}}
\label{pf.prop.bhargava.Pmo=gmp.equiv1}
\end{align}
(since $\prod_{j=1}^{i-1}\tup{c_i - c_j} = \tup{c_i-c_1} \tup{c_i-c_2} \cdots \tup{c_i-c_{i-1}}$
and $\prod_{j=1}^{i-1}\tup{x - c_j} = \tup{x-c_1} \tup{x-c_2} \cdots \tup{x-c_{i-1}}$).

Now, we have the following chain of logical equivalences:
\begin{align*}
& \tup{\text{$\tup{c_1, c_2, \ldots, c_m}$ is a greedy $m$-permutation of $C$}} \\
& \Longleftrightarrow
\left( \text{\eqref{eq.def.greedy-enum.geq} holds for each $i \in \set{1,2,\ldots,m}$} \right. \\
& \qquad \qquad \qquad \left. \text{ and each $x \in C \setminus \set{c_1, c_2, \ldots, c_{i-1}}$} \right)
\qquad \tup{\text{by Claim 1}} \\
& \Longleftrightarrow
\left( \text{\eqref{pf.prop.bhargava.Pmo=gmp.Pmo1} holds for each $i \in \set{1,2,\ldots,m}$} \right. \\
& \qquad \qquad \qquad \left. \text{ and each $x \in C \setminus \set{c_1, c_2, \ldots, c_{i-1}}$} \right)
\qquad \tup{\text{by \eqref{pf.prop.bhargava.Pmo=gmp.equiv1}}} \\
& \Longleftrightarrow
\tup{\text{$\tup{c_1, c_2, \ldots, c_m}$ is a $\tup{P,m}$-ordering of $C$}}
\qquad \tup{\text{by Claim 3}}.
\end{align*}
Hence, the $m$-tuple $\tup{c_1, c_2, \ldots, c_m}$
is a greedy $m$-permutation of $C$ if and only if
it is a $\tup{P, m}$-ordering of $C$.
This proves Proposition~\ref{prop.bhargava.Pmo=gmp}.
\end{proof}

Equipped with Proposition~\ref{prop.bhargava.Pmo=gmp},
we can now translate each result about greedy
$m$-permutations into the language of
$\tup{P, m}$-orderings as long as $\abs{C} \geq m$
(because Lemma~\ref{lem.bhargava.no-rep} shows that
any $\tup{P, m}$-ordering consists of distinct entries
in this case).\footnote{The case $\abs{C} < m$
is a degenerate case which can easily be
reduced to the case $\abs{C} \geq m$
by focussing only on the first $\abs{C}$ many entries
of the $\tup{P, m}$-ordering.
(All the other entries merely repeat the first
$\abs{C}$ many entries, in an arbitrary way,
so there is nothing of interest to say about them.)}
In particular,
Corollary~\ref{cor.bh-t1} becomes \cite[Theorem 1]{Bharga97},
while
Corollary~\ref{cor.bh-l2} becomes \cite[Lemma 2]{Bharga97}.
(More precisely, we obtain the analogues of
\cite[Theorem 1]{Bharga97} and \cite[Lemma 2]{Bharga97}
for $\tup{P, m}$-orderings instead of $P$-orderings.
But since the first $m$ entries of any $P$-ordering form
a $\tup{P, m}$-ordering, these analogues immediately yield
\cite[Theorem 1]{Bharga97} and \cite[Lemma 2]{Bharga97}.)

We note in passing that
the ``$P$-orderings of order $h$'' defined in
\cite[Section 2.2]{Bharga09} can also be regarded
as a particular case of greedy $m$-permutations
(up to the already mentioned technicalities);
we only need to modify the distance function $d'$.

\section{\label{sect.app-p}Appendix: Greediness of $\tup{1, 2, \ldots, m}$ for $p$-adic metrics}

In this section, we shall prove the claim made in
Example~\ref{exa.p-adic.greedy}.

We begin with a basic folklore lemma about inequalities.

\begin{lemma} \label{lem.app-p.1}
Let $I$ be a finite set.
Let $P$ be a totally ordered set.
For each $i \in I$, let $a_i$ and $b_i$ be two elements of $P$.
Assume that each $h \in P$ satisfies
\begin{align}
\abs{\set{i \in I \mid a_i \geq h}}
\geq
\abs{\set{i \in I \mid b_i \geq h}} .
\label{eq.lem.app-p.1.ass}
\end{align}

Let $f : P \to \RR$ be any weakly increasing
map\footnote{``Weakly increasing'' means that
$f\tup{p} \leq f\tup{q}$ for all $p, q \in P$
satisfying $p \leq q$.}.
Then,
\begin{equation}\label{l101}
  \sum_{i \in I} f\tup{a_i} \geq \sum_{i \in I} f\tup{b_i}.
\end{equation}
\end{lemma}

The quickest proof of Lemma~\ref{lem.app-p.1}
uses integration of step functions:

\begin{proof}[Proof of Lemma~\ref{lem.app-p.1}.]
We may assume that $P$ is finite (otherwise, simply
replace $P$ by the finite set
$\set{a_i \mid i \in I} \cup \set{b_i \mid i \in I}$).

Conditions and conclusion of Lemma~\ref{lem.app-p.1}
do not change if we replace all values $f\tup{x}$
of $f$ by $f\tup{x}+C$ for a real constant $C$.
This observation allows us to assume that
$f$ takes non-negative real values only.
Assuming this and denoting $M=\max_P f$, we have
\[
f\tup{x}=\int_0^M \chi_{[0,f\tup{x}]}\tup{\tau}
d\tau
\qquad \text{for all $x \in P$,}
\]
where $\chi_{[a,b]}$ denotes 
the characteristic function of the interval $[a,b]$. Therefore
\begin{align}
\sum_{i \in I} \underbrace{f\tup{a_i}}_{= \int_0^M \chi_{[0,f\tup{a_i}]}\tup{\tau} d\tau}
&=
\int_0^M \underbrace{\sum_{i \in I}
\chi_{[0,f\tup{a_i}]}\tup{\tau}}_{= \abs{\set{i\in I \mid f\tup{a_i}\geqslant \tau}}}
d\tau
\nonumber\\
&=
\int_0^M \abs{\set{i\in I \mid f\tup{a_i}\geqslant \tau}} d\tau.
\label{pf.lem.app-p.1.4}
\end{align}
But recall that $f$ is weakly increasing.
Hence, for any fixed real $\tau\in [0,M]$ and any $x \in P$,
the condition $f\tup{x}\geqslant \tau$
may be rewritten as $x\geqslant h\tup{\tau}$,
where
$h\tup{\tau} = \min \set{y \in P \mid f\tup{y}\geqslant \tau}$
(note that this minimum exists because
$\tau\leqslant M=\max_P f$).
Therefore, every $\tau\in [0,M]$ satisfies
\begin{align*}
\abs{\set{i\in I \mid f\tup{a_i}\geqslant \tau}}
&=
\abs{\set{i\in I \mid a_i\geqslant h\tup{\tau}}}
\\
&\geqslant \abs{\set{i\in I \mid b_i\geqslant h\tup{\tau}}}
\qquad \tup{\text{by \eqref{eq.lem.app-p.1.ass}}} \nonumber\\
&=
\abs{\set{i\in I \mid f\tup{b_i}\geqslant \tau}}.
\end{align*}
Integrating this inequality over $\left[0,M\right]$,
and rewriting the result using \eqref{pf.lem.app-p.1.4},
we get \eqref{l101}.
\end{proof}

\begin{verlong}
What follows is a second proof of Lemma~\ref{lem.app-p.1},
which was written for the original version of this preprint
and appeared in this form in arXiv:1909.01965v2.

\begin{proof}[Second proof of Lemma~\ref{lem.app-p.1}.]
Let $P'$ be the subset
$\set{a_i \mid i \in I} \cup \set{b_i \mid i \in I}$
of $P$.
This subset $P'$ is finite; let us thus write it in the form
$P' = \set{p_1 < p_2 < \cdots < p_m}$
for some $p_1, p_2, \ldots, p_m \in P$
(since $P$ is totally ordered).
Clearly, $a_i \in P'$ and $b_i \in P'$
for each $i \in I$.

For each $k \in \set{1, 2, \ldots, m}$, we
set $f_k = f\tup{p_k} \in \RR$.
From $p_1 < p_2 < \cdots < p_m$,
we obtain $f\tup{p_1} \leq f\tup{p_2} \leq \cdots \leq f\tup{p_m}$
(since $f$ is weakly increasing).
This rewrites as
$f_1 \leq f_2 \leq \cdots \leq f_m$
(since $f_k = f\tup{p_k}$ for all $k$).

Extend the $m$-tuple $\tup{f_1, f_2, \ldots, f_m}$
of real numbers to an $\tup{m+1}$-tuple $\tup{f_0, f_1, \ldots, f_m}$
by choosing an arbitrary real number $f_0$ that satisfies
$f_0 \leq f_k$ for all $k \in \set{1, 2, \ldots, m}$.
Thus, $f_0 \leq f_1 \leq \cdots \leq f_m$
(since $f_1 \leq f_2 \leq \cdots \leq f_m$).

For each $k \in \set{1, 2, \ldots, m}$, set $h_k = f_k - f_{k-1}$.
This number $h_k$ satisfies
\begin{equation}
h_k = f_k - f_{k-1} \geq 0
\label{pf.lem.app-p.1.hk0}
\end{equation}
(since $f_0 \leq f_1 \leq \cdots \leq f_m$, thus
$f_{k-1} \leq f_k$, thus $f_k - f_{k-1} \geq 0$)
and thus
\begin{align}
\underbrace{\abs{\set{i \in I \mid a_i \geq p_k}}}_{\substack{\geq \abs{\set{i \in I \mid b_i \geq p_k}} \\ \tup{\text{by \eqref{eq.lem.app-p.1.ass}, applied to $h = p_k$}}}} h_k
\geq
\abs{\set{i \in I \mid b_i \geq p_k}} h_k .
\label{pf.lem.app-p.1.ass2}
\end{align}
Moreover, each $j \in \set{0, 1, \ldots, m}$ satisfies
\begin{align}
\sum_{k = 1}^j \underbrace{h_k}_{= f_k - f_{k-1}}
& = \sum_{k = 1}^j \tup{f_k - f_{k-1}}
= \underbrace{f_j}_{\substack{= f\tup{p_j} \\ \tup{\text{by the definition of } f_j}}} - f_0 \nonumber\\
& \qquad \tup{\text{by the telescope principle}} \nonumber \\
&= f\tup{p_j} - f_0 .
\label{pf.lem.app-p.1.sumhk}
\end{align}

Now, let $i \in I$. Then,
$a_i \in P' = \set{p_1 < p_2 < \cdots < p_m}$.
Hence, there exists a unique $j \in \set{1, 2, \ldots, m}$
satisfying $a_i = p_j$.
Therefore,
\[
\sum_{\substack{j \in \set{1, 2, \ldots, m}; \\ a_i = p_j}} \tup{f\tup{a_i} - f_0}
= f\tup{a_i} - f_0 ,
\]
so that
\begin{align}
f\tup{a_i} - f_0
&= \sum_{\substack{j \in \set{1, 2, \ldots, m}; \\ a_i = p_j}} \tup{f\tup{\underbrace{a_i}_{= p_j}} - f_0} \nonumber\\
&= \sum_{\substack{j \in \set{1, 2, \ldots, m}; \\ a_i = p_j}} \tup{f\tup{p_j} - f_0} .
\label{pf.lem.app-p.1.1}
\end{align}

Now, forget that we fixed $i$.
We thus have proven \eqref{pf.lem.app-p.1.1}
for each $i \in I$.
Summing the equalities \eqref{pf.lem.app-p.1.1}
over all $i \in I$, we obtain
\begin{align}
\sum_{i \in I} \tup{f\tup{a_i} - f_0}
&= \underbrace{\sum_{i \in I} \sum_{\substack{j \in \set{1, 2, \ldots, m}; \\ a_i = p_j}}}_{= \sum_{j \in \set{1, 2, \ldots, m}} \sum_{\substack{i \in I; \\ a_i = p_j}}}
 \underbrace{\tup{f\tup{p_j} - f_0}}_{\substack{= \sum_{k = 1}^j h_k \\ \tup{\text{by \eqref{pf.lem.app-p.1.sumhk}}}}} \nonumber\\
&= \sum_{j \in \set{1, 2, \ldots, m}} \sum_{\substack{i \in I; \\ a_i = p_j}} \sum_{k = 1}^j h_k
= \underbrace{\sum_{j \in \set{1, 2, \ldots, m}} \sum_{k = 1}^j}_{= \sum_{k = 1}^m \sum_{j = k}^m} \underbrace{\sum_{\substack{i \in I; \\ a_i = p_j}} h_k}_{= \abs{\set{i \in I \mid a_i = p_j}} h_k } \nonumber\\
&= \sum_{k = 1}^m \sum_{j = k}^m \abs{\set{i \in I \mid a_i = p_j}} h_k .
\label{pf.lem.app-p.1.3}
\end{align}

Now, fix $k \in \set{1, 2, \ldots, m}$.
Recall that $P' = \set{p_1 < p_2 < \cdots < p_m}$.
Hence, the elements of $P'$ that are $\geq p_k$
are precisely $p_k, p_{k+1}, \ldots, p_m$.
In other words,
\begin{align}
\set{x \in P' \mid x \geq p_k}
= \set{p_k, p_{k+1}, \ldots, p_m} .
\label{pf.lem.app-p.1.4}
\end{align}
For each $i \in I$, we have the following chain
of logical equivalences:
\begin{align}
\tup{a_i \geq p_k}
&\Longleftrightarrow
\tup{a_i \in \set{x \in P' \mid x \geq p_k}}
\qquad \tup{\text{since $a_i \in P'$ holds always}} \nonumber\\
&\Longleftrightarrow
\tup{a_i \in \set{p_k, p_{k+1}, \ldots, p_m}}
\qquad \tup{\text{by \eqref{pf.lem.app-p.1.4}}} \nonumber\\
&\Longleftrightarrow
\tup{a_i = p_j \text{ for some } j \in \set{k, k+1, \ldots, m}} .
\end{align}
Thus,
\begin{align}
\set{i \in I \mid a_i \geq p_k}
& = \set{i \in I \mid a_i = p_j \text{ for some } j \in \set{k, k+1, \ldots, m}}
\nonumber \\
& = \bigcup_{j = k}^m \set{i \in I \mid a_i = p_j} .
\label{pf.lem.app-p.1.4equiv-equality}
\end{align}

We know that $p_k < p_{k+1} < \cdots < p_m$
(since $p_1 < p_2 < \cdots < p_m$);
thus, the elements $p_k, p_{k+1}, \ldots, p_m$
are distinct.
Hence, for any given $i \in I$, there is at most
one $j \in \set{k, k+1, \ldots, m}$ satisfying $a_i = p_j$.
In other words,
the sets $\set{i \in I \mid a_i = p_j}$
for $j \in \set{k, k+1, \ldots, m}$ are disjoint.
Hence, the size of their union equals the sum of their sizes;
i.e., we have
\[
\abs{\bigcup_{j = k}^m \set{i \in I \mid a_i = p_j}}
= \sum_{j = k}^m \abs{\set{i \in I \mid a_i = p_j}} .
\]
In view of \eqref{pf.lem.app-p.1.4equiv-equality},
this rewrites as
\begin{align}
\abs{\set{i \in I \mid a_i \geq p_k}}
= \sum_{j = k}^m \abs{\set{i \in I \mid a_i = p_j}} .
\label{pf.lem.app-p.1.6}
\end{align}

Now, forget that we fixed $k$.
We thus have proven \eqref{pf.lem.app-p.1.6}
for all $k \in \set{1, 2, \ldots, m}$.

The equality \eqref{pf.lem.app-p.1.3} becomes
\[
\sum_{i \in I} \tup{f\tup{a_i} - f_0}
= \sum_{k = 1}^m
\underbrace{\sum_{j = k}^m \abs{\set{i \in I \mid a_i = p_j}}}_{\substack{= \abs{\set{i \in I \mid a_i \geq p_k}} \\ \tup{\text{by \eqref{pf.lem.app-p.1.6}}}}} h_k
= \sum_{k = 1}^m \abs{\set{i \in I \mid a_i \geq p_k}} h_k .
\]
In view of
\[
\sum_{i \in I} \tup{f\tup{a_i} - f_0}
= \sum_{i \in I} f\tup{a_i} - \abs{I} f_0 ,
\]
this rewrites as
\begin{align}
\sum_{i \in I} f\tup{a_i} - \abs{I} f_0
= \sum_{k = 1}^m \abs{\set{i \in I \mid a_i \geq p_k}} h_k .
\label{pf.lem.app-p.1.at1}
\end{align}
The same argument (applied to the $b_i$ instead of the $a_i$) yields
\begin{align}
\sum_{i \in I} f\tup{b_i} - \abs{I} f_0
= \sum_{k = 1}^m \abs{\set{i \in I \mid b_i \geq p_k}} h_k .
\label{pf.lem.app-p.1.at2}
\end{align}
Hence, \eqref{pf.lem.app-p.1.at1} becomes
\begin{align*}
\sum_{i \in I} f\tup{a_i} - \abs{I} f_0
&= \sum_{k = 1}^m \underbrace{\abs{\set{i \in I \mid a_i \geq p_k}} h_k}_{\substack{\geq \abs{\set{i \in I \mid b_i \geq p_k}} h_k \\ \tup{\text{by \eqref{pf.lem.app-p.1.ass2}}}}} \\
&\geq \sum_{k = 1}^m \abs{\set{i \in I \mid b_i \geq p_k}} h_k
= \sum_{i \in I} f\tup{b_i} - \abs{I} f_0
\end{align*}
(by \eqref{pf.lem.app-p.1.at2}).
Adding $\abs{I} f_0$ to both sides of this inequality, we find
\[
\sum_{i \in I} f\tup{a_i} \geq \sum_{i \in I} f\tup{b_i} .
\]
This proves Lemma~\ref{lem.app-p.1}.
\end{proof}
\end{verlong}

We remark that Lemma~\ref{lem.app-p.1} can be generalized:
Instead of requiring $P$ to be a totally ordered set,
it suffices to assume that $P$ is a poset with the property that
if $a, b, c \in P$ satisfy $a \leq c$ and $b \leq c$,
then $a \leq b$ or $b \leq a$.
The Hasse diagram of such a poset $P$ is a forest if $P$ is finite.

In the following, \defn{$\NN$} shall mean the set $\set{0, 1, 2, \ldots}$.

\begin{corollary} \label{cor.app-p.3}
Let $p$ be a prime.
Let $m \in \NN$ and $s \in \ZZ \setminus \set{-1, -2, \ldots, -m}$.
Let $f : \NN \to \RR$ be a weakly increasing function.
Then,
\[
\sum_{j=1}^m f\tup{v_p\tup{s+j}}
\geq
\sum_{j=1}^m f\tup{v_p\tup{j}} .
\]
\end{corollary}

\begin{proof}[Proof of Corollary~\ref{cor.app-p.3}.]
The integers $s+1,s+2,\ldots,s+m$ are all nonzero
(since $s \in \ZZ \setminus \set{-1, -2, \ldots, -m}$).
Hence, $v_p\tup{s+1}, v_p\tup{s+2}, \ldots, v_p\tup{s+m}$
are well-defined elements of $\NN$.

Let $I$ be the finite set $\set{1,2,\ldots,m}$.

Fix $h \in \NN$.
Let $\mu = \floor{m/p^h}$.
Then, $\mu p^h \leq m$ but $\tup{\mu + 1} p^h > m$.
Hence, the set $I$ has exactly $\mu$ many elements divisible by $p^h$:
namely, $1 p^h, 2 p^h, \ldots, \mu p^h$.
Thus,
\begin{align}
\tup{\text{the number of } i \in I \text{ satisfying } p^h \mid i }
= \mu .
\label{pf.cor.app-p.3.1}
\end{align}

On the other hand, the set $I$ has at least $\mu$ many elements
$i \in I$ satisfying $p^h \mid s+i$
\ \ \ \ \footnote{\textit{Proof.}
The set $\set{s+1, s+2,\ldots,s+m}$ is an interval of $m$ consecutive
integers, and thus contains at least $\mu$ disjoint intervals
consisting of $p^h$ consecutive integers each
(since $m \geq \mu p^h$).
Each of the latter intervals contains a number divisible by $p^h$.
}.
Thus,
\begin{align}
& \tup{\text{the number of } i \in I \text{ satisfying } p^h \mid s+i } \nonumber\\
& \geq \mu = \tup{\text{the number of } i \in I \text{ satisfying } p^h \mid i }
\label{pf.cor.app-p.3.2}
\end{align}
(by \eqref{pf.cor.app-p.3.1}).

Now,
\begin{align*}
\abs{\set{i \in I \mid v_p\tup{s+i} \geq h }}
& = \tup{\text{the number of } i \in I \text{ satisfying } v_p\tup{s+i} \geq h } \\
& = \tup{\text{the number of } i \in I \text{ satisfying } p^h \mid s+i }
\end{align*}
(since an integer $z$ satisfies $v_p\tup{z} \geq h$ if and only if it satisfies $p^h \mid z$).
The same argument (applied to $0$ instead of $s$) yields
\begin{align*}
\abs{\set{i \in I \mid v_p\tup{i} \geq h }}
& = \tup{\text{the number of } i \in I \text{ satisfying } p^h \mid i } .
\end{align*}
In light of these two equalities, we can rewrite \eqref{pf.cor.app-p.3.2} as
\begin{align}
\abs{\set{i \in I \mid v_p\tup{s+i} \geq h }}
\geq
\abs{\set{i \in I \mid v_p\tup{i} \geq h }} .
\label{pf.cor.app-p.3.5}
\end{align}

Now, forget that we fixed $h$.
We thus have proved the inequality \eqref{pf.cor.app-p.3.5} for each $h \in \NN$.
Thus, Lemma~\ref{lem.app-p.1} (applied to $P = \NN$, $a_i = v_p\tup{s+i}$
and $b_i = v_p\tup{i}$) yields
\[
\sum_{i=1}^m f\tup{v_p\tup{s+i}}
\geq
\sum_{i=1}^m f\tup{v_p\tup{i}} .
\]
Renaming the index $i$ as $j$ in this inequality, we obtain
precisely the claim of Corollary~\ref{cor.app-p.3}.
\end{proof}

\begin{proof}[Proof of Example~\ref{exa.p-adic.greedy}.]
Let us only prove the claim for the ultra triple $\tup{E, w, d}$.
(The analogous statement about $\tup{E, w, d'}$ can be proven
in the same way, using a different choice of $f$;
namely, we would have to define $f$ by $f\tup{h} = h$ for all $h$.)
Thus, when we speak of perimeters in the following, we shall
mean perimeters with respect to $\tup{E, w, d}$.

So we need to show that $\tup{1, 2, \ldots, m}$ is a greedy
$m$-permutation of $E$ for the ultra triple $\tup{E, w, d}$.
In other words, we need to show that
\begin{align}
\PER\set{1, 2, \ldots, i} \geqslant \PER\set{1, 2, \ldots, i-1, x}
\label{pf.exa.p-adic.greedy.goal}
\end{align}
for all $i \in \set{1, 2, \ldots, m}$ and $x \in E \setminus \set{1, 2, \ldots, i-1}$.
To prove this, fix $i \in \set{1, 2, \ldots, m}$ and
$x \in E \setminus \set{1, 2, \ldots, i-1}$.
The definition of perimeter yields
\begin{align*}
\PER\set{1, 2, \ldots, i}
&= \sum_{e=1}^i \underbrace{w\tup{e}}_{=0} + \sum_{1 \leq e < f \leq i} d\tup{e, f}
= \sum_{1 \leq e < f \leq i} d\tup{e, f} \\
&= \sum_{1 \leq e < f \leq i - 1} d\tup{e, f} + \sum_{e = 1}^{i-1} d\tup{e, i}
\end{align*}
and similarly
\begin{align*}
\PER\set{1, 2, \ldots, i-1, x}
&= \sum_{1 \leq e < f \leq i - 1} d\tup{e, f} + \sum_{e = 1}^{i-1} d\tup{e, x} .
\end{align*}
In view of these equalities, we see that the inequality
\eqref{pf.exa.p-adic.greedy.goal}
boils down to
\begin{align}
\sum_{e = 1}^{i-1} d\tup{e, i} \geqslant \sum_{e = 1}^{i-1} d\tup{e, x} .
\label{pf.exa.p-adic.greedy.goal2}
\end{align}
Thus, it remains to prove \eqref{pf.exa.p-adic.greedy.goal2}.

Let $f : \NN \to \RR$ be the function that sends each $h \in \NN$ to $- p^{-h}$.
This function $f$ is clearly weakly increasing.
Hence, Corollary~\ref{cor.app-p.3} (applied to $m = i-1$ and
$s = -x$) yields
\begin{align}
\sum_{j=1}^{i-1} f\tup{v_p\tup{-x+j}}
\geq
\sum_{j=1}^{i-1} f\tup{v_p\tup{j}}
\label{pf.exa.p-adic.greedy.3}
\end{align}
(since $x \in E \setminus \set{1, 2, \ldots, i-1}$ leads to
$-x \in \ZZ \setminus \set{-1, -2, \ldots, -\tup{i-1}}$).

For each $e \in \set{1, 2, \ldots, i-1}$, we have $e \neq i$
and thus $i \neq e$
and therefore $d\tup{i, e} = p^{- v_p\tup{i-e}}$ (by the definition of $d$),
so that the definition of $f$ yields
\begin{align}
f\tup{v_p\tup{i-e}} = - \underbrace{p^{- v_p\tup{i-e}}}_{= d\tup{i, e} = d\tup{e, i}} = - d \tup{e, i} .
\label{pf.exa.p-adic.greedy.f1}
\end{align}
Also, for each $e \in \set{1, 2, \ldots, i-1}$, we have $e \neq x$
(since $x \in E \setminus \set{1, 2, \ldots, i-1}$)
and thus $d\tup{e, x} = p^{- v_p\tup{e-x}}$ (by the definition of $d$),
so that the definition of $f$ yields
\begin{align}
f\tup{v_p\tup{e-x}} = - \underbrace{p^{- v_p\tup{e-x}}}_{= d\tup{e, x}} = - d \tup{e, x} .
\label{pf.exa.p-adic.greedy.f2}
\end{align}
Now,
\begin{align*}
- \sum_{e = 1}^{i-1} d\tup{e, i}
&= \sum_{e = 1}^{i-1} \underbrace{\left(- d\tup{e, i}\right)}_{\substack{= f\tup{v_p\tup{i-e}} \\ \tup{\text{by \eqref{pf.exa.p-adic.greedy.f1}}}}}
= \sum_{e = 1}^{i-1} f\tup{v_p\tup{i-e}} = \sum_{j=1}^{i-1} f\tup{v_p\tup{j}} \\
& \qquad \tup{\text{here, we have substituted $j$ for $i-e$ in the sum}} \\
&\leq \sum_{j=1}^{i-1} f\tup{v_p\tup{-x+j}}
\qquad \tup{\text{by \eqref{pf.exa.p-adic.greedy.3}}} \\
&= \sum_{e = 1}^{i-1} f\tup{v_p\tup{\underbrace{-x+e}_{=e-x}}} \\
& \qquad \tup{\text{here, we have renamed the summation index $j$ as $e$}} \\
&= \sum_{e = 1}^{i-1} \underbrace{f\tup{v_p\tup{e-x}}}_{\substack{= - d \tup{e, x} \\ \tup{\text{by \eqref{pf.exa.p-adic.greedy.f2}}}}}
= \sum_{e = 1}^{i-1} \tup{- d\tup{e, x}}
= - \sum_{e = 1}^{i-1} d\tup{e, x} .
\end{align*}
Multiplying both sides of this inequality by $-1$,
we obtain
\[
\sum_{e = 1}^{i-1} d\tup{e, i} \geqslant \sum_{e = 1}^{i-1} d\tup{e, x} .
\]
In other words, \eqref{pf.exa.p-adic.greedy.goal2} is proven.
This completes the proof of Example~\ref{exa.p-adic.greedy}.
\end{proof}

We can actually prove a more general fact:

\begin{proposition} \label{prop.n-adic.greedy}
Let $\NN$, $c$, $\mathbf{r} = \tup{r_0, r_1, r_2, \ldots}$,
$v_{\mathbf{r}}\tup{x}$, $E$ and $d$ be as in
Example~\ref{exa.n-adic}.
Assume that $d\tup{a, b}$ is well-defined for each
$\tup{a, b} \in E \timesu E$.
Let $m \in \NN$.
Assume furthermore that $E$ contains $1, 2, \ldots, m$.
We define $w : E \to \RR$ by setting $w\tup{e} = 0$ for each $e \in E$.

Then, $\tup{1, 2, \ldots, m}$ is a greedy
$m$-permutation of $E$.
\end{proposition}

The proof of this is analogous to the above
proof of Example~\ref{exa.p-adic.greedy}.

\end{document}